\magnification=1200
\hsize 158truemm

\def\ref#1{\lbrack {#1}\rbrack}

\def\ekv#1#2{$${#2}\eqno(#1)$$}
\def\eekv#1#2#3{$$\eqalignno{&{#2}&({#1})\cr &{#3}\cr}$$}
\def\eeekv#1#2#3#4{$$\eqalignno{&{#2}&({#1})\cr &{#3}\cr &{#4}\cr}$$}
\def\eeeekv#1#2#3#4#5{$$\eqalignno{&{#2}&({#1})\cr &{#3}\cr &{#4}\cr
&{#5}\cr}$$}
\def\eeeeekv#1#2#3#4#5#6{$$\eqalignno{&{#2}&({#1})\cr
&{#3}\cr &{#4}\cr &{#5}\cr &{#6}\cr}$$}

\def\iint{\int\hskip -2mm\int}
\def\iiint{\int\hskip -2mm\int\hskip -2mm\int}

\font\liten=cmr10 at 8pt
\font\stor=cmr10 at 12pt

\def\bdd{bounded}
\def\bdy{boundary}
\def\canform{canonical transformation}

\def\hol{holomorphic}
\def\mfld{manifold}
\def\neigh{neighborhood}
\def\op{operator}
\def\og{orthogonal}
\def\pop{pseudodifferential operator}
\def\pseudor{pseudodifferential operator}
\def\schr{Schr{\"o}dinger operator}
\def\sa{self-adjoint}

\def\top{Toeplitz operator}
\def\uf{uniform}
\def\ufly{uniformly}
\def\vf{vector field}
\def\wrt{with respect to }
\def\Op{{\rm Op\,}}
\def\Re{{\rm Re\,}}
\def\Im{{\rm Im\,}}

\centerline{\stor Determinants of \pseudor s}
\centerline{\stor and complex deformations of
phase space}
\medskip
\centerline{\bf A. Melin\footnote{*}{\rm Dept.
of Mathematics, University of Lund, Box 118,
S-22100 Lund, Sweden} and J.
Sj{\"o}strand\footnote{**}{\rm Centre de
Math{\'e}matiques, Ecole Polytechnique,
F-91128 Palaiseau cedex, France and URM
7640, CNRS}}\footnote{}{\it
}
\footnote{}{\it Mathematics subject
classification 2000: \rm 31C10,
35P05, 37J45, 37K05, 47J20, 58J52}\bigskip
\par\noindent \par\noindent \it R{\'e}sum{\'e}.
\liten Consid{\`e}rons un op{\'e}rateur
h-pseudodiff{\'e}rentiel, dont le symbole
$p$ s'{\'e}tend holomorphiquement {\`a} un
voisinage tubulaire
de l'espace de phase r{\'e}el et converge
assez vite vers 1, pour que le d{\'e}terminant
soit bien d{\'e}fini.
Nous montrons que le logarithme du module du
d{\'e}terminant est major{\'e} par $(2\pi
h)^{-n}(I(\Lambda ,p)+o(1))$, $h\to 0$, o{\`u}
$I(\Lambda ,p)$ est l'int{\'e}grale de $\log\vert
p\vert $ sur
$\Lambda$, pour tout $\Lambda $ dans une classe de
d{\'e}formations de l'espace de phase r{\'e}el sur
lesquelles la  restriction de la forme
symplectique est r{\'e}elle et
non-d{\'e}g{\'e}ner{\'e}e. Nous montrons que
$I$ est une fonction Lipschitzienne de
$\Lambda$ et nous {\'e}tudions sa
diff{\'e}rentielle et parfois son hessien.
Sous des hypoth{\`e}ses suppl{\'e}mentaires
faibles, nous montrons qu'un point critique
$\Lambda $ de la fonctionnelle $I$ est de
mani{\`e}re infinit{\'e}simale un minimum {\`a}
l'ordre infini.
\smallskip
\par\noindent \par\noindent \it Abstract.
\liten Consider an h-\pseudor{}, whose symbol
$p$ extends holomorphically to a tubular
\neigh{} of the real phase space and converges
sufficiently fast to 1, so that the
determinant is welldefined. We show that the
logarithm of the modulus of this determinant
is bounded by $(2\pi h)^{-n}(I(\Lambda
,p)+o(1))$, $h\to 0$, where $I(\Lambda ,p)$ is
the integral of $\log\vert p\vert $ over
$\Lambda
$ and $\Lambda $ belongs to a class of
deformations of the real phase space on which
the restriction of the symplectic form is real
and non-degenerate. We show that $I$ is a
Lipschitz function of
$\Lambda $ and we study its differential and
sometimes its Hessian. Under weak additional
assumptions, we show that critical "points"
$\Lambda $ of the functional $I$ are
infinitesimal minima to infinite order.
\medskip
\centerline{\bf 0. Introduction.}
\medskip
\medskip\rm
\par In the theory of non-selfadjoint
operators, determinants play an important
role (see for instance [GoKr],
[MaMa]) and recent developments in the theory
of resonances have brought new interest in
operator determinants (see for instance [Me],
[Zw3], [Vo], for the use of determinants in
getting upper bounds on the density of
resonances, and [GuZw], [Zw1], [Sj4] for the
use in getting traceformulae).
\par There does not seem to be many works
devoted to estimates and asymptotics of
determinants by direct microlocal methods,
and the present work is an attempt in that
direction. An interesting feature is that the
determinant of an operator is independent of
the choice of norm of the Hilbert space where
the operator acts, and sometimes to some
extent even of the space itself. In the study
of resonances the possibility of changing
the ambient Hilbert space has played an
important role since the beginning of the
method of complex scaling and many variants of
that method have been used, including phase
space versions (see [AgCo], [BaCo], [HeSj],
[SjZw]). To get optimal estimates, one often
has to choose a Hilbert space norm adapted to
the problem (see for instance [HeSj], [Sj5],
[Ze]) and this choice is often the result of
phase space analysis.

\par The work [Da] gives examples of
paradoxes encountered when not looking
for the optimal Hilbert space. In that work,
Davies considers a non-selfadjoint \schr{ }
on ${\bf R}$ with a complex valued potential
and shows that one can obtain open sets of
"pseudo spectral" values $z\in{\bf C}$ for
which the norm of the resolvent (when it
exists) is exponentially large, while for this
class of operators, the true spectrum is
discrete and confined to a smaller
$h$-independent set. The explanation of this
(see also [Zw2]) is quite obvious to
specialists in partial differential equations:
If we consider more generally an
$h$-\pseudor{ } $p(x,hD_x)$ on ${\bf R}^n$
with  symbol
$p(x,\xi )$ in a suitable class, and if $\rho
_0=(x_0,\xi _0)\in{\bf R}^{2n}$ is a point
where the Poisson bracket satisfies ${1\over
i}\{ p,\overline{p}\} (\rho _0)>0$, then if
$z_0=p(\rho _0)$, we can construct a
WKB-solution $u(x;h)=a(x;h)e^{i\phi (x)/h}$,
with $\Im
\phi (x)\ge 0$ with equality precisely at
$x_0$ and where $a$ is an amplitude, such
that $\Vert u\Vert _{L^2}=1$, $\Vert
(P-z_0)u\Vert _{L^2}={\cal O}(h^\infty )$
(i.e. ${\cal O}_N(h^N)$ for every $N$). (See
[H{\"o}2] for this result in the high frequency
setting. The modifications for the proof in
the semi-classical framework are obvious.
Also, if $p(x,hD)$ is an analytic
$h$-\pseudor{ } in a suitable class, for
instance a differential operator with
analytic coefficients, then ${\cal
O}(h^\infty )$ can be sharpened to ${\cal
O}(e^{-1/(Ch)})$ for some $C>0$. See [Sj1]
for the additional information needed about
the WKB method with analytic symbols.) While
the closure of the set $p({\bf R}^{2n})$ will
contain all limit points of the
spectrum, when $h\to 0$, there is no reason to
expect equality in general. In the case of an
analytic $h$-\pseudor , say when $p$ extends
to a bounded holomorphic function in a tubular
neighborhood of ${\bf R}^{2n}$, one is tempted
(as in the method of complex scaling and
its variants) to look for some other Hilbert
space norm for which the pseudo-spectrum
becomes smaller. There are ways of doing such
modifications of the space that are associated
to deformations of the real phase space.
(See [Sj1], [HeSj] and the works on complex
scaling already cited). Instead of ${\bf
R}^{2n}$, we are led to consider some other
closed IR-\mfld{ }$\Lambda
\subset {\bf C}^{2n}$, i.e. a manifold
$\Lambda $ of real dimension $2n$
such that the restriction to $\Lambda $ of
the complex symplectic form
\ekv{0.1}
{\sigma =\sum_1^n d\xi _j\wedge dx_j,\ (x,\xi
)\in{\bf C}^{2n}={\bf C}_x^n\times {\bf
C}_\xi ^n,}
is real and non-degenerate (so that $\Lambda
$ is I-Lagrangian, i.e. Lagrangian for $\Im
\sigma $ and R-symplectic, i.e. symplectic
with respect to the restriction of $\Re
\sigma $). If $\Lambda $ is close to ${\bf
R}^{2n}$ or has some suitable
transversality property (as in Theorem 3.6
below), we can find a naturally associated
Hilbert space, on which $p(x,hD_x)$ acts a
bounded operator and it turns out that the
principal symbol is now ${p_\vert}_{\Lambda
}$ rather than ${p_\vert}_{{\bf R}^{2n}}$. The
idea would then be to try to find $\Lambda $
such that $p(\Lambda )$ is as small as
possible. The problem of finding such an
optimal
$\Lambda $ is probably very difficult and very
deep, and it is not excluded that one gets
some more complicated set than a smooth
manifold. Nevertheless, we think that the
problem should be attacked even though a
complete success may be remote or even out of
reach.

\par Recent
uses of Carleman estimates in semi-classical
problems by Lebeau-Robbiano [LeRo] and N. Burq
[Bu] are somewhat related to the ideas
developed here. Carleman estimates are
weighted estimates and when proving such
estimates, one effectively replaces
$p(x,hD)$ by a conjugation $e^{\phi
(x)/h}p(x,hD)e^{-\phi (x)/h}$ with symbol
$p_\phi (x,\xi )=p(x,\xi +i\phi '(x))$ for some
real-valued function $\phi (x)$. Geometrically,
this amounts to replacing the real phase space
by the
IR-\mfld{}
$\Lambda _\phi =\{ (x,\xi +i\phi '(x)); (x,\xi
)\in {\bf R}^{2n}\}$. When proving such
estimates one exploits the negativity of the
Poisson bracket ${1\over i}\{ p_\phi
,\overline{p_\phi }\}$ on the characteristic
set $p_\phi (x,\xi )=0$ or equivalently on
the set $p^{-1}(0)\cap \Lambda _\phi $. (This
method is originally due to H{\"o}rmander in a
different context.) This looks first somewhat
contradictary with what we shall do in the
present paper, namely look for IR-\mfld s on
which the Poisson bracket vanishes where $p$
does, but if we consider how the Carleman
estimates are used for instance in [Bu], we
see that the sharpest results will be obtained
in some kind of a limiting case where the
bracket above would be close to 0.

\par In this paper we do not try to study the
spectral problem, but we apply the same
philosophy to determinants of $h$-\pseudor s
$p(x,hD)$ for which the symbol tends to 1
sufficiently fast when $(x,\xi )\to \infty $,
so that the operator is of the form
$I+q(x,hD)$, with $q(x,hD)$ of trace class
(and so that $\det p(x,hD)$ is well-defined
([GoKr])). It should be noticed that in
practice the condition that $p\to 1$
suffiently fast at infinity can be replaced
by an ellipticity assumption near infinity.
Then one chooses some reference operator
 $\widetilde{p}(x,hD)$ which is elliptic
everywhere (and hence invertible for
sufficiently small $h$) and with the property
that $p(x,\xi )/\widetilde{p}(x,\xi )\to 1$
sufficiently fast when $(x,\xi )\to \infty $,
and the discussion will apply to the relative
determinant $\det
(p(x,hD)\widetilde{p}(x,hD)^{-1})$.

\par In section 1 and 3 we develop the
geometrical and analytical framework for
our results. We define a class of IR-\mfld s
contained in a tubular neighborhood of ${\bf
R}^{2n}$ (or of some other fixed linear
IR-\mfld .) If $\Lambda $ is such a \mfld ,
we define a corresponding $h$-dependent
Hilbert space $H(\Lambda )$ where for
simplicity we have imposed conditions at
infinity that imply that the space coincides
with $L^2({\bf R}^n)$, while the norm will
in general, be equivalent to the standard
$L^2$ norm only up to exponentially large
factors. If $p$ is a bounded holomorphic
function in a sufficiently large tube,
containing ${\bf R}^{2n}$ and $\Lambda $, we
show that $p(x,hD)$ (in Weyl quantization)
acts as a uniformly bounded operator in
$H(\Lambda )$. If $p(x,\xi )-1={\cal
O}(\langle (x,\xi )\rangle ^m)$, $m<-2n$, then
$p(x,hD)-1$ is of trace class and $\det
p(x,hD)$ is well-defined and independent of
whether we consider $p(x,hD)$ as an operator
in $L^2({\bf R}^n)$ or in $H(\Lambda )$. In
sections 2, 3 we show that
\ekv{0.2}
{
\log\vert \det p(x,hD)\vert \le (2\pi
h)^{-n}(I(\Lambda ,p)+o(1)),\ h\to 0, }
where
\ekv{0.3}
{
I(\Lambda ,p):={1\over 2}\int_\Lambda \log
(p(\rho )\overline{p(\rho )})\,\mu (d\rho ). }
Here $\mu (d\rho )=\vert \sigma ^n\vert
/(n!)$ is the symplectic volume element on
$\Lambda $. When $p$ is elliptic on $\Lambda
$ in the sense that $p(\rho )\ne 0$, $\rho
\in \Lambda $, we have equality in (0.2) and
in that case it is easy to see directly that
$I(\Lambda ,p)$ remains unchanged under small
deformations of $\Lambda $ (as long as the
global ellipticity is conserved). The major
problem is then to choose $\Lambda $ so that
$I(\Lambda ,p)$ becomes as small as possible
and this is the subject of the remainder of
the paper. As already pointed out we do not
solve this major problem but obtain several
results of some independent interest. The
discussion also applies to the
similar frameworks developed in [HeSj]
(see [Ze] for corresponding trace class
considerations) and [Sj2].

\par As reviewed in section 1, to smooth
deformations $I\ni t\mapsto \Lambda _t$ of
IR-\mfld s (where $I$ is some bounded open
interval containing $0$), we have an associated
smooth family of functions $I\ni t\mapsto
f_t\in C^\infty (\Lambda _t;{\bf R})$
well-defined up to a
$t$-dependent constant, such that if we let
$f_t$ also denote an almost holomorphic
extension to a neighborhood of $\Lambda _t$
and define a flow $\Phi _t$, by $\partial
_t\Phi _t(\rho )=\widehat{iH_{f_t}}(\rho )$
with
$H_{f_t}$ denoting the complex Hamilton
field, defined with the help of the complex
linear part of the differential of $f_t$ and
the complex symplectic form, and where
$\widehat{iH_{f_t}}=iH_{f_t}+\overline{iH_{f_t}}$
denotes the corresponding real \vf , then
$\Lambda _t=\Phi _t(\Lambda _0)$. The flow
$\Phi _t$ is not a complex \canform , because
of lack of holomorphy, but the restriction
$\kappa _t:\Lambda _0\to \Lambda _t$ becomes
a real
\canform . In the case of
deformations of the class of \mfld s used
in this paper, the corresponding $f_t$
belong to the space $C_b^\infty (\Lambda
_t;{\bf R})$ of smooth functions on $\Lambda
_t$ which are bounded at infinity together
with all their derivatives.
\par In section 4, we show that $I(\Lambda
,p)$ is a Lipschitz function with respect to
such deformations (where the
deformation parameter may also be
multi-dimensional, $t\in {\rm neigh\,}(0,{\bf
R}^k)$ (=some neighborhood of $0$ in ${\bf
R}^k$)), and that the corresponding
derivative (defined for almost all $t$) is
given by
\ekv{0.4}
{
{\rm "}\partial _tI(\Lambda
_t,p){\rm "}=-\int_{\Lambda _t}\langle
[d\,{\rm arg\,}p]_{\Lambda
_t},H_{f_t}\rangle \mu (d\rho )=\int_{\Lambda
_t}\langle [H_{{\rm arg\,}p}]_{\Lambda
_t},df_t\rangle \mu (d\rho ). }
Here a step in the proof is to see that the
differential form $d\,{\rm
arg\,}p_{\Lambda_t} $ and the Hamilton field
$H_{{\rm arg\,}p_{\Lambda_t}}$, where
$p_{\Lambda_t} ={p_\vert}_{\Lambda_t }$,
have $L^1$ coefficients (in a suitable
sense). By $[d\,{\rm arg\,}p]_{\Lambda
_t}$ and $[H_{{\rm arg\,}p}]_{\Lambda _t}$
we denote the $L^1$ extensions from $\Lambda
_t\setminus p_{\Lambda _t}^{-1}(0)$ to all of
$\Lambda _t$.  The last expression in (0.4)
can also be written as
\ekv{0.5}
{
-\int_{\Lambda _t}{\rm div\,}[H_{{\rm
arg\,}p}]_{\Lambda _t}(\rho )\, f_t(\rho )\,\mu
(d\rho ), }
where
\ekv{0.6}
{
{\rm div\,}[H_{{\rm arg\,}p}]_{\Lambda _t}
}
is a distribution of order $\le 1$ with
support in $p_{\Lambda _t}^{-1}(0)$ and
with integral 0. (The divergence of a
Hamilton field is zero.) The problem of
minimizing
$I(\Lambda ,p)$ can be thought of as a
variational problem, and we are then
interested in critical points, i.e we would
like to find an IR-\mfld , such that the
distribution (0.6) is zero.

\par We do not solve the variational
problem in this paper, but in section 7, we
show under some additional fairly weak
assumptions on $\Lambda =\Lambda _0$, that if
$I\ni t\mapsto\Lambda _t$ is a smooth
deformation such that $\Lambda _0$ is
critical (in the sense that (0.6)
vanishes), then for every $N\in{\bf N}$,
there is a constant $C_N$ such that
\ekv{0.7}
{
I(\Lambda _t,p)\ge I(\Lambda _0,p)-C_N \vert
t\vert ^N. }
So,  the
critical points (with some weak additional
assumptions) are infinitesimal minima.
This result depends on two observations.
The first one is that if $f$ is
independent of $t$ and extends to a
bounded holomorphic function in a tube,
then we can define the IR-\mfld s $\Lambda
_t=\exp (\widehat{tiH_f})(\Lambda )$ for
$t\in{\rm neigh\,}(0,{\bf C})$, and it is
quite easy to see that $I(\Lambda _t,p)$
becomes a subharmonic function of $t$. But if
$f_\Lambda $ is real, then we see that
$\Lambda _t$ only depends on $\Re t$, so
$I(\Lambda _t,p)$ becomes a convex function
of $t$. Hence, if $\Lambda _0$ is critical,
$I(\Lambda _0,p)$ has to be a minimum for the
particular family $I(\Lambda _t,p)$.
The above observation can be partially
extended to the case when $f$ is only
$C^\infty $ on $\Lambda _0$, and we use it in
section 5 to compute second derivatives of
$I(\Lambda ,p)$ with respect to the deformation
parameter, for more general deformations. The
second observation is that for a general
deformation
$I\ni t\mapsto \Lambda _t$ of IR-\mfld s, we
can approximate $\Lambda _t$ to infinite
order when $t\to 0$ by the result of an
autonomous deformation acting on $\Lambda
_0$. (The deformation will depend on $t$ but
will be approximately autonomous.) Another
more technical ingredient in the proof is the
approximation of $I(\Lambda ,p)$ by
\ekv{0.8}
{
I_\epsilon (\Lambda ,p)={1\over 2} \int \log
{p(\rho )\overline{p(\rho )}+\epsilon ^2\over
1+\epsilon ^2}\, \mu (d\rho ), }
when $\epsilon \searrow 0$. This
approximation is used at many places, and
some more refined considerations are
developed in section 6.

\par In section 8, we consider the case when
the differentials of the real and imaginary
parts of $p_\Lambda $ are linearly
independent at every point of $p_\Lambda
^{-1}(0)$. In this case we see that
$I(\widetilde{\Lambda },p)$ becomes a smooth
function of
$\widetilde{\Lambda }$ in a neighborhood of
$\Lambda $ and the differential is now a
Radon measure acting on $f_t$:
\ekv{0.9}
{
\partial _tI(\Lambda
_t,p)=2\pi \int_{\Lambda_t \cap
p^{-1}(0)}f_t {i\over 2}\{
p_{\Lambda _t},\overline{p_{\Lambda _t}}\}
\lambda _t(d\rho ) =2\pi \int_{\Lambda_t \cap
p^{-1}(0)}f_t{1\over (n-1)!}\sigma ^{n-1}.}
Here $\lambda _t(d\rho )$ denotes the
Liouville measure on $p_{\Lambda
_t}^{-1}(0)$. From this result we see that
$\Lambda $ is critical iff ${i\over 2}\{
p_\Lambda ,\overline{p_\Lambda }\}\equiv 0$
on $p_\Lambda ^{-1}(0)$. Notice that this
property implies that if we let $p(x,hD)^*$
denote the adjoint with respect to the
scalar product of $H(\Lambda )$, then
there exist $h$-\pseudor s $A,B,C$ of order
$0$ in $h$, such that for the commutator of
$p$ and $p^*$:
\ekv{0.10}
{{1\over h}
[p(x,hD),p(x,hD)^*]=Ap(x,hD)+Bp(x,hD)^*+hC.}
In other words, in order to minimize
$I(\Lambda ,p)$ (and the corresponding
determinant in view of (0.2)), we should
choose $\Lambda $ so that $p(x,hD):H(\Lambda
)\to H(\Lambda )$ tends to be a
normal operator.

\par In section 8, we also consider the
case when $p$ depends holomorphically on a
spectral parameter, and investigate the
minimization problem on an infinitesimal
level. (In a future paper we return to this
situation in the case of dimension 2, and get
much more complete results, including the
determination of all eigenvalues of the
operator in some fixed $h$-independent domain
in the complex plane.)

\par In section 9, we make some further
considerations in the case when $p_\Lambda $
is of principal type ($dp_\Lambda \ne 0$ on
$p_\Lambda ^{-1}(0)$). On one hand we see
that there are situations when $\partial_t
I(\Lambda _t,p)$ is not a Radon measure in
$f_t$, and in the case when $p_{\Lambda_0} $
is real, we see that $I(\Lambda _t,p)$ is in
general not differentiable at $t=0$, because
of a jump in the derivative at that point.
We compute the amount of the jump.

\par In section 10, we consider two
examples of bounds for relative determinants
with a spectral parameter
\ekv{0.11}
{\det((p(x,hD)-z)(\widetilde{p}(x,hD)-z)^{-1}),}
i.e. we study $I(\Lambda
,(p-z)(\widetilde{p}-z)^{-1})$. If $z$ varies
in a region where $\widetilde{p}(x,hD)$ has
no spectrum then $\partial _z\partial
_{\overline{z}}$ applied to the determinant
(0.11) is a constant times the positive
measure obtained by putting a Dirac measure
at each eigenvalue of $p(x,hD)$. It is
therefore of interest to consider $\partial
_z\partial _{\overline{z}}I(\Lambda
,(p-z)(\widetilde{p}-z)^{-1})$ which is
easily seen to be $\ge 0$. (In this section
we let $\Lambda $ be constant, while at the
end of section 7, we give an argument that
indicates that we still have positivity if
we let $\Lambda =\Lambda _z$ depend on $z$
in such way that $\Lambda _z$ becomes a
minimizer.) The first example is when
$p_\Lambda $ is real, we see here that
$\partial _z\partial _{\overline{z}}I$
behaves in good agreement with what can be
obtained for (0.11)  from the Weyl asymptotics
of the real eigenvalues. The second example
concerns the case when $p_\Lambda $ vanishes
to precisely the second order at a point. This
is related to resonances for a
semi-classical \schr , generated by a
non-degenerate critical point of the the
potential. Here we also get good agreement
with the known exactly computable case of
harmonic oscillators.

\par Many estimates in the case when $\Lambda $
is of limited regularity, did not get room
in the present paper and we plan to collect
these results in a separate work  (different
from the one mentioned above).

The contents of the paper is:

 \par 1. IR
\mfld s and their deformations.
\par 2. Determinants of $h$-\pseudor s on
${\bf R}^n$.
\par 3. $h$-\pseudor s and IR
deformations of ${\bf R}^{2n}$.
\par 4. H{\"o}lder properties of
$I(\Lambda ,p)$, $I_\epsilon (\Lambda ,p)$,
the differential w.r.t. $\Lambda $.
\par 5. Second derivative under
non-autonomous flows.
\par 6. Continuity and
convergence for the differential of $\Lambda
\mapsto I(\Lambda ,p)$.
\par 7. Minimality to infinite order
of critical points.
\par 8. The codimension 2 case.
\par 9. The case when $p_{\Lambda
_0}$ is of principal type.
\par 10. Examples.
\bigskip

\def\pseudor{pseudodifferential operator }

\def\Op{{\rm Op\,}}

\centerline{\bf 1. IR
manifolds and their deformations.}
\medskip
\par In this section we recall some
more or less well-known facts. Consider
first locally in ${\bf R}^{2n}$
 a smooth family $\Lambda _t$, $t\in
I\subset{\bf R}$ of smooth Lagrangian
manifolds, so that $\Lambda _t$ is
Lagrangian for every $t$ in the interval
$I$ and of the form $\kappa _t({\bf
R}^n)$, where $\kappa _t(y)$ is smooth in
$(t,y)$ and has injective $y$-differential.
Let
$\nu _t(\kappa _t(y))=\partial _t(\kappa
_t)(y)$ be the corresponding deformation
field, first defined as a section of
$T{({\bf R}^{2n})_\vert }_{\Lambda _t}$. If
$\widetilde{\kappa }_t$ is a second family
with the same properties relative to the
family $\Lambda _t$ and $\widetilde{\nu
}_t$ the corresponding vectorfield, then
$\nu _t-\widetilde{\nu }_t\in T(\Lambda
_t)$, so we can define invariantly a
deformation field $\nu _t$ as a section of
the normal bundle $N(\Lambda _t)=({T({\bf
R}^{2n}) _\vert}_{\Lambda _t})/T(\Lambda
_t)$.

\par Let $\sigma = \sum d\xi _j\wedge dx_j$
be the canonical 2 form on ${\bf R}^{2n}$.
We add $t$ as a variable and consider
$$L:=\{ (t,x,\xi );\, (x,\xi )\in\Lambda
_t\}\subset {\bf R}\times {\bf R}^{2n}.$$
Differential forms on ${\bf R}^{2n}$ will
be considered also as differential forms
on ${\bf R}\times {\bf R}^{2n}$ in the
natural way. Let $\omega $ be a smooth one
form on ${\bf R}^{2n}$ with $d\omega
=\sigma $. (We can for instance take
$\omega =\sum \xi _jdx_j$.) Since
${d\omega _\vert}_{\Lambda _t}=0$, we have
$\omega =df_t$ locally on $\Lambda _t$,
where
$f_t$ is smooth and well defined up to a
$t$ dependent constant. Adjusting the
constants we may assume that $f(t,x,\xi
)=f_t(x,\xi )$ is smooth on $L$. On $L$ we
see that the restriction of $\omega $ to
each submanifold $t={\rm const.}$ is equal
to the corresponding restriction of $df$
and hence
\ekv{1.1}{{\omega_\vert}_{L}=df+q_t(x,\xi
)dt,} for some function $q=q_t(x,\xi )$
which is smooth on $L$. If we replace
$\omega $ by another primitive
$\widetilde{\omega }$ of $\sigma $, then
$\widetilde{\omega }=\omega +dw$ for some
smooth function
$w$, and $f$ above is replaced by
$\widetilde{f}$ with
$\widetilde{f}_t=f_t+({w_\vert}_{\Lambda
_t})$. This does not change $q_t$, so we
get a function $q_t(x,\xi )$ on $L$ which
is well defined up to some arbitrary
smooth function $C(t)$. We fix such a
choice of $q_t$ and define
$$\Lambda =\{(t,\tau ,x,\xi );\, (x,\xi
)\in\Lambda _t ,\, \tau =-q_t(x,\xi )\}.$$
Then $${(\tau dt+\omega )_\vert}_{\Lambda
}\simeq{(-q_tdt+\omega )_\vert}_{L}=df,$$
so
$${(d\tau \wedge dt+\sigma)_\vert}_{\Lambda
}  =d{(\tau dt+\omega )_\vert}_{\Lambda
}= 0.$$
Hence $\Lambda $ is Lagrangian. Let
$\widetilde{q}(t,x,\xi )$ be an
extension of
$q_t(x,\xi )$ to ${\bf R}\times {\bf
R}^{2n}$. Then $\tau +\widetilde{q}$
vanishes on $\Lambda $ so its Hamilton
field ${\partial \over \partial
t}+H_{\widetilde{q}_t}$ is tangent to
$\Lambda $. This means that
$H_{\widetilde{q}_t}=\nu _t$ is the
deformation field defined earlier. Notice
that the choice of extension
$\widetilde{q}_t$ of $q_t$ affects the
Hamilton field only by a component which
is tangent to $\Lambda _t$, so with some
slight abuse of notation, we can say that
\ekv{1.2}
{
\nu _t=H_{q_t}.
}
Here it is understood that the RHS is
defined only as a section of the normal
bundle of $\Lambda _t$. We formulate the
main result so far as a lemma:
\medskip
\par\noindent \bf Lemma 1.1. \it Let
$\Lambda _t\subset{\bf R}^n\times {\bf
R}^n$, $t\in{\rm
neigh\,}(0,{\bf R})$, be a smooth family of
Lagrangian \mfld s, in the sense that $L$
above is a smooth sub\mfld{}. Then {\sl
locally} we can find a smooth real-valued
function
$q_t(x,\xi )$, $(x,\xi )\in\Lambda _t$ (on
$L$) such that if $\widetilde{q}_t(x,\xi
)$ is any smooth real extension of $q_t$
to ${\bf R}\times {\bf R}^n\times {\bf
R}^n$, then ${\partial \over \partial
t}+H_{\widetilde{q}_t}$ is tangent to
$\Lambda $ defined earlier. In other
words, $H_{\widetilde{q}_t}$ is a
local representative for the deformation
field of the family $\{ \Lambda _t\}$.
Moreover, $q_t$ is uniquely determined on
each $\Lambda _t$ up to a $t$-dependent
constant.\rm\medskip

\par If $\Lambda _t$ are defined globally
as closed manifolds of ${\bf R}^n\times
{\bf R}^n$ and simply connected, then the
$q_t$ can be defined globally on  $\Lambda
_t$ and are unique up to a $t$-dependent
constant.

\par We shall apply this discussion to
$I$-Lagrangian sub\mfld s of ${\bf
C}^{2n}={\bf C}_x^n\times {\bf C}_\xi ^n$
and we first review some
differential-geometric and symplectic
notions. On this space we have the complex
symplectic 2-form:
\ekv{1.3}
{\sigma =\sum d\xi _j\wedge dx_j,}
which is a non-degenerate closed
(2,0)-form. The corresponding real 2-forms
${\rm Re\,}\sigma ={1\over 2}(\sigma
+\overline{\sigma })$, ${\rm Im\,}\sigma
={1\over 2i}(\sigma -\overline{\sigma })$
are closed and non-degenerate and hence
give rise to real symplectic structures on
${\bf C}^{2n}$. A general vectorfield of
type (1,0) is a vector field of the type
\ekv{1.4}
{\nu =\sum_1^n(a_j(x,\xi ){\partial \over
\partial x_j}+b_j(x,\xi ){\partial \over
\partial \xi _j}),}
where $a_j$ and $b_j$ are complex
functions defined on some open subset of
${\bf C}^{2n}$. To $\nu $
we associate the real vectorfield
\ekv{1.5}
{
\widehat{\nu }:=2{\rm Re\,}\nu =\nu
+\overline{\nu }, }
which is simply the vector
$(a_1,..,a_n,b_1,..,b_n)$ when identifying
${\bf C}^n\times {\bf C}^n$ with the
underlying real \mfld{} ${\bf R}^{2n}\times
{\bf R}^{2n}$. Notice that $\widehat{\nu }$
is the unique real vectorfield with the
property that  $\nu (f)=\widehat{\nu }(f)$
at every point where the function $f$ is
differentiable and satisfies
$\overline{\partial }f=0$. (Here we use the
standard notation:
$df=\partial f+\overline{\partial }f$,
where $\partial f$
and $\overline{\partial }f$ denote the
complex linear and antilinear parts of the
differential.) If $a_j$, $b_j$
are sufficiently smooth, then locallyand
for $t$ small enough, we can define $\Phi
_t(\rho )=\exp (t\widehat{\nu })(\rho )$
and we notice that the components $(\Phi
_t(\rho ))_{x_j}$, $(\Phi _t(\rho ))_{\xi
_j}$ solve the system of ODEs
\ekv{1.6}
{
{d\over dt}(\Phi _t(\rho ))_{x_j}=a_j(\Phi
_t(\rho )),\ {d\over dt}(\Phi _t(\rho
))_{\xi _j}=b_j(\Phi _t(\rho )),\ \Phi
_0(\rho )=\rho . }

\par Let $f$ be a complex-valued
$C^1$-function on some open subset of ${\bf
C}^{2n}$. We define the complex
Hamilton field of $f$, to be the unique
complex vectorfield of type (1,0) which
satisfies the pointwise relation
\ekv{1.7}
{
\langle \sigma ,H_f\wedge t\rangle
=-\langle \partial f,t\rangle ,\ t\in
T({\bf C}^{2n})\otimes {\bf C}.  }
We have the usual formula:
\ekv{1.8}
{
H_f=\sum_1^n({\partial f\over \partial \xi
_j}{\partial \over \partial x_j}-{\partial
f\over \partial x_j}{\partial \over
\partial \xi _j}). }
If $g$ is a real-valued $C^1$ function on
an open subset of ${\bf C}^{2n}$, we
define the real Hamilton fields
$H_g^\alpha $, for $\alpha ={\rm
Re\,}\sigma $, ${\rm Im\,}\sigma $, by the
usual relation:
\ekv{1.9}
{
\langle \alpha ,H_g^\alpha \wedge t\rangle
=-\langle dg,t\rangle ,\ t\in T({\bf
C}^{2n}).  }

\par There are some useful relations
between the three types of Hamilton
fields, when we have some additional
information about $f$. The most important
case is the one when $df$ is complex
linear at some given point $\rho _0$:
\ekv{1.10}
{\overline{\partial }f(\rho _0)=0,}
and in the following calculations, we
restrict the attention to that point.
Since $\sigma $
is of type (2,0) and $\overline{H_f}$ is
of type (0,1), we have
\ekv{1.11}
{
\langle \sigma ,\overline{H_f}\wedge
t\rangle =0,\ t\in T_{\rho _0}({\bf C}^{2n}
)\otimes {\bf C}. }
>From (1.7), (1.10), we get
\ekv{1.12}
{
\langle \sigma ,H_f\wedge t\rangle
=-\langle df,t\rangle , \ t\in T_{\rho
_0}({\bf C}^{2n})\otimes {\bf C}. }
Restrict the relations to real tangent
vectors ($t\in T_{\rho _0}({\bf C}^{2n})$)
and take sums and differences:
\ekv{1.13}
{
\langle \sigma ,2{\rm Re\,}H_f\wedge
t\rangle =-\langle df,t\rangle , }
\ekv{1.14}
{
\langle \sigma ,2{\rm Im\,}H_f\wedge
t\rangle =-{1\over i}\langle df,t\rangle . }
Taking the real and imaginary parts of
these relations, we get
\ekv{1.15}
{\langle {\rm Re\,}\sigma ,2{\rm
Re\,}(H_f)\wedge t\rangle =-\langle d\,{\rm
Re\,}f,t\rangle , }
\ekv{1.16}
{\langle {\rm Im\,}\sigma ,2{\rm
Re\,}(H_f)\wedge t\rangle =-\langle d\,{\rm
Im\,}f,t\rangle , }
\ekv{1.17}
{\langle {\rm Re\,}\sigma ,2{\rm
Im\,}(H_f)\wedge t\rangle =-\langle d\,{\rm
Im\,}f,t\rangle , }
\ekv{1.18}
{\langle {\rm Im\,}\sigma ,2{\rm
Im\,}(H_f)\wedge t\rangle =\langle d\,{\rm
Re\,}f,t\rangle , }
so at points where $\overline{\partial
}f=0$:
\eekv{1.19}
{\widehat{H_f}=2{\rm Re\,}H_f= H_{{\rm
Re\,}f}^{{\rm Re\,}\sigma }
=H_{{\rm Im\,}f}^{{\rm Im\,}\sigma }}
{\widehat{H_{if}}=-2{\rm Im\,}H_f=H_{{\rm
Re\,}f}^{{\rm Im\,}\sigma }=-H_{{\rm
Im\,}f}^{{\rm Re\,}\sigma }.}

\par Later, it will be useful to have a
relation also in the case when $df(\rho
_0)$ is real. From (1.7) we get for $t\in
T_{\rho _0}({\bf C}^{2n})$:
$$\langle \sigma ,iH_f\wedge t\rangle
=-\langle i\partial _f,t\rangle ,
$$
$$
\langle \sigma ,\overline{iH_f}\wedge
t\rangle =0,
$$
$$
2\langle \sigma ,{\rm Re\,}(iH_f)\wedge
t\rangle =-\langle i\partial f,t\rangle ,
$$
$$
2\langle \overline{\sigma },{\rm
Re\,}(iH_f)\wedge t\rangle =\langle
i\overline{\partial }f,t\rangle ,
$$
$$
\langle {\rm Im\,}\sigma ,4{\rm
Re\,}(iH_f)\wedge t\rangle =-\langle
df,t\rangle ,
$$
so
\ekv{1.20}
{
H_f^{{\rm Im\,}\sigma }=4{\rm
Re\,}(iH_f)=2\widehat{iH_f} }
at every point where $df$ is real.

\par A smooth sub\mfld{ } $\Lambda \subset
{\bf C}^{2n}$ is called I-Lagrangian if it
is Lagrangian (and hence of real dimension
$2n$) with respect to ${\rm Im\,}\sigma $.
It is called R-symplectic if it is a
symplectic manifold with respect to the
restriction ${{\rm Re\,}\sigma
_\vert}_{\Lambda }$, i.e. if this 2-form
is non-degenerate. We say that $\Lambda $
is an IR-\mfld{ } if $\Lambda $ is both
I-Lagrangian and R-symplectic. It is easy to
see that an IR-\mfld{ } is totally real of
maximal dimension (m.t.r from now on, see
[H{\"o}We]) and for such a manifold $\Lambda $
we know that every $f\in C^\infty (\Lambda
)$ has a smooth extension $\widetilde{f}$ to
some neighborhood of $\Lambda $ such that
$\overline{\partial }\widetilde{f}$
vanishes to infinite order on $\Lambda $.
Moreover this extension is unique modulo
functions which vanish to infinite order on
$\Lambda $. We say that $\widetilde{f}$ is
almost holomorphic (a.h.) at $\Lambda $. We
recall that the notion of almost holomorphic
extensions was introduced by L. H{\"o}rmander
[H{\"o}1].
\par Let $\Lambda _t$ be a smooth family
of IR-manifolds. Let $q_t$ be the
corresponding smooth family as in Lemma
1.1, now with ${\rm Im\,}\sigma $ as the
real symplectic form. Let
$f_t=\widetilde{q}_t+ir_t$ be an
a.h. extension of $q_t$, so that
$\widetilde{q}_t$ is a smooth real
extension of $q_t$. Then (1.19) applies and
we get at the points of $\Lambda _t$ a
corresponding deformation field,
\ekv{1.21}
{\nu_t=\widehat{iH_{f_t}} =
H_{\widetilde{q}_t}^{{\rm Im\,}\sigma
}=-H^{{\rm Re\,}\sigma }_{r_t}.} Since the
differential of
$f_t$ is uniquely determined by that of $q_t$ at
every point of $\Lambda _t$, we see that
(1.21) gives a unique section in
${T({\bf C}^{2n})_\vert}_{\Lambda _t}$ and
not just a normal vector. Integrating $\nu
_t$, we get local diffeomorphisms $\kappa
_{t,s}:\Lambda _s\to\Lambda _t$ and from
the last expression in (1.21), it follows
that $\kappa _{t,s}^*({{\rm Re\,}\sigma
_\vert}_{\Lambda _t})={{\rm Re\,}\sigma
_\vert}_{\Lambda _s}$, so $\kappa _{t,s}$
are real canonical transformations, i.e.
canonical transformations between the real
symplectic \mfld s $\Lambda _s$ and
$\Lambda _t$.

\par Finally we notice that $\nu _t$ in
(1.21) is the realization of the deformation
field as a section in $JT(\Lambda _t)$,
where $J$ denotes the mapping of $T{\bf
C}^n$ into itself, induced by
multiplication by the imaginary unit $i$.
Indeed, $\widehat{i\nu }=J\nu $ for every
(1,0)-type vectorfield $\nu $ and since
$f_t$  is real on
$\Lambda _t$ , we see that
$\widehat{H_{f_t}}$ is the vector field
tangent to $\Lambda _t$, given by
\ekv{1.22}
{\widehat{H_{f_t}} =H^{{\sigma
_\vert}_{\Lambda
_t}}_{{{f_t}_\vert}_{\Lambda _t}}.}
\bigskip

\centerline{\bf 2. Determinants of
\pseudor s on ${\bf R}^n$.}
\medskip
\par If $m\in C^\infty ({\bf
R}^{2n};]0,\infty [)$, we denote by $S(m)$
the space of all $a\in C^\infty ({\bf
R}^{2n})$ such that
\ekv{2.1}
{\partial _x^\alpha \partial _\xi ^\beta
a(x,\xi )={\cal O}_{\alpha ,\beta }(m(x,\xi
)),}
uniformly on ${\bf R}^{2n}$ for all
multiindices $\alpha $ and $\beta $. We
will always assume that $m$ is an order
function in the sense that
$m(x,\xi )\le C_0\langle (x,\xi )-(y,\eta
)\rangle ^{N_0}m(y,\eta )$ for some fixed
positive constants $C_0,N_0$, where we use
the notation
$\langle (x,\xi )\rangle =(1+\vert (x,\xi
)\vert ^2)^{1/2}$. For a symbol
$a=a(x,\xi ;h)$ depending on
$h\in ]0,h_0]$, $h_0>0$, we say that $a\in
S(m)$, if (2.1) holds uniformly w.r.t.
$h$. Most of the time such symbols will be
of the type $a=a_0(x,\xi )+hr(x,\xi ;h)$
with $a_0,r\in S(m)$. We then call $a_0$
the principal part or the principal symbol
of $a$.

\par If $a\in S(m)$, we denote by $\Op (a)$
the corresponding $h$-Weyl quantization:
\ekv{2.2}
{\Op (a)u(x)={1\over (2\pi
h)^n}\iint e^{{i\over h}(x-y)\cdot \theta
}a({x+y\over 2},\theta ;h)u(y)dyd\theta ,}
and recall (see for instance [DiSj]) that
$\Op (a)$ is continuous
${\cal S}({\bf R}^n)\to{\cal S}({\bf
R}^n)$,\break ${\cal S}'({\bf R}^n)\to{\rm
S}'({\bf R}^n)$. If $m=1$, then $\Op (a)$
is ${\cal O}(1):L^2\to L^2$. If $m$ is
integrable, then $\Op (a)$ is of trace
class and the corresponding trace class
norm satisfies $\Vert \Op (a)\Vert _{{\rm
tr}}={\cal O}(1)h^{-n}$. The trace is
given by
\ekv{2.3}
{
{\rm tr\,}\Op (a)={1\over (2\pi
h)^n}\iint a(x,\theta ;h)dxd\theta . }
Let $m\in S(m)$ be integrable and let
$p\in S(1)$ be of the form
\ekv{2.4}
{p=1+a,\ a\in S(m),}
with principal symbol $p_0=1+a_0$, where
$a_0$ is the principal symbol of $a$. Then
$\det (\Op (p))$ is well defined, and if
$I\ni t\mapsto a^t\in S(m)$ is a mapping of
class $C^1$ defined on some interval
$I$, and we put $p^t=1+a^t$, then we have
\ekv{2.5}
{{\partial _t\det\Op (p^t)\over \det \Op
(p^t)}={\rm tr\,}(\Op (p^t)^{-1}\Op
(\partial _tp^t)),}
whenever $\Op (p^t)$ has an $L^2$ bounded
inverse $({\rm Op\,}(p^t))^{-1}$.
See also the remark at the end of this
section for some related observations. If we
assume that
$p$ in (2.4) is elliptic in the sense that
$p_0(x,\xi )\ne 0$ on ${\bf R}^{2n}$, then we can find
a map as above with $I=[0,1]$, such that
$p^t$ is elliptic for all $t$ and with
$p^0=1$, $p^1=p$. In fact, we only have
to notice that $\log p_0$ is globally
well defined in $S(m)$, and put
$p^t=p_0^t+t(p-p_0)$, $p_0^t=e^{t\log p_0}$.
For
$h>0$ small enough, we then know that $\Op
(p^t)$ is invertible with inverse $\Op
(q^t)$, where
$q^t=1+b^t$, $b^t\in S(m)$,
$q_0^t=1/p_0^t$. Combining this with
(2.5), (2.3) and standard
$h$-pseudodifferential calculus, we get
\eekv{2.6}
{
{\partial _t\det\Op (p^t)\over \det\Op
(p^t)}={1\over (2\pi h)^n}(\iint {\partial
_tp_0^t\over p_0^t}dxd\xi +{\cal
O}(h))}
{\hskip 22mm ={1\over (2\pi h)^n}(\iint
\partial _t\log (p_0^t)dxd\xi +{\cal O}(h)). }
Integrating, we get
\ekv{2.7}
{\log\det \Op (p)={1\over (2\pi
h)^n}(\iint\log (p_0)dxd\xi +{\cal O}(h)).}
In this identity, we use the natural
branch of the logarithm obtained by
continuous deformation with $\log (1)=0$.
In the RHS, we have a natural deformation
from $(x,\xi )=\infty $, while a priori,
we need to use  the specific deformation
$\Op (p^t)$ in order to define the LHS. It
is clear however that these two branches
of $\log (p_0)$ coincide.

\par We now drop the ellipticity
assumption, and derive an upper bound on
$$\log\vert \det\Op (p)\vert ={1\over
2}\log\det (\Op (p)^*\Op (p)),$$
where $\Op (p)^*$ denotes the adjoint
operator. If $\epsilon >0$, the preceding
discussion applies with $\Op (p)$
replaced by $\Op (p)^*\Op (p)+\epsilon
^2\Op (q)$, if $q\in S(m)$, $\Op (q)\ge
0$, $q_0>0$ near $p_0^{-1}(0)$. We get
for $h>0$ small enough:
\eekv{2.8}
{
\log\vert \det\Op (p)\vert \le {1\over
2}\log\det (\Op (p)^*\Op (p)+\epsilon
^2\Op (q)) }
{\hskip 3cm ={1\over (2\pi h)^n}(\iint
{1\over 2}\log (\vert p_0\vert ^2+\epsilon
^2q_0)dxd\xi +{\cal O}_\epsilon (h)).}

\par For later convenience, we observe
that if $q_0=1$ near $p_0^{-1}(0)$ or if
$\log \vert p_0\vert $ is integrable, then
\ekv{2.9}
{
\iint \log (\vert p_0\vert ^2+\epsilon
^2q_0)dxd\xi -\iint\log ({\vert p_0\vert
^2+\epsilon ^2\over 1+\epsilon ^2})dxd\xi
\to 0,\ \epsilon \to 0. }
Notice that
\ekv{2.10}
{
\lim_{\epsilon \to 0}\iint\log (\vert
p_0\vert ^2+\epsilon ^2q_0)dxd\xi
=\iint \log \vert p_0(x,\xi )\vert ^2
dxd\xi
\in
[-\infty ,+\infty [. }

\par From (2.8), (2.10), we get
\ekv{2.11}
{\log\vert \det\Op (p)\vert \le {1\over
(2\pi h)^n}(\iint\log \vert p_0(x,\xi
)\vert dxd\xi +o(1)),\ h\to 0,}
provided that the integral to the right is
$\ne -\infty $. (It is easy to see that so
is the case when $p_0$ never vanishes to
infinite order at any point. See Lemma 4.2
below.) Note that in the elliptic case, we
get by taking real parts in (2.7):
$$\log\vert \det\Op (p)\vert ={1\over
(2\pi h)^n}(\iint\log\vert p_0\vert dxd\xi
+{\cal O}(h)).$$
\smallskip
\par\noindent \it Remark. \rm Let ${\cal
T}_1$  denote the space of trace class
operators on ${\cal H}=L^2({\bf R}^n )$.
The mapping
$$
 A \mapsto {\psi}(A)=(\det (I+A))(I+A)^{-1}
$$ which is defined near $0\in {\cal T}_1$
extends to an analytic mapping
$$ {\psi}: {\cal T}_1\to  {\cal L}({\cal H})
$$ where the space to the right is the
space of continuous linear operators on
${\cal H}$. (The analyticity of ${\psi}$ at
$A_0\in {\cal T}_1$ follows easily from the
corresponding property in the
finite-dimensional situation if one
approximates $A_0$ in ${\cal T}_1$ by a
finite rank element.) It is  natural to
consider  ${{\cal L}}({\cal H})$
 as the dual of ${\cal T}_1$ obtained via
the pairing
${\cal T}_1 \times {\cal L} ({\cal H}) \ni
(A, B) \mapsto {\rm tr \,}(AB)$. Since
${\psi}(A)(I+A)=(I+A){\psi}(A)=\det
(I+A)\cdot I$ it follows in particular that
$\det (I+A)$ is analytic in ${{\cal
T}}_1$.  Its differential at $A_0 $ is an
element in ${\cal L}({\cal H})$, and we
claim that it is equal to ${\psi} (A_0)$
(cf. \ (2.5)).  For reasons of analyticity
in $A_0$ it suffices to verify this when
$I+A_0$ is invertible. Replacing $A$ by
$(I+A_0)A$ in our considerations and using
the multiplicative property of the
determinant we see that our assertion
follows from the fact that ${\rm tr\,}(A)$
is the linear part in the Taylor expansion
of $\det (I+A)$ at the origin.

\medskip If $\cal A $ denotes the Banach
space of all functions
$a(x,{\xi})$ in ${\bf R}^{2n}$ such that
$$
\sum _{| \alpha|+|\beta |\le N}\iint
|\partial _x ^{\alpha}\partial _{\xi}
^{\beta}a(x,{\xi})|\, dx\, d{\xi} <\infty,
$$ where $N$ is a sufficiently large
positive integer, then the Weyl
quantization $a \mapsto {\rm Op\,}(a)$ is a
continuous linear mapping from $\cal A$ to
${\cal T}_1$. It follows that $\log
|\det({\rm Op\, }(p))|$, where $p=1+a$, is
a plurisubharmonic function in $\cal A$.
The inequality (2.11) provides us with a
bound from above for that function in terms
of $\iint \log |p(x,{\xi})|\, dx\, d{\xi}$,
which is also a plurisubharmonic function
of $a$.
\bigskip

\centerline{\bf 3 $h$-\pseudor s and IR
deformations of ${\bf R}^{2n}$.}
\medskip
\par Let $\Phi _0(x)$ be a real quadratic
form on ${\bf C}^n$ which is strictly
plurisubharmonic (st.pl.s.h.). Let
$$H_{\Phi _0}={\rm Hol\,}({\bf C}^n)\cap
L^2({\bf C}^n; e^{-2\Phi _0/h}L(dx)),$$
where $L(dx)$ denotes the Lebesgue measure on
${\bf C}^n\simeq {\bf R}^{2n}$ and ${\rm
Hol\,}({\bf C}^n)$ is the space of all
holomorphic functions on ${\bf C}^n$. We
recall (see for instance the beginning of
[Sj2]) that there exists a unitary operator (a
generalized Bargman transform) $W: L^2({\bf
R}^n)\to H_{\Phi _0}$, given by
\ekv{3.1}
{
Wu(x)=Ch^{-3n/4}\int e^{{i\over h}\phi
_0(x,y)}u(y) dy, }
where $\phi _0$ is quadratic and
holomorphic on ${\bf C}_x^n\times {\bf
C}_y^n$ with
\ekv{3.2}
{
\det \partial _x\partial _y\phi _0(x,y)\ne
0,\ {\rm Im\,}\partial _y^2\phi_0 >0, }
and such that the complex canonical
tranformation
\ekv{3.3}
{
\kappa _W:\, (y,-\partial _y\phi
_0(x,y))\mapsto (x,\partial _x\phi _0(x,y))  }
maps ${\bf R}^{2n}$ onto
$$
\Lambda _{\Phi _0}:= \{ (x,{2\over
i}{\partial \Phi _0(x)\over \partial x});\,
x\in{\bf C}^n\}.  $$
If we define ${\bf R}^{2n}\ni (y(x),\eta
(x))=\kappa _W^{-1}(x,{2\over i}\partial
_x\Phi _0(x)) $, then it follows that
\ekv{3.4}
{
-(\partial _y\phi _0)(x,y(x))=\eta (x),\
(\partial _x\phi _0)(x,y(x))={2\over
i}\partial _x\Phi _0(x), }
\ekv{3.5}
{-{\rm Im\,}\phi _0(x,y(x))=\Phi
_0(x),\ \Phi _0(x)+{\rm Im\,}\phi _0(x,y)\sim
+\vert y-y(x)\vert ^2.}

\par We also recall that if $a\in C_b^\infty
(\Lambda _{\Phi _0})$ (meaning that $a\in
C^\infty (\Lambda _{\Phi _0})$ is bounded
together with all its derivatives to all
orders), then
\ekv{3.6}
{
{\rm Op\,}(a)u(x)={1\over (2\pi h)^n}\iint
e^{ {i\over h}(x-y)\cdot \theta }a({x+y\over
2},\theta ) u(y)dyd\theta  }
can be defined as a bounded operator
$H_{\Phi _0}({\bf C}^n)\to H_{ \Phi _0}({\bf
C}^n)$, by choosing the only possible
integration contour: $\theta ={2\over
i}\partial _x\Phi _0({x+y\over 2})$.
More generally, if $m\in{\bf R}$ and $a\in
S(\langle (x,\xi )\rangle ^m)$ in the sense
that $\langle (x,\xi )\rangle ^{-m}\nabla ^ka$
is bounded for every $k\in{\bf N}$, then we
can still define ${\rm Op\,}(a)$ by the
oscillatory integral (3.6) to be a
bounded operator $H_{\Phi _0}(\langle
x\rangle ^k)\to H_{\Phi _0}(\langle x\rangle
^{k-m})$ where $H_{\Phi _0}(\langle x\rangle
^k)={\rm Hol\,}({\bf C}^n)\cap L^2({\bf
C}^n;\langle x\rangle ^{2k}e^{-2\Phi
_0/h}L(dx))$ is equipped with the natural
norm. (We have bounds on the operator norms
which are uniform in $h$. See [Sj2].) The
standard metaplectic invariance for the Weyl
quantization holds in this setting and we get
\ekv{3.7}
{
W^{-1}{\rm Op\,}(a)W={\rm Op\,}(b),\ b=a\circ
\kappa _W. }

\par We also recall that the orthogonal
projection $\Pi _0:L^2_{\Phi _0}\to H_{\Phi
_0}$ is given by
\ekv{3.8}
{
\Pi _0u(x)={C\over h^n}\int e^{{2\over h}\Psi
_0(x,y)}u(y) e^{-{2\over h}\Phi
_0(y)}L(dy), }
where $\Psi _0(x,y)$ is the unique quadratic
form on ${\bf C}^n_x\times {\bf C}_y^n$ which
is holomorphic in $x$ and anti-holomorphic in
$y$ (from now on hol-a-hol) and satisfies
\ekv{3.9}
{
\Psi _0(x,x)=\Phi _0(x).
}
(We recall from [Sj1,2] that this is obtained
by writing the identity operator as an
$h$-\pseudor{ }
$$u\mapsto {\widetilde{C}\over h^n}\iint
e^{{2\over h}(\Psi _0(x,\overline{\theta }
)-\Psi _0(y,\overline{\theta }))}u(y)
dyd\theta ,$$ and choosing the integration
contour $\theta =\overline{y}$.)

\par Let $\widetilde{W}\subset\subset
W\subset\subset {\bf C}^n$ be convex open
\neigh s of $0$. Let $\Phi \in C^\infty ({\bf
C}^n;{\bf R})$ with $\nabla \Phi (x)-\nabla
\Phi _0(x)\in C_b^\infty ({\bf C}^n):=S({\bf
C}^n,1)$, so that $\nabla \Phi -\nabla \Phi
_0$ and its derivatives are all bounded. Define
$$\Lambda _\Phi =\{ (x,{2\over i}{\partial
\Phi \over \partial x}(x)); x\in {\bf C}^n\}$$
and assume that
\ekv{3.10}
{
{2\over i}{\partial \Phi \over \partial
x}(x)-{2\over i}{\partial \Phi _0\over
\partial x}(x)\in\widetilde{W},\ x\in {\bf
C}^n, }
so that $\Lambda _\Phi \subset \Lambda _{\Phi
_0}+\{ 0\} \times \widetilde{W}$.

\par Let $p(x,\xi )$ be holomorphic in
$\Lambda _{\Phi _0}+\{ 0\} \times W$ with
$p(x,\xi )={\cal O}(\langle (x,\xi )\rangle
^m)$ for some $m\in{\bf R}$. We
define ${\rm Op\,}(p)$ as in (3.6) and get a
\ufly{} \bdd{} operator
\ekv{3.11}
{{\rm Op\,}(p):\, H_{\Phi _0}(\langle
x\rangle ^k)\to H_{\Phi _0}(\langle x\rangle
^{k-m}),}
where we recall that
\ekv{3.12}
{
H_{\Phi _0}(\langle x\rangle ^k):={\rm
Hol\,}({\bf C}^n)\cap L^2({\bf C}^n;\langle
x\rangle ^{2k}e^{-\Phi _0(x)/h}L(dx)). }

This only requires that $p_{\Lambda _{\Phi
_0}}:={p_\vert}_{\Lambda _{\Phi _0}}\in
S(\Lambda _{\Phi _0},\langle (x,\xi
\rangle ^m)$. We now exploit the holomorphy
assumtion and make contour deformations. Let
$$\theta _\Phi (x,y)=\int_0^1{2\over
i}{\partial \Phi \over \partial x}((1-t)x+ty)
dt,$$
and notice that
$$\theta _{\Phi _0}(x,y)={2\over i}{\partial
\Phi _0\over \partial x}({x+y\over 2}),$$
so that
\ekv{3.13}
{\theta _{\Phi }(x,y)-{2\over i}{\partial
\Phi _0\over \partial x}({x+y\over
2})\in\widetilde{W},\ x,y\in{\bf C}^n.}
As in [Sj1,2], we can replace
the contour $\theta =\theta _{\Phi
_0}(x,y)$ in (3.6) (with $a=p$) by
\ekv{3.14}
{
\theta =\theta _\Phi (x,y)+{i\over
C}{\overline{x-y}\over \langle x-y\rangle },
} where $C>0$ is large enough. It follows
that
\ekv{3.15}
{
{\rm Op\,}(p)={\cal O}(1):\, H_\Phi (\langle
x\rangle ^k)\to H_\Phi (\langle x\rangle
^{k-m}), }
uniformly in $h$, where $H_\Phi (\langle
x\rangle ^k)$ is defined as in (3.12).

\par Now add the assumption that
\ekv{3.16}
{\partial _{\overline{x}}\partial _x\Phi \ge
{\rm Const.}>0,}
so that $\Lambda _{\Phi }$
is an IR \mfld{}. We shall describe ${\rm
Op\,}(p)$ as a kind of \top{} in the spirit
of [BoSj] and start by studying the
asymptotics of the \og{} projection
\ekv{3.17}
{\Pi _\Phi :L^2(e^{-2\Phi /h})\to H_\Phi (1).}
Let $\psi (x,y)\in C^\infty ({\bf C}^n_x\times
{\bf C}^n_y)$ be almost holomorphic in $x$
and almost anti-holomorphic in $y$ at the
diagonal ${\rm diag\,}({\bf C}^n\times {\bf
C}^n)$, such that
\ekv{3.18}
{\psi (x,x)=\Phi (x),\ \partial
_{\overline{x}}\psi ,\partial _y\psi ={\cal
O}_N(\vert x-y\vert ^N),\forall N,}
\ekv{3.19}
{\nabla ^2\psi \in C_b^\infty ,\ \partial
_{\overline{x}}\psi ,\partial _y\psi \in
C_b^\infty .}
For the last part of (3.19), we write $\Phi
=\Phi _0+f$, where $\nabla f\in C_b^\infty $
and take $\psi =\Psi _0+\widetilde{f}$, where
$\widetilde{f}$ is an almost
holomorphic-anti-holomorphic extension of $f$
(viewed as a function on ${\rm diag\,}({\bf
C}^n\times {\bf C}^n)$) with $\partial
_{\overline{x}}\widetilde{f}, \partial _y
\widetilde{f}\in C_b^\infty
$. It is wellknown and easy to check that
\ekv{3.20}
{\Phi (x)+\Phi (y)-2{\rm Re\,}\psi (x,y)\sim
\vert x-y\vert ^2,}
uniformly for $\vert x-y\vert \le 1/C$, if
$C>0$ is large enough.

\par Put
$$Z_j=h\partial _{\overline{x}_j}+\partial
_{\overline{x}_j}\Phi,\, Z_j^*=-h\partial
_{x_j}+\partial _{x_j}\Phi ,\,
{^t\hskip -2pt Z}_j=-h\partial
_{\overline{x}_j}+\partial
_{\overline{x}_j}\Phi ,\,
{^t\hskip -2pt Z}_j^*=h
\partial_{x_j}+\partial_{x_j}\Phi  . $$
If $F(x,y)=2\psi (x,y)-\Phi (x)-\Phi (y)$,
we get
$$Z_j(x,hD_x)\circ
e^{F(x,y)/h}=e^{F(x,y)/h}\circ (h\partial
_{\overline{x}_j}+2\partial
_{\overline{x}_j}\psi (x,y)),$$
$$Z_j^*(x,hD_x)\circ e^{F/h}=e^{F/h}\circ
(-h\partial _{x_j}+2\partial _{x_j}(\Phi
(x)-\psi )),$$
$${^t\hskip -2pt Z}_j(y,hD_y)\circ
e^{F/h}=e^{F/h}\circ (-h\partial
_{\overline{y}_j}+2\partial
_{\overline{y}_j}(\Phi (y)-\psi )),$$
$${^t\hskip -2pt Z}_j^*(y,hD_y)\circ e^{F/h}
=e^{F/h}\circ (h\partial _{y_j}+2\partial
_{y_j}\psi ).$$
We will use these relations only in a region
$\vert x-y\vert \le 1/C$, where we notice
that the 0th order coefficients of the
operators in the right most factors,
($2\partial _{\overline{x}_j}\psi $,
$2\partial _{x_j}(\Phi (x)-\psi )$ etc) all
belong to $C^\infty _b$.
\smallskip
\par\noindent \it Definition. \rm Let $m(x)$
be an order function on ${\bf
C}^n$ (as after (2.1)). An integral operator
\ekv{3.21}
{
Ku(x)=\int k(x,y;h)u(y) dy
}
will be called $m$-negligible if $k$ is
negligible in the sense that
\ekv{3.22}
{
\vert Z_x^\alpha (Z_x^*)^\beta ({^t\hskip
-2pt Z}_y)^\gamma ({^t\hskip
-2pt Z}_y^*)^\delta k(x,y;h)\vert \le
C_{N,\alpha ,\beta ,\gamma ,\delta }
h^N\langle x-y\rangle ^{-N}m(x),\
(x,y)\in{\bf C}^{2n} } for all $N\in {\bf N}$
and all multiindices $\alpha ,\beta ,\gamma
,\delta \in{\bf N}^n$. Here we use the
subscripts $x$, $y$ to indicate in which
variables the operators $Z,Z^*$ etc should
act.
\smallskip
\par\noindent \it Exemple. \rm Let $a(x,y)$
be smooth with support in a sufficiently
small \neigh{} of ${\rm diag\,}({\bf
C}^{2n})$ of the form
$\vert x-y\vert
\le 1/{\cal O}(1)$ and satisfy $\nabla
_x^\alpha
\nabla _y^\beta a={\cal O}_N(1)m(x)(\vert
x-y\vert ^N+h^N)$ for all $N$ and all
multiindices $\alpha ,\beta $, where $m$ is
an order function. Then
$$k(x,y;h)=e^{{1\over h}(2\psi (x,y)-\Phi
(x)-\Phi (y))}a(x,y)$$
is $m$-negligible. In fact, the preceding
computations show that
$$Z_x^\alpha (Z_x^*)^\beta ({^t\hskip -2pt
Z}_y)^\gamma ({^t\hskip -2pt Z}_y^*)^\delta
(k(x,y;h))=e^{{1\over h}(2\psi (x,y) -\Phi
(x)-\Phi (y))}a_{\alpha ,\beta ,\gamma ,\delta
},$$ where $a_{\alpha ,\beta ,\gamma ,\delta }$
satsifies the same estimates as $a$ above. It
then suffices to use that
$$\vert e^{{1\over h}(2\psi (x,y)-\Phi
(x)-\Phi (y))}\vert \le e^{-{1\over Ch}\vert
x-y\vert ^2}.$$
\medskip
\par\noindent \bf Lemma 3.1. \it Let $K$ be
an integral operator as in (3.21). Then $K$
is m-negligible iff
\ekv{3.23}
{
m^{-1}{\rm ad}_x^\epsilon Z^\alpha
(Z^*)^\beta K Z^\gamma  (Z^*)^\delta ={\cal
O}(h^N): L^2\to L^2, }
for all $\alpha ,\beta ,\gamma ,\delta
,\epsilon \in {\bf N}^n$, $N\in {\bf N}$.
\medskip
\par\noindent \bf Proof. \rm The kernel of
the operator in (3.23) is
\ekv{3.24}
{
m(x)^{-1}(x-y)^\epsilon Z_x^\alpha (Z_x^*)^\beta ({^t\hskip -2pt
Z}_y)^\gamma ({^t\hskip -2pt Z}_y^*)^\delta
(k(x,y;h)),
}
and if $K$ is $m$-negligible, this kernel is
${\cal O}(h^N\langle x-y\rangle ^{-N})$ for
all $N\ge 0$, so (3.23) follows.

\par For the opposite implication, let
$x_0\in{\bf C}^n$ and put $\xi _0=2i\partial
_x\Phi (x_0)$. Then with $f_j=\partial
_{\overline{x}_j}\Phi (x)-\partial
_{\overline{x}_j}\Phi (x_0)$, we get
$$Z_j=e^{-i{\rm Re\,}(\xi _0\cdot (\cdot ))/h}
\circ (h\partial _{\overline{x}_j}+f_j(x))
\circ e^{i{\rm Re\,}(\xi _0\cdot (\cdot
))/h},$$
$$Z_j^*=e^{-i{\rm Re\,}(\xi _0\cdot (\cdot
))/h}
\circ (-h\partial _{x_j}+\overline{f_j(x)})
\circ e^{i{\rm Re\,}(\xi _0\cdot (\cdot
))/h}.$$ Notice that $f_j$ and all its
derivatives
\wrt $x$ are uniformly bounded \wrt $x_0$,
in a domain $\vert x-x_0\vert \le {\cal
O}(1)$.

\par Similarly, put $\eta _0=-2i\partial
_y\Phi (y_0)$. Then we get
$${^t\hskip -2pt Z_j}= e^{-i{\rm Re\,}(\eta
_0\cdot (\cdot ))/h}\circ  (-h\partial
_{\overline{y}_j}+f_j(y))\circ e^{i{\rm
Re\,}(\eta _0\cdot (\cdot ))/h},$$
$${^t\hskip -2pt Z}_j^*= e^{-i{\rm Re\,}(\eta
_0\cdot (\cdot ))/h}\circ  (h\partial
_{y_j}+\overline{f_j}(y))\circ e^{i{\rm
Re\,}(\eta _0\cdot (\cdot ))/h},$$

\par  Now assume (3.23). Then for
$(x_0,y_0)\in {\bf C}^n\times {\bf C}^n$,
$\vert x-x_0\vert ,\vert y-y_0\vert \le 1$,
we reduce $Z$, $Z^*$, .. as above and
conclude by letting first $\vert \alpha \vert
+\vert \beta \vert +\vert \gamma \vert +\vert
\delta \vert \le 4n+1$, that $m(x)^{-1}
e^{i{\rm Re\,}(\xi _0\cdot
x)/h}k(x,y;h)e^{i{\rm Re\,}(\eta _0\cdot
y)/h}={\cal O}(h^N)$, i.e.
$m(x)^{-1}k(x,y;h)={\cal O}(h^N)$. In fact,
we multiply $k(x,y;h)$ by $\chi
(x-x_0)e^{i\Re (\xi _0\cdot x)/h}\chi
(y-y_0)e^{i\Re (\eta _0\cdot y)/h}$,
where $\chi \in C_0^\infty $. Since the $f_j$
are \ufly{} \bdd{}, the corresponding
operator will then satisfy (3.23) with the
$Z_j$ replaced by their leading parts.
Repeating this argument with
$K$ replaced by
${\rm ad}_x^\epsilon Z^\alpha (Z^*)^\beta
KZ^\gamma (Z^*)^\delta $, we see that
$${1\over m(x)}(x-y)^\epsilon Z_x^\alpha
(Z_x^*)^\beta ({^t\hskip -2pt Z}_y)^\gamma
({^t\hskip -2pt Z}_y^*)^\delta k(x,y)={\cal
O}(h^N).$$
(We first get this for $\vert x-x_0\vert
,\vert y-y_0\vert \le 1$, but uniformly \wrt
$x_0,y_0$, which are arbitrary.) (3.22)
follows.\hfill{$\#$}
\medskip

\par As a second example of negligible
operators, we look at the off-diagonal
contribution to ${\rm Op\,}(p)$ above, with
the contour (3.14). Along this contour we have
\ekv{3.25}
{
d\theta dy=\pm \det \Big( {\partial (\theta
_\Phi (x,y)+{1\over C}{\overline{x-y}\over
\langle x-y\rangle })\over \partial
\overline{y}}\Big) d\overline{y}dy=J(x,y)L(dy).
} It is clear that $J\in C_b^\infty $. The
realization of $e^{-\Phi /h}{\rm
Op\,}(p)e^{\Phi /h}$ becomes
\ekv{3.26}
{
e^{-\Phi /h}{\rm Op\,}(p) e^{\Phi /h}u(x)=
{1\over (2\pi h)^n}\int e^{{1\over
h}(F(x,y)-{1\over C}{\vert x-y\vert ^2\over
\langle x-y\rangle })}q(x,y;h)u(y)L(dy), }
where
\ekv{3.27}
{F(x,y)=-\Phi (x)+\Phi (y)+2\int_0^1
{\partial \Phi \over \partial
x}((1-t)x+ty)dt\cdot (x-y),}
\ekv{3.28}
{
q(x,y;h)=J(x,y)p({x+y\over 2},\theta _\Phi
(x,y)+{i\over C}{\overline{x-y}\over \langle
x-y\rangle }). }
We notice that
\ekv{3.29}
{
\nabla _x^\alpha \nabla _y^\beta
q(x,y;h)={\cal O}_{\alpha ,\beta }(\langle
{x+y\over 2}\rangle ^m). }

\par $F$ is purely imaginary, and if we let
$F_0$ denote the corresponding function,
defined with $\Phi _0$ instead of $\Phi $,
then $F_0$ is quadratic and
$$
\nabla ^k(F-F_0)={\cal O}(\langle x-y\rangle )
$$
for every $k$. (Recall that $\nabla (\Phi
-\Phi _0)\in C_b^\infty $.) It follows that
\ekv{3.31}
{
\nabla ^k((\nabla F)(x,y)-(\nabla F)(x,x)),\,
\nabla ^k(\nabla F(x,y)-\nabla F(y,y))\,
={\cal O}(\langle x-y\rangle ). }
We compute for $x=y$:
\ekv{3.32}
{\partial _xF=\partial _x\Phi ,\ \partial
_{\overline{x}}F=-\partial
_{\overline{x}}\Phi ,\ \partial
_yF=-\partial _y\Phi ,\ \partial
_{\overline{y}}F=\partial _{\overline{y}}\Phi
,}
so if $G(x,y)=F(x,y)-{1\over C}{\vert
x-y\vert ^2\over \langle x-y\rangle }$ is the
exponent in (3.26), we get
\eekv{3.33}
{
\partial _xG(x,y)=\partial _x\Phi
(x)+r_1(x,y),\ \partial
_{\overline{x}}G(x,y)=-\partial
_{\overline{x}}\Phi (x)+r_2(x,y), }
{\partial _yG(x,y)=-\partial _y\Phi
(y)+r_3(x,y),\ \partial
_{\overline{y}}G(x,y)=\partial
_{\overline{y}}\Phi (y)+r_4(x,y),}
with
\ekv{3.34}
{\nabla ^kr_j(x,y)={\cal O}(\langle
x-y\rangle ).}
Consequently
\eeeekv{3.35}
{Z_{j,x}(e^{G(x,y)/h}q)=e^{G(x,y)/h}(h\partial
_{\overline{x}_j}+r_{2,j}(x,y))q,}
{Z_{j,x}^*(e^{G(x,y)/h}q)=e^{G(x,y)/h}(-h
\partial _{x_j}-r_{1,j}(x,y))q,}
{{^t\hskip -2pt
Z}_{j,y}(e^{G(x,y)/h}q)=e^{G(x,y)/h}(-h\partial
_{\overline{y}_j}-r_{4,j}(x,y))q, }
{{^t\hskip -2pt Z}_{j,y}^*(e^{G(x,y)/h}q)=e^{G
(x,y)/h}(h\partial _{y_j}+r_{3,j(x,y)})q.}
Now (3.29), (3.34) give
$$Z_x^\alpha (Z_x^*)^\beta ({^t\hskip -2pt
Z}_y)^\gamma ({^t\hskip -2pt Z}_y^*)^\delta
(e^{{1\over h}G(x,y)}q)=e^{{1\over h}G}{\cal
O}(1)\langle x-y\rangle ^{\vert \alpha \vert
+\vert \beta \vert +\vert \gamma \vert +\vert
\delta \vert }\langle {x+y\over 2}\rangle
^m.$$
Using finally that $\vert e^{{1\over
h}G(x,y)}\vert =e^{-{1\over Ch}{\vert
x-y\vert ^2\over \langle x-y\rangle }}$,
we conclude that
\ekv{3.36}
{u\mapsto {1\over (2\pi h)^n}\int e^{{1\over
h}(F(x,y)-{\vert x-y\vert ^2\over C\langle
x-y\rangle })}(1-\chi (x-y))q(x,y;h)u(y)L(dy)}
is $\langle \cdot \rangle ^m$-negligible, if
$\chi \in C_0^\infty ({\bf C}^n)$ is equal to
$1$ near 0.
\smallskip
\par\noindent \it Formal construction of the
\og{} projection $\Pi _\Phi $. \rm Let $\psi
$ be as in (3.18--20). Let \break $a(x,y;h)\in
S({\bf C}^{2n},1)$ have its support in $\vert
x-y\vert \le 1/C$, with $C$ large enough,
and such that
\ekv{3.37}
{\partial _{\overline{x}}a(x,y;h),\, \partial
_ya(x,y;h) ={\cal O}(\vert x-y\vert ^\infty
+h^\infty ).}
Consider the \op{}
\ekv{3.38}
{Au(x)=( \pi h)^{-n}\int e^{{2\over
h}\psi (x,y)}a(x,y;h)u(y) e^{-2\Phi
(y)/h}L(dy),}  which is ${\cal O}(1):
L^2(e^{-2\Phi /h}L(dx))\to L^2(e^{-2\Phi
/h}L(dx)) $ and satisfies $\overline{\partial
}\circ A={\cal O}(h^\infty ):L^2(e^{-2\Phi
/h})\to L^2(e^{-2\Phi /h})$. Equivalently for
the reduced operators, we have
$$e^{-\Phi /h}Ae^{\Phi /h}={\cal O}(1):
L^2\to L^2,\ Z_j e^{-\Phi /h}Ae^{\Phi
/h}={\cal O}(h^\infty ):L^2\to L^2.$$
The calculations in the first example above
show that
\ekv{3.39}
{Z^\alpha (Z^*)^\beta e^{-\Phi /h}Ae^{\Phi
/h}Z^\gamma (Z^*)^\delta ={\cal O}(h^{\vert
\alpha \vert +\vert \beta \vert +\vert \gamma
\vert +\vert \delta \vert }):L^2\to L^2,}
and applying ${\rm ad}_x^\epsilon $ to this
operator, we gain another power $h^{\epsilon
/2}$ in the estimates (3.39). Notice that the
composition of $e^{\Phi /h}Ae^{\Phi /h}$ with
an $m$-negligible operator is $m$-negligible
in view of Lemma 3.1.

Formally $A$ is a \pop , for if we put
$\widetilde{\psi }(x,y)=\psi
(x,\overline{y})$,
$\widetilde{a}(x,y)=a(x,\overline{y})$, then
$\widetilde{\psi }$, $\widetilde{a}$ are
almost holomorphic on the anti-diagonal :
$y=\overline{x}$, and
$$\eqalign{Au(x)&=({i\over 2\pi h})^n\iint
e^{{2\over h}(\widetilde{\psi
}(x,\overline{y})-\widetilde{\psi
}(y,\overline{y}))}\widetilde{a}(x,
\overline{y};h)u(y)dyd\overline{y}\cr
&=({i\over 2\pi h})^n\iint_{\theta
=\overline{y}}e^{{2\over h}(\widetilde{\psi
}(x,\theta )-\widetilde{\psi }(y,\theta
))}\widetilde{a}(x,\theta ;h)u(y)dyd\theta
.}$$ With this in mind, we compute
\ekv{3.40}
{(e^{-\Phi /h}Ae^{\Phi /h})^2u(x)=(\pi
h)^{-n}\int I(x,y;h)u(y)L(dy),}
\ekv{3.41}
{I(x,y;h)=(\pi h)^{-n}\int e^{(-\Phi (x)+2\psi
(x,z)-2\Phi (z)+2\psi (z,y)-\Phi
(y))/h}a(x,z;h)a(z,y;h)L(dz).}
Because of the support property of $a$, the
support of $I$ is contained in a set where
$\vert x-y\vert $ is small.

\par We have
\ekv{3.42}
{
{\rm Re\,}(-\Phi (x)+2\psi (x,z)-2\Phi
(z)+2\psi (z,y)-\Phi (y))\sim -(\vert
x-z\vert ^2+\vert z-y\vert ^2). }
Further, when $x=y$, the function
\ekv{3.43}
{z\mapsto -\Phi (x)+2\psi (x,z)-2\Phi
(z)+2\psi (z,y)-\Phi (y)}
has the non-degenerate critical point $z=x=y$.
We then know from [MeSj] that $I(x,y;h)$ has
an asymptotic expansion, when $h\to 0$ and to
understand the exponent appearing there, we
should look for the critical point of the
almost \hol{} extension to the
complexification of ${\bf C}^n_z$. For that,
it is convenient to identify ${\bf C}_z^n$
with the anti-diagonal $\{ (z,w)\in{\bf
C}^{2n};\, w=\overline{z}\}$ and write the
function (3.43) as
\ekv{3.44}
{
z\mapsto -\Phi (x)+2\widetilde{\psi
}(x,\overline{z})-2\widetilde{\psi
}(z,\overline{z})+2\widetilde{\psi
}(z,\overline{y})-\Phi (y). }
A natural choice of almost \hol{} extension
is then
\ekv{3.45}
{
(z,w)\mapsto -\Phi (x)+2\widetilde{\psi
}(x,w)-2\widetilde{\psi
}(z,w)+2\widetilde{\psi
}(z,\overline{y})-\Phi (y). }
We get the critical point $z=x$,
$w=\overline{y}$ (mod ${\cal O}(\vert
x-y\vert ^\infty ))$ and the corresponding
critical value
\ekv{3.46}
{
-\Phi (x)+2\widetilde{\psi
}(x,\overline{y})-\Phi (y)=-\Phi (x)+2\psi
(x,y)-\Phi (y). }
The stationary phase method gives
\ekv{3.47}
{I(x,y;h)=e^{{1\over h}(-\Phi (x)+2\psi
(x,y)-\Phi (y))}b(x,y;h)+r(x,y;h),}
where $r={\cal O}(h^\infty )$ and $b\sim
b_0(x,y)+b_1(x,y)h+..$ in $S(1)$, and each
term in this asymptotic expansion is
determined by the behaviour of
$(z,w)\mapsto
\widetilde{a}(x,w;h)\widetilde{a}(z,
\overline{y};h)$ near the critical point
$w=\overline{y}$, $z=x$, so $\partial
_{\overline{x}}b,\partial _yb={\cal O}(\vert
x-y\vert ^\infty +h^\infty )$, and
$b_0(x,y)=f(x,y)a_0(x,y)^2$, where
\ekv{3.48}
{{1\over C}\le f(x,x)\le C.}
We also notice that if $a\sim a_0+ha_1+...$,
then $b_j=2fa_0a_j+q_j$, where $q_j$ depends
on $a_0,...,a_{j-1}$ only. If we choose $\psi
$ such that $\overline{\psi (x,y)}=\psi
(y,x)$, as we may, and assume that
$\overline{a(x,y;h)}=a(y,x;h)$ we achieve
that $\overline{I(x,y;h)}=I(y,x;h)$.

\par From (3.39), (3.45), the symbol
properties of $b$ and (3.47), we deduce that
\ekv{3.49}
{Z_x^\alpha (Z_x^*)^\beta {^t\hskip -2pt
Z}_y^\gamma ({^t\hskip -2pt Z}_y^*)^\delta
r(x,y;h)={\cal O}(h^{N(\alpha ,\beta ,\gamma
,\delta )}).}
Combining this with the fact that $r={\cal
O}(h^\infty )$, we get
\ekv{3.50}
{
Z_x^\alpha (Z_x^*)^\beta {^t\hskip -2pt
Z}_y^\gamma ({^t\hskip -2pt Z}_y^*)^\delta
r(x,y;h)={\cal O}(h^\infty ).
}
In fact, using the reduction in the second part
of the proof of Lemma 3.1, we see that for
$\vert x-x_0\vert , \vert y-y_0\vert \le 1$,
$x_0=y_0$, we have:
\ekv{\widetilde{3.49}}
{
\partial _{\overline{x}}^\alpha \partial
_x^\beta \partial _{\overline{y}}^\gamma
\partial _y^\delta (e^{{i\over h}({\rm
Re\,}(\xi _0\cdot x+\eta _0\cdot y))}r)={\cal
O}(h^{\widetilde{N}(\alpha ,\beta ,\gamma
,\delta )}), }
while (3.50) is equivalent to
\ekv{\widetilde{3.50}}
{
\partial _{\overline{x}}^\alpha \partial
_x^\beta \partial _{\overline{y}}^\gamma
\partial _y^\delta (e^{{i\over h}({\rm
Re\,}(\xi _0\cdot x+\eta _0\cdot y))}r)={\cal
O}(h^\infty ). }
By interpolation inequalities for
derivatives, ($\widetilde{3.49}$) and the
fact that $r={\cal O}(h^\infty )$, we get
($\widetilde{3.50}$) and hence $(3.50)$ as
claimed.

\par The result (3.50) can be reformulated,
as saying that
\eekv{3.51}
{
(e^{-\Phi /h}Ae^{\Phi /h})^2u(x)=}
{(\pi h)^{-n}\int
e^{{1\over h}(-\Phi (x)+2\psi (x,y)-\Phi
(y))}b(x,y;h)u(y)L(dy)+e^{-\Phi /h}Re^{\Phi
/h}u(x), }
where $e^{-\Phi /h}Re^{\Phi /h}$ is
1-negligible and
\ekv{3.52}
{Ru(x)=(\pi h)^{-n}\int r(x,y;h)u(y)L(dy),\
{\rm supp\,}(r)\subset \{(x,y);\, \vert
x-y\vert
\le {1\over {\cal O}(1)}\} .}

\par We can now construct $a$ as above with
$a(x,x;h)$ real, $1/C\le a_0(x,x)\le C$, such
that $b=a+{\cal O}(\vert x-y\vert ^\infty
+h^\infty )$. Let $\widetilde{\Pi }_\Phi $ be
the correponding \op{} as in (3.38). Then
\ekv{3.53}
{\widetilde{\Pi }_\Phi ^2=\widetilde{\Pi
}_\Phi +R,}
for a new \op{} $R$ with $e^{-\Phi
/h}Re^{\Phi /h}$ 1-negligible. Further, we
may arrange so that $\psi
(y,x)=\overline{\psi (x,y)}$,
$a(y,x)=\overline{a(x,y)}$,
\ekv{3.54}
{
\widetilde{\Pi }_\Phi ^*=\widetilde{\Pi
}_\Phi , }
where the star indicates that we take the
adjoint in $L^2(e^{-2\Phi /h}L(dx))$.

\par We pause in order to recall some general
estimates for $\Pi _\Phi $ by means of the
$\overline{\partial }$-\op{}, and we will
follow the appendix of [Sj2] with some
routine elaborations. Put
$X_j=h^{-1/2}Z_j=h^{1/2}\partial
_{\overline{x}_j}+h^{-1/2}\partial
_{\overline{x}_j}\Phi $, $\overline{\partial }
_\Phi =e^{-\Phi /h}h^{1/2}\overline{\partial
}e^{\Phi /h}=\sum X_j\otimes
d_{\overline{x}_j}^\wedge$, where $\omega
^\wedge$ denotes the operation of left
exterior multiplication by the 1-form $\omega
$. Then
$\overline{\partial }_\Phi $ is a complex and we
put $$\Delta _\Phi =\overline{\partial }_\Phi^*
\overline{\partial }_\Phi +\overline{\partial
}_\Phi \overline{\partial }_\Phi ^*=(\sum
X_j^*X_j)\otimes
1+\sum_{j,k}[X_j,X_k^*]d\overline{x}_j^\wedge
dx_k^\rfloor ,$$
where $[X_j,X_k^*]=2\partial
_{\overline{x}_j}\partial _{x_k}\Phi \in
C_b^\infty $ and $\omega ^\rfloor$ denotes the
transpose of left exterior multiplication
by $\omega ^\wedge$. Here we use the scalar
product
$\langle x,y\rangle ={1\over 2}\Re x\cdot
\overline{y}$ on ${\bf C}^n$ extended to the
complexified space, so that $\langle
dx_j,dx_k\rangle =\langle
d\overline{x}_j,d\overline{x}_k\rangle =0$,
$\langle dx_j,d\overline{x}_k\rangle =\delta
_{j,k}$.
$\Delta _\Phi
$ conserves the degree of differential forms
and if we identify (0,1)-forms with functions
with values in ${\bf C}^n$, the restriction to
such forms is given by
$$\Delta _\Phi ^{(1)}=(\sum X_j^*X_j)\otimes
1+2\Phi ''_{\overline{x}x}.$$
Let $H^0=L^2$ (here in the vector-valued
version $L^2({\bf C}^n;{\bf C}^n)$) and let
$H^1\subset L^2$ be the Hilbert space with
norm
$$\Vert u\Vert _{H^1}^2=\sum \Vert X_ju\Vert
^2+\sum \Vert X_j^*u\Vert ^2+\Vert u\Vert
^2.$$ (More precisely, we define $H^1$ as
the closure of the Schwartz space ${\cal S}$
for the norm above and check that this space
coincides with $\{ u\in L^2;\, X_ju,
X_j^*u\in L^2\}$.)
By construction, the map
$$H^1\ni u\mapsto ((X_ju)_j,(X_j^*u)_j,u)\in
(L^2)^{2n+1}$$
is an isometry. Let $H^{-1}\subset{\cal
S}'$ be the dual of $H^1$. Then the adjoint
map
$$(L^2)^{2n+1}\ni
((u_j)_j,(v_j)_j,u_0)\mapsto v=u_0+\sum
X_ju_j+\sum X_j^*v_j\in H^{-1} $$
is surjective of norm $1$ with a right
inverse of norm 1. In particular,
$$\Vert v\Vert
_{H^{-1}}^2=\inf_{v=u_0+\sum X_ju_j+\sum
X_j^*v_j}\Vert u_0\Vert ^2+\sum\Vert u_j\Vert
^2+\sum \Vert v_j\Vert ^2.$$
We have $H^1\subset H^0\subset H^{-1}$ with
corresponding inequalities for the norms.

\par An easy calculation shows that
\ekv{3.55}
{
[\Delta _\Phi ^{(1)},X_j]=\sum f_kX_k+g,\
[\Delta _\Phi ^{(1)},X_j^*]=\sum
\widetilde{f}_kX_k^*+\widetilde{g}, }
with $f_k,\widetilde{f}_k,g,\widetilde{g}$
\bdd{} in $C_b^\infty $ when $h\to 0$.

\par In [Sj2] we recalled the estimate
$$\Vert u\Vert _{H^1}\le {\cal O}(1)\Vert
\Delta _\Phi ^{(1)}u\Vert _{H^{-1}},\ u\in{\cal
S},$$
implying that $\Delta _\Phi ^{(1)}:H^1\to
H^{-1}$ is bijective with a uniformly \bdd{}
inverse for small $h$. Applying this to
$X_ju$ we get
$$\Vert X_ju\Vert _{H^1}\le {\cal O}(1)
(\Vert X_j\Delta _\Phi ^{(1)}u\Vert
_{H^{-1}}+\Vert [\Delta _\Phi
^{(1)},X_j]u\Vert_{H^{-1}}).
$$
Using (3.55) and the fact that multiplication
by a $C_b^\infty $ function is bounded
$H^j\to H^j$, $j=-1,0,1$, we get
$$\Vert X_ju\Vert _{H^1}\le {\cal O}(1)(\Vert
\Delta _\Phi ^{(1)}u\Vert  +{\cal O}(1)\Vert
u\Vert _{H^1})\le {\cal O}(1)\Vert \Delta
_\Phi ^{(1)}u\Vert .$$
Similarly,
$$\Vert X_j^*u\Vert _{H^1}\le {\cal
O}(1)\Vert \Delta _\Phi ^{(1)}u\Vert .$$
For $s\in{\bf N}$, let $H^s$ be the Hilbert
space of $L^2$ functions for which the norm
\ekv{3.56}
{\Vert u\Vert _{H^s}^2=\sum_{\vert \alpha
\vert +\vert \beta \vert \le s}\Vert X^\alpha
(X^*)^\beta u\Vert ^2}
is finite. (Again one checks that ${\cal S}$
is a dense subspace.) Let $H^{-s}\subset{\cal
S}'$ be the dual space. Then
$$..\subset H^2\subset H^1\subset H^0\subset
H^{-1}\subset H^{-2}\subset ..$$ with
corresponding inequalities for the norms. The
last estimates show that
\ekv{3.57}
{\Vert u\Vert _{H^2}\le {\cal O}(1)\Vert
\Delta _\Phi ^{(1)}u\Vert .}
Continuing this argument, we get
\ekv{3.58}
{\Vert u\Vert _{H^{s+2}}\le{\cal O}(1)\Vert
\Delta _\Phi ^{(1)}u\Vert _{H^s},}
for $-1\le s\in{\bf Z}$. By regularization
and bijectivity for $s=1$, we see that
$\Delta _\Phi ^{(1)}:H^{s+2}\to H^s$ is
bijective with a uniformly \bdd{} inverse for
$s\ge -1$. By duality this extends to all
$s\in{\bf Z}$. (It helps to use that $(\Delta
_\Phi ^{(1)})^{-1}:{\cal S}\to{\cal S}$, as
we recalled in [Sj2].)

\par We have ${\rm ad}_x\Delta _\Phi
^{(1)}={\cal O}(h^{1/2}):H^{s+1}\to H^s$,
${\rm ad}_x^2\Delta _\Phi ^{(1)}:{\cal
O}(h):H^s\to H^s$, ${\rm ad}_x^k\Delta _\Phi
^{(1)}=0$, $k\ge 3$. It
follows that
${\rm ad}_x(\Delta _\Phi
^{(1)})^{-1}=-(\Delta _\Phi ^{(1)})^{-1}{\rm
ad}_x(\Delta _\Phi ^{(1)})(\Delta _\Phi
^{(1)})^{-1}={\cal O}(h^{1/2}):H^s\to H^{s+3},$
and more generally,
\ekv{3.59}
{{\rm ad}_x^k(\Delta _\Phi ^{(1)})^{-1}={\cal
O}(h^{k/2}): H^s\to H^{s+2+k}.}
Recall also (see for instance [Sj2]), that
\ekv{3.60}
{e^{-\Phi /h}\Pi _\Phi e^{\Phi
/h}=1-\overline{\partial }_\Phi ^*(\Delta
_\Phi ^{(1)})^{-1}\overline{\partial }_\Phi .}

\par Since $\overline{\partial }_\Phi e^{-\Phi
/h}\widetilde{\Pi }_\Phi e^{\Phi /h}$ is
1-negligible, we see from (3.59), (3.60)
that
$$
e^{-\Phi /h}\Pi _\Phi e^{\Phi /h}e^{-\Phi
/h}\widetilde{\Pi }_\Phi e^{\Phi
/h}=e^{-\Phi /h}\widetilde{\Pi }_\Phi
e^{\Phi /h}+{\rm negligible}.
$$
By duality,
$$
e^{-\Phi /h}\widetilde{\Pi } _\Phi e^{\Phi
/h}e^{-\Phi /h}\Pi _\Phi e^{\Phi
/h}=e^{-\Phi /h}\widetilde{\Pi }_\Phi
e^{\Phi /h}+{\rm negligible}.
$$
Here "negligible" means some (1-)negligible
operator.
Put
\ekv{3.61}
{\widehat{\Pi }_\Phi =\Pi _\Phi
\widetilde{\Pi }_\Phi \Pi _\Phi ,}
so that
\ekv{3.62}
{
e^{-\Phi /h}\widehat{\Pi }_\Phi e^{\Phi
/h}=e^{-\Phi /h}\widetilde{\Pi }_\Phi e^{\Phi
/h}+{\rm negligible}. }
Then
\eekv{3.63}
{\widehat{\Pi }_\Phi ={\cal
O}(1):L^2(e^{-2\Phi /h}L(dx))\to H_\Phi ,\
\widehat{\Pi }_\Phi ^*=\widehat{\Pi }_\Phi ,}
{\widehat{\Pi }_\Phi ^2=\widehat{\Pi }_\Phi
+e^{\Phi /h}({\rm negligible})e^{-\Phi /h}.}
It follows from (3.63) that the spectrum of
$\widehat{\Pi }_\Phi $ is concentrated to an
${\cal O}(h^\infty )$ neighborhood of $\{
0,1\}$. Let $\check{\Pi }_\Phi $ be the
spectral projection associated to the part of
the spectrum close to 1:
\ekv{3.64}
{
\check{\Pi }_\Phi ={1\over 2\pi
i}\int_\gamma  (z-\widehat{\Pi }_\Phi
)^{-1}dz, }
where $\gamma $ is the positively oriented
\bdy{} of the disc $D(1,1/2)$.
\medskip
\par\noindent \bf Proposition 3.2. \it
$\check{\Pi }_\Phi =\widehat{\Pi }_\Phi
+e^{\Phi /h}({\rm negligible})e^{-\Phi /h}.$
\medskip
\par\noindent \bf Proof. \rm As an
approximation for $(z-\widehat{\Pi }_\Phi
)^{-1}$, we try
$$R(z)={1\over z}(1-\widehat{\Pi }_\Phi
)+{1\over z-1}\widehat{\Pi }_\Phi .$$
Then, dropping the subscript $\Phi $,
\ekv{3.65}
{(z-\widehat{\Pi })R(z)=1-({1\over z}-{1\over
z-1})\widehat{\Pi }(1-\widehat{\Pi }).}
Here $\widehat{\Pi }(1-\widehat{\Pi
})=e^{\Phi /h}({\rm negligible})e^{-\Phi
/h}$, and we next verify that that if
$K$ is of the form $e^{\Phi /h}({\rm
negligible})e^{-\Phi /h}$, then
$(1-K)^{-1}=1+L$, where $L$ has the same
property. Indeed,
$L=K+K^2+K^3+..=K+K(1-K)^{-1}K$ and we see
that $L$ has the required property. It
follows that \ufly{} for $z\in\gamma $:
$$(z-\widehat{\Pi })^{-1}=R(z)(1-({1\over
z}-{1\over z-1})\widehat{\Pi }(1-\widehat{\Pi
}))^{-1}=R(z)+e^{\Phi /h}({\rm
negligible})e^{-\Phi /h}.$$
Plugging this into (3.64), we get the
proposition.\hfill{$\#$}
\medskip

\par Notice that we also have
\ekv{3.66}
{\check{\Pi }_\Phi =\widetilde{\Pi }_\Phi
+e^{\Phi /h}({\rm negl.})e^{-\Phi /h},}
and that by construction, $\check{\Pi }_\Phi
$ is an \og{} projection in
$L^2(e^{-2\Phi /h}L(dx))$ with image
contained in $H_\Phi $, so that $\check{\Pi
}_\Phi \prec \Pi _\Phi $ in the sense
that ${\cal R}(\check{\Pi }_\Phi )\subset
{\cal R}(\Pi _\Phi )$. In particular, $\Pi
_\Phi -\check{\Pi }_\Phi $ is an \og{}
projection.
\medskip
\par\noindent\bf Proposition 3.3. \it
$\check{\Pi }_\Phi =\Pi _\Phi $.
\medskip
\par\noindent \bf Proof. \rm In the case
$\Phi =\Phi _0$, we have $\check{\Pi }_{\Phi
_0}=\Pi _{\Phi _0}+{\cal O}(h^\infty )$,
since our construction will reproduce the
explicitly known $\Pi _{\Phi _0}$. Hence
$\check{\Pi }_{\Phi _0}=\Pi _{\Phi _0}$.
\par Assume first that $\Phi -\Phi _0$ is
\bdd{} so that $\Phi -\Phi _0\in C_b^\infty $.
Put $\Phi _t=t\Phi +(1-t)\Phi _0$. Then
$L^2(e^{-2\Phi _t/h}),H_{\Phi _t}$ are
independent of $t$ as spaces and the norms
vary continously with $t$. It is therefore
clear that $\Pi _{\Phi _t}$ is
norm-continuous in $t$ in the sense that
$\Vert \Pi _{\Phi _s}-\Pi _{\Phi _t}\Vert \to
0$, $s\to t$. We can make the
construction of $\check{\Pi }_{\Phi _t}$ so
that this operator also becomes norm-continous
in $t$. We then have a norm-continous family
of \og{} projections $\Pi _{\Phi
_t}-\check{\Pi }_{\Phi _t}$ which vanishes
for $t=0$, and hence for all $t$. Hence we
get the proposition when $\Phi -\Phi _0$ is
\bdd .

\par In the general case, let $R(t)$ be a
smooth increasing function: $[0,1[\to
[R(0),+\infty [$ with $R(0)\gg 0$, $R(t)\to
\infty $, $t\to 1$, and put
\ekv{3.67}
{\Phi ^t(x)=\Phi _0(x)+(\Phi (x)-\Phi
_0(x))\chi ({x\over R(t)}),\ 0\le t\le 1,}
with the convention that $\Phi ^1=\Phi $.
Here $\chi \in C_0^\infty (B(0,2);[0,1])$,
$\chi (x)=1$ for $x\in B(0,1)$. Then $\nabla
^2\Phi ^t$ belongs to a \bdd{} set in
$C_b^\infty $, $\nabla \Phi ^t-\nabla \Phi
_0$ is \ufly{} \bdd . $\Phi ^t$ is \ufly{}
st.pl.sh. also \wrt{} $t$.
Moreover $\nabla ^k\Phi _t\to \nabla
^k\Phi $, $t\to 1$, \ufly{} on every compact,
and $\Phi ^t-\Phi _0$ is \bdd{} for every
fixed $t<1$. We then know that $\check{\Pi
}_{\Phi _t}=\Pi _{\Phi ^t}$ for $t<1$.

\par On the other hand, the projections
$e^{-\Phi ^t/h}\check{\Pi }_{\Phi ^t}e^{\Phi
^t/h}$, $e^{-\Phi ^t/h}\Pi _{\Phi ^t}
e^{\Phi ^t/h}$ are strongly continuous on
$[0,1]$, in the first case by examining the
construction and restricting to ${\cal S}$,
in the second case by looking at the
representation (3.60), that we also restrict
to ${\cal S}$. It follows that the two
operators coincide also for $t=1$, so we get
the proposition in the general
case.\hfill{$\#$}.
\medskip

\par The proposition and (3.66) imply that
\ekv{3.68}
{\Pi _\Phi =\widetilde{\Pi }_\Phi +e^{\Phi
/h}({\rm negligible})e^{-\Phi /h},}
which completes the promised asymptotic study
of $\Pi _\Phi $.

\par We next prepare for trace computations.
Write (3.38) for $\widetilde{\Pi }_\Phi $:
\ekv{3.69}
{
\widetilde{\Pi }_\Phi u(x)= {1\over (\pi
h)^n}\int e^{{2\over h}(\psi (x,y)-\Phi (y))}
f(x,y;h)u(y)L(dy), }
where $f\sim \sum_0^\infty f_j(x,y)h^j$ in
$S(1)$, ${\rm supp\,}f\subset \{ (x,y);\,
\vert x-y\vert \le 1/C\}$, $\partial
_{(\overline{x},y)}f={\cal O}(\vert x-y\vert
^\infty )$, and put $\widetilde{\psi
}(x,y)=\psi (x,\overline{y})$, $\widetilde{f}
(x,y;h)=f(x,\overline{y};h)$. Recall that
$\widetilde{\Pi }_\Phi $ is formally the
identity \op{} viewed as a \pop{}, and that
$L(dy)=(i/2)^ndyd\overline{y}$:
$$1u(x)=({i\over 2\pi h})^n\iint e^{{2\over
h}(\widetilde{\psi }(x,\theta )-\widetilde{\psi
}(y,\theta ))}\widetilde{f}(x,\theta
;h)u(y)dy d\theta .$$
Use the Kuranishi trick:
$$2(\widetilde{\psi }(x,\theta
)-\widetilde{\psi }(y,\theta ))=i(x-y)\cdot
\widetilde{\theta }(x,y,\theta ),\
\widetilde{\theta }(x,x,\theta )={2\over
i}{\partial \widetilde{\psi }\over \partial
x}(x,\theta ),$$
and
$$1u(x)={1\over (2\pi h)^n}\iint e^{{i\over
h}(x-y)\cdot \widetilde{\theta }
}u(y)dyd\widetilde{\theta },$$
to see that
$$i^n\widetilde{f}_0(x,\theta )dxd\theta
=dxd\widetilde{\theta }(x,x,\theta ).$$
Put $\theta =\overline{x}$:
$$i^n\widetilde{f}_0(x,\overline{x})dx
d\overline{x}={{dxd\xi }_\vert
}_{\Lambda _\Phi },
$$
since $\widetilde{\theta
}(x,x,\overline{x})={2\over i}{\partial \Phi
\over \partial x}(x)$. This can also be
written
\ekv{3.70}
{{1\over \pi ^n}f_0(x,x)L(dx)={1\over (2\pi
h)^n}{dxd\xi _\vert}_{\Lambda _\Phi },}
where we identify ${\bf C}_x^n$ with $\Lambda
_\Phi $ by means of the natural projections
and ${dxd\xi _\vert}_{\Lambda _\Phi }=\mu
(d(x,\xi ))$ is the symplectic volume form on
$\Lambda _\Phi $.

\par Consider the operator ${\rm Op\,}(p)$ in
(3.11), (3.15). Noticing that $\Pi _\Phi $,
$\widetilde{\Pi }_\Phi $ act as \bdd{} \op{}s
on $H_\Phi (\langle x\rangle ^k)$, we write
${\rm Op\,}(p)={\rm Op\,}(p)\Pi _\Phi $,
recall that $\Pi _\Phi =\widetilde{\Pi }_\Phi
+e^{\Phi /h}({\rm negl.})e^{-\Phi /h}$, and
develop ${\rm Op\,}(p)\widetilde{\Pi }_\Phi $
by means of stationary phase (after pealing
off negligible off-diagonal parts), to see
that
\ekv{3.71}
{
{\rm Op\,}(p)={\rm top}(q)+e^{\Phi
/h}(\langle \cdot \rangle
^m\hbox{-negl.})e^{-\Phi /h}, }
\ekv{3.72}
{
{\rm top\,}(q)u(x)={1\over (\pi h)^n}\int
e^{{2\over h}(\psi (x,y)-\Phi
(y))}f(x,y;h)q(x,y;h)u(y)L(dy). }
Here $q(x,y;h)\in S(\langle \cdot \rangle
^m)$ has its support in $\vert x-y\vert \le
1/C$, $q\sim \sum_0^\infty  q_j(x,y)h^j$ in
$S(\langle \cdot \rangle ^m)$, $\partial
_{\overline{x},y)}q_j={\cal O}(\vert x-y\vert
^\infty \langle x\rangle ^m)$, and
$q_0(x,x)=p(x,{2\over i}{\partial \Phi \over
\partial x}(x))$.

\par Let $\widetilde{S}(\langle \cdot \rangle
^m)$ be the space of $q(x,y;h)\in S(\langle
\cdot \rangle ^m)$ with support in $\vert
x-y\vert \le 1/C$, such that $\partial
_{\overline{x},y}q={\cal O}((\vert x-y\vert
^\infty +h^\infty )\langle x\rangle ^m)$. Let
$\widetilde{S}_{{\rm cl}}(\langle \cdot
\rangle ^m)$ be the subspace of $q$ with
$q\sim\sum_j^\infty q_j(x,y)h^j$ in
$S(\langle \cdot \rangle ^m)$, where
$\partial _{\overline{x},y}q_j={\cal O}(\vert
x-y\vert ^\infty \langle x\rangle ^m)$.

\par If $q(x,y)\in C^\infty ({\bf C}^{2n})$,
we define $q^*(x,y)=\overline{q(y,x)}$. The
\sa{}ness of $\widetilde{\Pi }_\Phi $ is
equivalent to $f^*=f$ (assuming that we have
chosen $\psi $ with $\psi ^*=\psi $). We
have
\ekv{3.73}
{
{\rm top\,}(q)^*={\rm top\,}(q^*).
}
Moreover, by stationary phase,
\ekv{3.74}
{
{\rm top\,}(q)^*{\rm top\,}(q)={\rm
top\,}(r)+e^{\Phi /h}(\langle \cdot \rangle
^{2m}\hbox{-negl.})e^{-\Phi /h}, }
where $r\in S(\langle \cdot \rangle ^{2m})$,
$r^*=r$, $r\sim\sum_1^\infty r_jh^j$,
$\partial _{(\overline{x},y)}r_j={\cal
O}(\vert x-y\vert ^\infty \langle x\rangle
^{2m})$, $r_0=q_0^*q_0$. (3.74) is a special
case of a more general result for the
composition of two Toeplitz operators ${\rm
top\,}(a)\circ {\rm top\,}(b)$ that can be
formulated in the obvious way.

\par Recall that if $k\in C_0^\infty ({\bf
C}^{2n})$, ${\rm diam\,}{\rm supp\,}(k)\le
C_0$ and $Ku(x)=\int k(x,y)u(y)L(dy)$, then
the trace class norm of $K$ is $\le
C_1(C_0,n)\Vert \widehat{k}\Vert _{L^1}$,
where $\widehat{k}$ is the Fourier transform
of $k$. Moreover the trace of $K$ is given by
${\rm tr\,}(K)=\int k(x,x)L(dx)$.
\medskip
\par\noindent \bf Lemma 3.4. \it Let $m<-2n$.
If $R$ is $\langle \cdot \rangle
^m$-negligible, then $R$ is of trace class
with trace class norm $\le {\cal O}(h^\infty
)$.\medskip
\par\noindent \bf Proof. \rm Let $r(x,y;h)$
be the kernel of $R$. For $(x_0,y_0)\in{\bf
C}^{2n}$, we define $\xi _0=2i\partial _x\Phi
(x_0)$, $\eta _0=-2i\partial _y\Phi (y_0)$ as
in the proof of Lemma 3.1. Then for $\vert
x-x_0\vert ,\vert y-y_0\vert \le 1$, we have
$$\nabla _x^k\nabla _y^\ell ( e^{i{\rm Re\,}\xi
_0\cdot x/h}r(x,y;h) e^{i{\rm Re\,}\eta
_0\cdot y/h})={\cal O}(1) h^N\langle
x_0-y_0\rangle ^{-N}\langle x_0\rangle ^m,$$
for every $N\i{\bf N}$
and all multiindices $k,\ell$. (The estimates
are \uf{} in $(x_0,y_0)$.) If $\chi \in
C_0^\infty (B(0,1))$, it follows that $\chi
(\cdot -x_0)R\chi (\cdot -y_0)$ is of trace
class norm $\le {\cal O}(1)h^N\langle
x_0-y_0\rangle ^{-N}\langle x_0\rangle ^m$.
Using partitions of unity separately in the
$x$ and $y$ variables, we get the
lemma.
\hfill{$\#$}
\medskip
\par\noindent \bf Lemma 3.5. \it Let $m<-2n$,
and let $q\in S(\langle \cdot \rangle ^m)$,
with support in $\vert x-y\vert \le 1/C$.
Then ${\rm top\,}(q)$ is of trace class as
an \op{} on $L^2(e^{-2\Phi /h}L(dx))$ and
the corresponding trace class norm is ${\cal
O}(h^{-\rho })$, for every $\rho >5n/2$. Further,
$${\rm tr\,}{\rm top\,}(q)={1\over (2\pi
h)^n}(\int_{\Lambda _\Phi }q(x,\xi )\mu
(d(x,\xi ))+{\cal O}(h)\int_{\Lambda _\Phi
}\vert q(x,\xi )\vert \mu (d(x,\xi ))).$$
\medskip

\par\noindent
\bf Proof. \rm Let $\chi$ be a cut-off
function with support in the unit ball and let
$b(x,y)$ be the integral kernel of the
operator $B=\chi (\cdot \, -x_0)e^{-{\Phi}/h}
{\rm top\,}(q)e^{{\Phi}/h}\chi (\cdot \,
-y_0)$. This has the same trace class norm
$\Vert B \Vert _{\rm tr}$ as the operator
$B_1$ with integral kernel $b_1(x,{\xi})
=e^{i\Re
  ({\xi}_0\cdot x)/h}b(x,y)e^{i\Re
({\eta}_0\cdot y)/h}$.
  If $ 2n< {\tau} < 2n+1$, then
  $$\Vert B_1 \Vert_{\rm tr}\le
C\Vert\widehat{b_1}\Vert_{L^1}\le C_{\tau}
\Vert \langle D \rangle ^{\tau} b_1 \Vert
  _{L^2}\le C_\tau (\Vert\langle D
\rangle^{2n}b_1\Vert_{L^2})^{2n+1-{\tau}}(\Vert
\langle D \rangle ^{2n+1}b_1
\Vert _{L^2})
  ^{{\tau}-2n}.$$
  The estimates above together with (3.20)
show that
  $$
   \Vert \langle D \rangle ^jb_1 \Vert _{L^2}
   = {\cal O}(1) h^{-n-j} \langle  x_0 \rangle
^m
   (\iint \chi (x-x_0)\chi
(y-y_0)e^{-c|x-y|^2/h }\,
  L(dx)L(dy))^{1/2}
  $$ where $j=2n,2n+1$ and $c$ is a positive
constant. The right hand-side  may be
estimated from above by a constant times
$h^{-n-j+n/2}\langle x_0 \rangle ^m
s(x_0-y_0)$ where $s$ is a rapidly decreasing
function. A combination of these estimates
with $j=2n$ and $j=2n+1$ gives
  $$
   \Vert B \Vert _{\rm tr} ={\cal O}(1)
h^{-n/2-{\tau}}
   \langle x_0 \rangle ^m s(x_0-y_0).
   $$ The assertion about the trace class norm
of ${\rm top\,}(q)$ follows then by a
partition of unity if one chooses
${\tau}={\rho}-n/2$ (assuming that ${\rho} <
5n/2+1$).

\par The second statement follows from (3.70).
\hfill{$\#$}\medskip

\par With $q\in \widetilde{S}(\langle \cdot
\rangle ^m)$, we put
\ekv{3.75}
{{\rm Top\,}(q)=\Pi _\Phi {\rm top\,}(q)\Pi
_\Phi .}
Then
\ekv{3.76}
{{\rm Top\,}(q)={\rm top\,}(q)+e^{\Phi
/h}(\langle \cdot \rangle
^m\hbox{-negl.})e^{-\Phi /h}.} Notice that if
$m<-2n$, then the trace of
${\rm Top\,}(q)$ is  independent of whether
we view our operator as acting in $H_\Phi $
or in $L^2(e^{-2\Phi /h}L(dx))$.
\par If we choose $q$ as in (3.71), then
\ekv{3.77}
{{\rm Op\,}(p)={\rm Top\,}(q)+e^{\Phi
/h}(\langle \cdot \rangle
^m\hbox{-negl.})e^{-\Phi /h}.}

\par We can now adapt the discussion in
section 2. Assume that $p(x,\xi )=1+a(x,\xi
)$ with $a={\cal O}(\langle (x,\xi )\rangle
^m)$ in $\Lambda _{\Phi _0}+\{ 0\} \times W$
for, $m<-2n$. If we first assume that $p\ne
0$ on $\Lambda _\Phi $, we see that ${\rm
Op\,}(p):H_\Phi \to H_\Phi $ has a \bdd{}
inverse
\ekv{3.78}
{{\rm Op\,}(p)^{-1}={\rm Top\,}(r)+e^{\Phi
/h}({\rm negl.})e^{-\Phi /h},}
where $r\in\widetilde{S}_{{\rm cl}}(1)$,
$r_0(x,x)=1/p(x,{2\over i}\partial _x\Phi
(x))$. The same holds for ${\rm Op\,}(p^t)$,
if we define the deformation $p^t$ from 1 to
$p$ as in section 2. Using the calculus
above, we get the analogue of (2.7):
\ekv{3.79}
{
\log\det {\rm Op\,}(p)={1\over (2\pi
h)^n}(\int_{\Lambda _\Phi}\log p(\rho )\mu
(d\rho )+{\cal O}(h)),  }
where $\mu $ is the symplectic volume density
on $\Lambda _\Phi $. The subsequent argument
of section 2 can also be carried over, so if
we now allow $p$ to vanish on $\Lambda
_\Phi $, we get
\ekv{3.80}
{
\log \vert \det {\rm Op\,}(p)\vert \le
{1\over (2\pi h)^n}(\int_{\Lambda _\Phi }\log
\vert p(\rho )\vert \mu (d\rho )+o(1)),\ h\to
0. }

\par We can now state the main result of this
section.

\medskip
\par\noindent \bf Theorem 3.6. \it Let
$p(x,\xi )$ be holomorphic in an open tubular
neighborhood of ${\bf R}^{2n}$ of the form
${\bf R}^{2n}+iW$, $0\in W\subset\subset {\bf
R}^{2n}$ and satisfying $p(x,\xi )-1={\cal
O}(\langle (x,\xi )\rangle ^m)$ there for
some $m<-2n$. Let $F\subset {\bf C}^{2n}$ be
a complex Lagrangian space (i.e. a complex
subspace of dimension $n$ on which $\sigma $
vanishes) which is strictly negative in the
sense that ${1\over 2i} \sigma
(t\wedge\overline{t})<0$, $\forall t\in
F\setminus \{ 0\}$. Let
$\widetilde{W}\subset\subset W$ be an open
convex \neigh{} of $0\in{\bf R}^{2n}$ and let
$\Lambda \subset {\bf
R}^{2n}+i\widetilde{W}$ be an $IR$-\mfld{}
of the form $\{ \rho +\ell (\rho );\, \rho
\in{\bf R}^{2n}\}$, where $\ell\in C_b^\infty
({\bf R}^{2n};F)$. Then ${\sigma
_\vert}_{\Lambda _t}$ is real, where
$\Lambda _t=\{ \rho +t\ell (\rho );\, \rho
\in{\bf R}^{2n}\}$, $0\le t\le 1$,
 and we assume this form to be non-degenerate,
uniformly on each $\Lambda _t$, uniformly for
$t\in[0,1]$. Then
\ekv{3.81}
{
\log \vert \det p^w(x,hD_x)\vert \le {1\over
(2\pi h)^n}(o(1)+\int_\Lambda \log \vert
p\vert \mu (d\rho )),\ h\to 0. }
\medskip
\par\noindent \bf Proof. \rm Let $\kappa
:{\bf C}^{2n}\to {\bf C}^{2n}$ be a complex
linear canoncial transformation with
$\kappa (F)=\{ (x,\xi )\in{\bf C}^{2n};\,
x=0\}$. The negativity of $F$ implies
that this space is of the form
$$\xi ={\partial f\over \partial x},\
x\in{\bf C}^n$$
where $f$ is a holomorphic quadratic form
with ${\rm Im\,}f''<0$. Then it is clear
that $\kappa $ is of the form (3.3) where
$\phi _0$ satisfies (3.2) (noticing that
$\phi _0(0,y)=-f(y)$), and hence that
$\kappa =\kappa _W$ with $W$ as in (3.1).
Then $\kappa ({\bf R}^{2n})=\Lambda _{\Phi
_0}$, where $\Phi _0$ is a st.pl.s.h.
quadratic form. The assumption about the
form of $\Lambda $ implies that $\kappa
(\Lambda )=\Lambda _{\Phi }$ with
$\nabla \Phi -\nabla \Phi _0\in C_b^\infty
$, and more generally, $\kappa (\Lambda
_t)=\Lambda _{\Phi _t}$, $\Phi _t=t\Phi
+(1-t)\Phi _0$. Since ${\sigma
_\vert}_{\Lambda _t}$
is uniformly non-degenerate, we see that
${\partial ^2\Phi_t \over
\partial\overline{x}\partial x}$ is
uniformly non-degenerate and hence
strictly positive, since so is the case for
$t=0$.

\par We recall the metaplectic invariance
(3.7), so that
$W^{-1}p^w(x,hD)W=q^w(x,hD)$, where
$q=p\circ \kappa ^{-1}$. Clearly $\det
p^w(x,hD)=\det q^w(x,hD)$, where we first
consider $q^w(x,hD)$ as a bounded
operator $H_{\Phi _0}\to H_{\Phi _0}$. If
we assume that $\Phi -\Phi _0$ is
bounded, then (3.81) follows from (3.80)
(with $p$ replaced by $q$), since $H_\Phi $
is then the same space as
$H_{\Phi _0}$. In the general case, we
can approach $\Phi $ by a sequence of
$\widetilde{\Phi }$ such that
$\widetilde{\Phi }-\Phi _0 $ is bounded,
$\nabla \widetilde{\Phi }\to \nabla \Phi
$ in $C_b^\infty $,  $\widetilde{\Phi
}$ is uniformly st.pl.s.h. and $\Lambda
_{\widetilde{\Phi }}$ stays inside the
convex tube $\kappa ({\bf
R}^{2n}+i\widetilde{W})$. We then have
(3.80) uniformly (with $p, \Lambda _\Phi$
replaced by $q, \Lambda _{\widetilde{\Phi }}$)
and passing to the limit we get it also for
$\Phi $. \hfill{$\#$}\medskip

\par In [Sj2]
function spaces are defined associated to
IR-manifolds which are obtained by certain
global holomorphic deformations and it is
showed there how $h$-\pseudor s act on these
spaces and that this action can be pulled back
to an action on $H_{\Phi _0}$. Under
essentially the same assumptions on $P$, we
still have the conclusion (3.81).

\bigskip

\centerline{\bf 4. H{\"o}lder properties of
$I(\Lambda ,p)$, $I_\epsilon (\Lambda ,p)$;
the differential w.r.t. $\Lambda $.}
\medskip
\par Let $\Phi _0$ and $\Lambda _{\Phi _0}$
be as in section 3. Let
$\widetilde{W}\subset\subset W$ be open
bounded neighborhoods of $0$ in $\Lambda
_{\Phi _0}$. Let $p$ be holomorphic in
$\Lambda _{\Phi _0}+iW$ and satisfy
\ekv{4.1}
{
p(x,\xi )-1={\cal O}(\langle (x,\xi \rangle
^m), \hbox{ for some }m<-2n. }
The IR \mfld s $\Lambda $ under consideration
in the remainder of this paper will satisfy:
\ekv{4.2}
{
\Lambda \subset \Lambda _{\Phi
_0}+i\widetilde{W}\hbox{ is closed}. }
\ekv{4.3}
{\Lambda \hbox{ is diffeomorphic to }{\bf
R}^{2n}.}
\eekv{4.4}
{\hbox{Outside a bounded set, we have }\Lambda
=\Lambda _\Psi , }
{\hbox{with }\Psi -\Phi _0\in C_b^\infty
,\ (\partial _{\overline{x}}\partial _x\Psi
)^{-1}={\cal O}(1).}
(From the discussion in section 3 and in this
section, it will be clear that the assumption
on $\Psi -\Phi _0$ can be weakend to: $\Psi
'-\Phi _0'\in C_b^\infty $.)
\par The last assumption allows us to
identify $\Lambda $ near infinity with a
neighborhood of infinity of ${\bf C}^n$ by
means of ${\bf C}^n\ni x\simeq (x,{2\over
i}{\partial \Psi \over \partial x }(x))\in
\Lambda $, and we can therefore define the
space $C_b^\infty (\Lambda )$.
\par Sometimes we let $\Lambda =\Lambda _t$
depend smoothly on $t\in{\rm neigh\,}(0;{\bf
R}^k)$. We then require $\Lambda _t$ to
fulfill (4.2)--(4.4) with uniformity in
(4.2), (4.4). Moreover, if we represent
$\Lambda =\Lambda _{\Psi _t}$ outside a
bounded set, then we assume that for every
$N\ge 0$,
\ekv{4.5}
{\nabla _t^N\Psi _t''\in C_b^\infty
,\hbox{ uniformly in }t.}
 \par We will write $p_\Lambda $ or
sometimes only $p$ for the restriction
${p_\vert}_{\Lambda }$. The zeros of $p$ in
$\Lambda _{\Phi _0}+i\widetilde{W}$ are
confined to a bounded set and there exists
$m_0\in{\bf N}$, such that for every $\rho
\in\Lambda _{\Phi _0}+i\widetilde{W}$,
there exists $m\in \{ 0,1,..,m_0\}$, such
that $\nabla ^mp (\rho )\ne 0$. If $\Lambda
$ is an IR \mfld{ } as above, we conclude
(since $\Lambda $ is maximally totally
real) that for every $\rho \in\Lambda $, we
have $\nabla ^mp_\Lambda (\rho )\ne 0$ for
some $0\le m\le m_0$. Let $\rho _0\in\Lambda
$ be a zero of $p_\Lambda $ and choose $m$
as above. Then we can find a real smooth
vectorfield $\nu $ on ${\rm neigh\,}(\rho
_0,\Lambda )$ such that $\nu ^mp_\Lambda
(\rho _0)\ne 0$. Choose local coordinates
$(x_1,x_2,..,x_{2n})=(x_1,x')$ for $\Lambda $
centered at
$\rho _0$, such that $\nu ={\partial \over
\partial x_1}$. By Malgrange's preparation
theorem ([Ma]), we have in a neighborhood of
$x=0$:
\ekv{4.6}
{p_\Lambda (x)=q(x)(x_1^m+\sum_1^m
a_j(x')x_1^{m-j})=q(x)\prod_{j=1}^m(x_1-\lambda
_j(x')),}
where $a_j(x')$, $q(x)$ are smooth with
$a_j(0)=0$, $q(0)\ne 0$, and $\lambda
_j(x')$ are the zeros of the second factor
in the middle expression. (The ordering of
the roots does not matter.)  If we let
$\Lambda =\Lambda _t$ depend smoothly on
$t\in{\rm neigh\,}(0,{\bf R}^k)$ as above,
and let $\rho _0$ be a zero of $p_{\Lambda
_0}$, and if we choose $x_1,..,x_{2n}$ to be
local coordinates on $\Lambda _t$ depending
smoothly on $t$ and centered at $\rho _0$
when $t=0$, then (4.6) extends to $p_{\Lambda
_t}$ with $q=q(t,x)$, $a_j=a_j(t,x')$ being
smooth in $(t,x)$, and with $\lambda
_j=\lambda _j(t,x')$. From (4.6), we get
easily:
\medskip
\par\noindent \bf Lemma 4.1. \it For the
choice of coordinates above and for $\epsilon
_0>0$ small enough, we have
\ekv{4.7}
{
\lambda (\{ x_1\in [-\epsilon _0,\epsilon
_0];\, \vert p_\Lambda (x)\vert \le \delta \}
)\le C\delta ^{1/m}, }
for $\vert x'\vert \le \epsilon _0$ and
$0<\delta \le 1$. Here $\lambda $ denotes the
one dimensional Lebesgue measure. In the
parameter dependent case, the estimate is
uniform in $t$ for $\vert t\vert \le \epsilon
_0$ (small enough).
\medskip
\par\noindent \bf Proof. \rm It suffices to
observe that $\{ x_1\in [-\epsilon
_0,\epsilon _0];\, \vert p_\Lambda (x)\vert
\le \delta \}$ is contained in the union of
the intervals $[{\rm Re\,}\lambda
_j(x')-C\delta ^{1/m},{\rm Re\,}\lambda
_j(x')+C\delta ^{1/m}]$.\hfill{$\#$}
\medskip

\par If $\mu (d\rho )=\mu _\Lambda (d\rho )$
denotes the symplectic volume element on
$\Lambda $, it follows that
\ekv{4.8}
{
\mu (\{ \rho \in \Lambda ;\, \vert p_\Lambda
(\rho )\vert \le \delta \} )\le C\delta
^{1/m_0}, }
again locally uniformly in $t$ in the
parameter dependent case.

\par For $0<\epsilon \le 1$, we put
\ekv{4.9}
{
I_\epsilon (\Lambda ,p)={1\over
2}\int_\Lambda  \log {\vert p_\Lambda (\rho
)\vert ^2+\epsilon ^2\over 1+\epsilon
^2}\,\,\mu (d\rho ). }
Notice that the integral converges in view of
(4.1), (4.4). We also define $I(\Lambda
,p)=I_0(\Lambda ,p)$:
\ekv{4.10}
{
I(\Lambda ,p)=\int_\Lambda  \log \vert
p_\Lambda (\rho )\vert \,\, \mu (d\rho ). }
Here the convergence over a neighborhood of
infinity follows from (4.1), (4.4) but the
possible presence of zeros of $p_\Lambda $
requires a verification of the convergence of
\ekv{4.11}
{
\int_{\vert p_\Lambda \vert ^{-1}([0,{1\over
2}])}\log \vert p_\Lambda (\rho )\vert \,
\mu (d\rho )=\int_0^{1/2}\log t \,\,dV(t), }
where $V(t)=\mu (\vert p_\Lambda \vert
^{-1}([0,t]))={\cal O}(t^{1/m_0})$, by (4.8).
Integration by parts shows that the last
integral is equal to
$$[(\log t)V(t)]_0^{1/2}-\int_0^{1/2}{1\over
t}V(t) dt.$$

\par From the dominated convergence theorem,
we see that
\ekv{4.12}
{
I_\epsilon (\Lambda ,p)\to I(\Lambda ,p),\
\epsilon \to 0. }
Using (4.8), we can estimate the rate of
convergence.
\medskip
\par\noindent \bf Lemma 4.2. \it We have
\ekv{4.13}
{
I_\epsilon (\Lambda ,p)-I(\Lambda ,p)={\cal
O}(\epsilon ^{1/m_0}),\ \epsilon \to 0. }
\medskip
\par\noindent \bf Proof. \rm For $\epsilon
>0$, we have
\ekv{4.14}
{
\partial _\epsilon I_\epsilon (\Lambda
,p)=\int_\Lambda [{\epsilon \over \vert
p_\Lambda (\rho )\vert ^2+\epsilon
^2}-{\epsilon \over 1+\epsilon ^2}]\mu
(d\rho ). }
The contribution to the integral from the
region where $\vert p_\Lambda (\rho )\vert
>1/2$ is ${\cal O}(\epsilon )$ and the
contribution from the bounded region where
$\vert p_\Lambda (\rho )\vert \le 1/2$ is
$$\eqalign{
&{\cal O}(\epsilon )+\int_{\vert p_\Lambda
\vert ^{-1}([0,{1\over 2}])}{\epsilon \over
(\epsilon ^2+\vert p_\Lambda (\rho )\vert
^2)}\,\mu (d\rho )= {\cal O}(\epsilon
)+\int_0^{1/2}{\epsilon \over \epsilon
^2+t^2}dV(t)=\cr
&{\cal O}(\epsilon ) +[{\epsilon \over
\epsilon
^2+t^2}V(t)]_0^{1/2}+\int_0^{1/2}{2\epsilon
t\over (\epsilon ^2+t^2)^2}V(t)dt=\cr
&{\cal
O}(\epsilon )+{\cal
O}(1)\int_0^{1/2}{2t\epsilon \over (\epsilon
^2+t^2)^2}t^{1/m_0}dt={\cal O}(\epsilon
^{{1\over m_0}-1}).}$$
Hence $\partial _\epsilon I_\epsilon (\Lambda
,p)={\cal O}(\epsilon ^{-1+1/m_0})$ and the
Lemma follows by
integration.\hfill{$\#$}\medskip

\par We next study the differential of
$I_\epsilon $ and $I$ with respect to
$\Lambda $. Let ${\rm neigh\,}(0,{\bf R})\ni
t\mapsto
\Lambda _t$ be a smooth family of IR-\mfld s
as above, satisfying (4.1)--(4.5). We then
have a corresponding generator $q_t\in
C_b^\infty (\Lambda _t;{\bf R})$ such that a
deformation field of the family is given by
\ekv{4.15}
{
\nu _t=H_{\widetilde{q}_t}^{{\rm Im\,}\sigma
}, }
when $\widetilde{q}_t$ is a smooth real
extension of $q_t$ to a full neighborhood of
$\Lambda _t$. If $f_t$ denotes an almost
holomorphic extension of $q_t$, we take
$\widetilde{q}_t={\rm Re\,}f_t$, so that
(cf. (1.19)):
\ekv{4.16}
{
\nu _t=\widehat{iH_{f_t}},\hbox{ on }\Lambda
_t. }
With
$$F_\epsilon (p)={1\over 2}\log
({p\overline{p}+\epsilon ^2\over 1+\epsilon
^2}),$$
we have
$$I_\epsilon (\Lambda _t,p)=\int_{\Lambda
_t}F_\epsilon (p)\mu (d\rho ),$$
and
\ekv{4.17}
{
\partial _tI_\epsilon (\Lambda
_t,p)=\int_{\Lambda _t}({\partial F_\epsilon
\over \partial p}\nu _t(p)+{\partial
F_\epsilon \over \partial
\overline{p}}\overline{\nu _t(p)})\mu (d\rho
), }
where we also used that $\nu _t$ in (4.16)
induces canonical transformations $\kappa
_{t,s}:\Lambda _s\to\Lambda _t$, which
conserve the symplectic volume element (see
(1.19)). Since
$p$ is holomorphic, we have $\nu
_t(p)=iH_{f_t}p$, where in the last
expression, we may view
$H_{f_t}=\widehat{H_{f_t}}$ as the real
Hamilton field on $\Lambda _t$ of
${{f_t}_\vert }_{\Lambda _t}=q_t$ with
respect to ${\sigma _\vert}_{\Lambda _t}$.
>From (4.17), we get
\eekv{4.18}
{
\partial _tI_\epsilon (\Lambda
_t,p)=\int_{\Lambda _t}({\partial F_\epsilon
\over \partial p}(p) iH_{f_t}p-{\partial
F_\epsilon \over \partial
\overline{p}}(p)iH_{f_t}\overline{p})\mu
(d\rho ) }
{
={i\over 2}\int_{\Lambda _t}\langle
{\overline{p}\over p\overline{p}+\epsilon
^2}dp-{p\over p\overline{p}+\epsilon
^2}d\overline{p},H_{f_t}\rangle \mu (d\rho )
={i\over 2}\int_{\Lambda _t}\langle
{p\overline{p}\over p\overline{p}+\epsilon
^2}({dp\over p}-{d\overline{p}\over
\overline{p}}),H_{f_t}\rangle \mu (d\rho ). }
Here
\ekv{4.19}
{
{1\over 2i}({dp\over p}-{d\overline{p}\over
\overline{p}})=d\,{\rm arg\,}p, }
where $p=p_\Lambda \ne 0$, so we get
\ekv{4.20}
{
\partial _tI_\epsilon (\Lambda
_t,p)=-\int_{\Lambda _t}\langle
{p\overline{p}\over p\overline{p}+\epsilon
^2}d({\rm arg\,}p_{\Lambda
_t}),H_{f_t}\rangle \mu (d\rho ), }
which can also be written
\ekv{4.21}
{
\partial _tI_\epsilon (\Lambda
_t,p)=\int_{\Lambda _t}\langle
{p\overline{p}\over p\overline{p}+\epsilon
^2}H_{{\rm arg\,}p_{\Lambda _t}},df_t\rangle
\mu (d\rho ). }
In the last two integrals the integration is
restricted to $\Lambda _t\setminus
p^{-1}(0)$, and the Hamilton fields are the
ones for the real symplectic structure on
$\Lambda _t$.

\medskip
\par\noindent \bf Proposition 4.3. \it The
coefficients of the differential form
$d\,{\rm arg\,}p_{\Lambda _t}$ (defined on
$\Lambda _t\setminus p^{-1}(0)$) belong to
$L^1(\Lambda _t,\mu )$.
\medskip
\par\noindent \bf Proof. \rm The
integrability near infinity follows from
(4.1) (and we there identify $\Lambda _t$ with
a neighborhood of infinity in ${\bf C}^n$ as
indicated after (4.4)), so we only have to
establish the local integrability near every
point in
$p^{-1}(0)\cap\Lambda _t$. Let $\rho _0$ be
such a point and let $m\in\{ 1,2,..\}$ be
the order of vanishing, so that $\nabla
^kp(\rho _0)=0$, $0\le k\le m-1$, $\nabla
^mp(\rho _0)\ne 0$. Choose $2n$ linearly
independent and commuting \vf s $\nu
_1,..,\nu _{2n}$ in a neighborhood of $\rho
_0$ such that $\nu _j^mp(\rho _0)\ne 0$. Let
$x_1,..,x_{2n}$
be the local coordinates centered at $\rho
_0$ with $\nu _j={\partial\over \partial
x_j}$. Let
$A_j={\partial \over \partial x_j}{\rm
arg\,}p_{\Lambda _t}$ for $x\in{\rm
neigh\,}(0)\setminus p_{\Lambda _t}^{-1}(0)$,
so that $d\,{\rm arg\,}p_{\Lambda
_t}=\sum_1^{2n}A_j(x)dx_j$. We have to show
that $A_j$ is locally integrable near $0$ and
may choose $j=1$ for simplicity. Applying
(4.6) to $p_\Lambda =p_{\Lambda _t}$, we see
that
$$\int_{\{ x_1\in [-\epsilon _0,\epsilon
_0];\, p_{\Lambda _t}(x)\ne 0\} }\vert
{\partial \over \partial x_1}{\rm
arg\,}p_{\Lambda _t}(x)\vert dx_1\le {\rm
Const.},$$
for $\vert x'\vert \le \epsilon _0$, and the
local integrability follows.\hfill{$\#$}
\medskip

\par Let $[d\,{\rm arg\,}p]_{\Lambda _t}$
denote the differential form on $\Lambda _t$
with $L^1$ coefficients whose restriction to
$\Lambda _t\setminus p^{-1}(0)$ is equal to
$d\,{\rm arg\,}p_{\Lambda _t}$. Similarly
let $[H_{{\rm arg\,}p}]_{\Lambda _t}$ be the
$L^1$ \vf{} on $\Lambda _t$ determined by
$H_{{\rm arg\,}p_{\Lambda _t}}$ on $\Lambda
_t\setminus p^{-1}(0)$. Write
\ekv{4.22}
{
\partial _tI_\epsilon (\Lambda
_t,p)=-\int_{\Lambda _t}\langle
{p\overline{p}\over p\overline{p}+\epsilon
^2} [d\,{\rm arg\,}p]_{\Lambda
_t},H_{f_t}\rangle \mu (d\rho ), }
\ekv{4.23}
{
\partial _tI_\epsilon (\Lambda
_t,p)=\int_{\Lambda _t}\langle
{p\overline{p}\over p\overline{p}+\epsilon
^2}[H_{{\rm arg\,}p}]_{\Lambda
_t},df_t\rangle \mu (d\rho ). }
$H_{f_t}$, $df_t$ have bounded coefficients
near infinity (using the identification
with ${\bf C}^n$ near infinity indicated
after (4.4)), and from Proposition 4.3 and the
dominated convergence theorem, we get the
first part of the following result:
\medskip
\par\noindent \bf Proposition 4.4. \it We
have $\partial _tI_\epsilon (\Lambda _t,p)\to
{\rm "}\partial _tI(\Lambda _t,p){\rm "}$,
$\epsilon \to 0$, where
\ekv{4.24}
{
{\rm "}\partial _tI(\Lambda
_t,p){\rm "}=-\int_{\Lambda _t}\langle
[d\,{\rm arg\,}p]_{\Lambda
_t},H_{f_t}\rangle \mu (d\rho )=\int_{\Lambda
_t}\langle [H_{{\rm arg\,}p}]_{\Lambda
_t},df_t\rangle \mu (d\rho ). }
Moreover, $t\mapsto I(\Lambda _t,p)$ is a
Lipschitz function and the a.e.{ }derivative
is given by ${\rm "}\partial _tI(\Lambda
_t,p){\rm "}$.\medskip\rm

\par Here the last statement follows from the
fact that $t\mapsto \Vert [d\, {\rm
arg\,}p]_{\Lambda _t}\Vert _{L^1}$ is
locally bounded, as we can see by adding a
smooth $t$ dependence in (4.6).

\par A Hamilton field is divergence free in
symplectic coordinates, so ${\rm
div\,}[H_{{\rm arg\,}p}]_{\Lambda _t}$
is a distribution of order $\le 1$ with
support in $p^{-1}(0)\cap \Lambda _t$. From
(4.24), we get
\ekv{4.25}
{
{\rm "}\partial _t(\Lambda _t,p){\rm
"}=-\int_{\Lambda _t}f_t\,{\rm div\,}[H_{{\rm
arg\,}p}]_{\Lambda _t}\mu (d\rho ). }
Notice that the integral does not change if
we add a $t$-dependent constant to $f_t$. When
$dp_{\Lambda },\,
\overline{dp_{\Lambda }}$ are pointwise
linearly independent on
$p_\Lambda ^{-1}(0)$, we shall obtain in
section 8  a simple expression for this
distribution and see that it is a Radon
measure.

\bigskip

\centerline{\bf 5. Second derivative under
non-autonomous flows.}
\medskip
\par Let $I\ni t\mapsto
\Lambda _t$ be a smooth deformation of
IR-manifolds as in the preceding section, where
$I$ is a small open interval which contains
$0$. Let
$f_t\in C_b^\infty (\Lambda _t;{\bf R})$ be a
corresponding smooth family of generating
functions (unique up to a $t$ dependent
constant) so that if we extend $f_t$ to an
almost holomorphic function in a neighborhood
of $\Lambda _t$, then a corresponding
deformation field is given at $\Lambda _t$ by
\ekv{5.1}
{
\nu _t=\widehat{iH_{f_t}}=2{\rm
Re\,}(iH_{f_t}),}
and we recall that $\nu _t$ generates a
family $\kappa _{t,s}:\Lambda _s\to \Lambda
_t$ of canonical diffeomorphisms.

\par Let $p$ be a holomorphic function as in
the preceding section and let
$F$ be a smooth real-valued function defined on
a complex neighborhood of $\overline{p(\Lambda
_I)}$, with $F(1)=0$, where $\Lambda
_I=\cup_{t\in I}\Lambda _t$. Put
\ekv{5.2}
{J(t)=\int_{\Lambda _t}F(p)\mu (d\rho ).}

With $p_t=p\circ \kappa _t$, $\kappa
_t=\kappa _{t,0}$, we get
$$J(t)=\int_{\Lambda _0}F(p_t)\mu (d\rho ),$$
and
$$\eqalign{\partial _tJ(t)&=\int_{\Lambda
_0}[{\partial F\over \partial p}(p_t)\partial
_tp_t+{\partial F\over \partial
\overline{p}}(p_t)\overline{\partial
_tp_t}]\mu (d\rho )\cr
&=\int_{\Lambda _t}[{\partial F\over
\partial p}(p)iH_{f_t}p-{\partial F\over
\partial
\overline{p}}i\overline{H_{f_t}p}]
\mu (d\rho ).}$$
Since $p$ is holomorphic near $\Lambda
_t$, we can consider (as in the preceding
section) $H_{f_t}$ as the real Hamilton field
$H_{f_t}^{{\sigma _\vert}_{\Lambda _t}}$.
Viewing in the same way $H_p$,
$H_{\overline{p}}$ as complex Hamilton fields
on $\Lambda _t$, we get
$$\eqalign{
\partial _tJ(t)&=i\int_{\Lambda _t}[{\partial
F\over \partial p}(p)H_{f_t}(p)-{\partial
F\over \partial
\overline{p}}(p)H_{f_t}(\overline{p})]\mu
(d\rho )\cr
&={1\over i} \int_{\Lambda _t}[{\partial
F\over \partial p}(p)H_p(f_t)-{\partial
F\over \partial
\overline{p}}(p)H_{\overline{p}}(f_t)]\mu
(d\rho ).
 }$$
Since the transpose of $H_p$
with respect to $\mu (d\rho )$
is equal to $-H_p$, we get
\ekv{5.3}
{
\partial _tJ(t)=\int_{\Lambda _t}G_t(\rho
)f_t(\rho )\mu (d\rho ),\ G_t(\rho
)=i(H_p({\partial F\over \partial
p})-H_{\overline{p}}({\partial F\over
\partial \overline{p}})). }
Notice that $G_t(\rho )$ is defined on
$\Lambda _t$ with $H_p, H_{\overline{p}}$
defined as complex vectorfields on $\Lambda
_t$, so $G_t(\rho )$ depends both on $\rho $
and on $T_\rho \Lambda _t$.

\par Let $I\ni t\mapsto \widetilde{\Lambda
}_t$ be a second family with the same
properties and let $\widetilde{f}_t$ be a
corresponding smooth family of generating
fucntions, that we also extend almost
holomorphically. We assume that
\ekv{5.4}
{\widetilde{\Lambda }_0=\Lambda _0,}
\ekv{5.5}
{\widetilde{f}_0=f_0.}
Equivalently, we have
\ekv{5.6}
{
{\rm dist\,}(\widetilde{\Lambda }_t,\Lambda
_t)={\cal O}(t^2). }
Possibly after shrinking $I$ around $0$, we
can define
\ekv{5.7}
{\widetilde{J}(t)=\int_{\widetilde{\Lambda
}_t}F(p)\mu (d\rho ),}
and analogously to (5.3), we have
\ekv{5.8}
{
\partial
_t\widetilde{J}(t)=\int_{\widetilde{\Lambda
}_t}\widetilde{G}_t(\rho
)\widetilde{f}_t(\rho )\mu (d\rho ). }

\par Let $\Phi _t=\widetilde{\kappa}_t\circ
\kappa _t^{-1}:\Lambda _t\to \widetilde{\Lambda
}_t$, where $\widetilde{\Lambda
}_t=\widetilde{\kappa}_t(\Lambda _0)$, so
that $\Phi _t$ is symplectic and
\ekv{5.9}
{
\Phi _t(\rho )-\rho ={\cal O}(t^2).
}
Since we also have ${\rm dist\,}(T_{\Phi
_t(\rho )}(\widetilde{\Lambda }_t),T_\rho
(\Lambda _t))={\cal O}(t^2),$
we get $\widetilde{G}_t(\Phi _t(\rho
))=G_t(\rho )+{\cal O}(t^2),$ and hence
\ekv{5.10}
{
\partial _t\widetilde{J}(t)-\partial
_tJ(t)=\int_{\Lambda _t}G_t(\rho ) g_t(\rho
)\mu (d\rho )+{\cal O}(t^2), }
where
\ekv{5.11}
{
g_t=\widetilde{f}_t\circ \Phi _t-f_t,
}
and we can assume that $g_t(\rho )\in
C^\infty (I\times {\rm neigh\,}(\Lambda _0))$
with $g_t(\cdot )$ almost holomorphic on
$\Lambda _t$. We can also assume that
\ekv{5.12}
{
g_0(\rho )\equiv 0.}
(5.10) implies that
\eekv{5.13}
{
\partial _t^2\widetilde{J}(0)-\partial
_t^2J(0)=\int_{\Lambda _0}G_0(\rho
)(\partial _tg)_{t=0}(\rho )\mu (d\rho ) }
{
=i\int_{\Lambda _0}[H_p({\partial F\over
\partial p})-H_{\overline{p}}({\partial
F\over \partial \overline{p}})]((\partial
_t\widetilde{f}_t)_{t=0}-(\partial
_tf_t)_{t=0})(\rho )\mu (d\rho ), }
and where we used (5.9) in the last step. Here
we can undo the previous integration by parts
and get
\ekv{5.14}
{
\partial _t^2\widetilde{J}(0)-\partial
_t^2J(0)={1\over i}\int_{\Lambda _0}
[{\partial F\over \partial p}H_p(\partial
_t\widetilde{f}-\partial _tf)_{t=0}-{\partial
F\over \partial
\overline{p}}H_{\overline{p}}(\partial
_t\widetilde{f}-\partial _tf)_{t=0}]\mu
(d\rho ). }
Notice that $(\partial _tf_t)_{t=0}=\partial
_t(f(t,\kappa _t(\rho )))_{t=0}$ is real, and
similarly for $\partial _t\widetilde{f}_t$.

\par In order to exploit the last relation,
we need to compute the second derivative in
the case of "autonomous flows", and we shall
temporarily consider the case of a special
family of deformations. Let $\widetilde{f}\in
C_b^\infty (\Lambda _0;{\bf R})$ and denote
by $\widetilde{f}$
also an almost holomorphic extension of
$\widetilde{f}$. For $t\in {\rm neigh\,}(0,{\bf
C})$, consider the IR-\mfld s
$\widetilde{\Lambda }_t=\widetilde{\kappa }
_t(\Lambda _0)$, where
$\widetilde{\kappa } _t=\exp H_{{\rm
Re\,}(t\widetilde{f})}^{{\rm Im\,}\sigma }$.
To infinite order on $\Lambda _0$, we have
$$H_{{\rm Re\,}(t\widetilde{f})}^{{\rm
Im\,}\sigma} \equiv
\widehat{itH_{\widetilde{f}}}\equiv -H_{{\rm
Im\,}(t\widetilde{f})}^{{\rm Re\,}\sigma },$$
and it follows that $\widetilde{\kappa }
_t^*\sigma =\sigma +{\cal O}(t^\infty )$ on
$\Lambda _0$. Let
$\partial _t$, $\partial _{\overline{t}}$
denote holomorphic and antiholomorphic
derivatives. We also see that
$\partial _{\overline{t}}\widetilde{\kappa }
_t(\rho )={\cal O}(t^\infty )$ on $\Lambda
_0$. We have
$$\widetilde{J}(t):=\int_{\widetilde{\Lambda
}_t}F(p)\mu (d\rho )=\int_{\Lambda
_0}F(p_t)\mu (d\rho )+{\cal O}(t^\infty ),\
p_t:=p\circ \widetilde{\kappa } _t,$$
$$\eqalign{
\partial _t\widetilde{J}(t)=\int_{\Lambda
_0}\partial _t(F(p_t))\mu (d\rho )+{\cal
O}(t^\infty )&=\cr
 \int_{\Lambda _0}({\partial F\over \partial
p}(p_t)\partial _tp_t+{\partial F\over
\partial \overline{p}}(p_t)\overline{\partial
_{\overline{t}}p_t})\mu (d\rho )+{\cal
O}(t^\infty )&=\int_{\Lambda _0}{\partial
F\over \partial p}(p_t)\partial _tp_t\, \mu
(d\rho )+{\cal O}(t^\infty ),
 }$$
$$
\partial _{\overline{t}}\partial
_t\widetilde{J}(t)=\int_{\Lambda _0}{\partial
^2 F\over \partial \overline{p}\partial
p}(p_t)\overline{\partial _tp_t}\partial
_tp_t\, \mu (d\rho )+{\cal O}(t^\infty ).
$$
For $t=0$ we have $\partial
_tp_t=iH_{\widetilde{f}}(p)$ and hence
\ekv{5.15}
{
(\partial _{\overline{t}}\partial
_t\widetilde{J}(t))_{t=0}=\int_{\Lambda _0}{
\partial ^2F\over \partial
\overline{p}\partial
p}(p)(H_{\widetilde{f}}p)(\overline{
H_{\widetilde{f }}(p)})\,\mu (d\rho ). }

\par Next we compute
$$\eqalign{
\partial _{{\rm
Im\,}t}\widetilde{J}(t)&=i(\partial
_t\widetilde{J}(t)-\partial
_{\overline{t}}\widetilde{J}(t))=i\int_{
\widetilde{\Lambda } _0} ({\partial F\over
\partial p}(p_t)\partial _tp_t-{\partial F\over
\partial \overline{p}}(p_t)\overline{\partial
_tp_t})\mu (d\rho )+{\cal O}(t^\infty )\cr
&=i\int_{\widetilde{\Lambda }_t}({\partial
F\over
\partial p}(p)iH_{\widetilde{f}}p-{\partial F
\over \partial
\overline{p}}(p)\overline{iH_{\widetilde{f}}p})\,\mu
(d\rho )+{\cal O}(t^\infty ). }$$
Modulo ${\cal O}(t^\infty )$, we have
$H_{\widetilde{f}}p=H_{\widetilde{f}}^
{{\sigma _\vert}_{\tilde{\Lambda }_t}}p$,
${\rm Im\,}\widetilde{f}=0$, on
$\widetilde{\Lambda }_t$, so we get
\eekv{5.16}
{
\partial _{{\rm
Im\,}t}\widetilde{J}(t)=-\int_{
\widetilde{\Lambda }_t}({\partial F\over
\partial p}(p)H_{{\rm
Re\,}\tilde{f}}^{{\sigma
_\vert}_{\tilde{\Lambda }_t}}p+{\partial
F\over \partial
\overline{p}}H_{{\rm
Re\,}\tilde{f}}^{{\sigma
_\vert}_{\tilde{\Lambda }_t}}\overline{p})\mu
(d\rho )+{\cal O}(t^\infty ) }
{
\hskip 2truecm =-\int_{\widetilde{\Lambda
}_t} H_{{\rm
Re\,}\tilde{f}}^{{\sigma
_\vert}_{\tilde{\Lambda }_t}}(F(p))\mu
(d\rho )+{\cal O}(t^\infty )={\cal
O}(t^\infty ). }
(The idea behind this is that in the case
when $\widetilde{f}$ is holomorphic, then
$\widetilde{\Lambda }_t=\exp
\widehat{itH_{\tilde{f}}}(\Lambda _0)$ only
depends on ${\rm Re\,}t$, and hence we have
the same for $\widetilde{J}(t)$.) The last
relation can be differentiated any number of
times w.r.t. ${\rm Re\,}t$, ${\rm Im\,}t$ and
combining this with (5.15), we get
\ekv{5.17}
{
(\partial _{{\rm
Re\,}t}^2\widetilde{J}(t))_{t=0}=
4(\partial _{\overline{t}}\partial
_t\widetilde{J}(t))_{t=0}=4\int_{\Lambda _0}
{\partial ^2F\over \partial
\overline{p}\partial
p}(p)H_{\tilde{f}}(p)
\overline{H_{\tilde{f}}(p)}\mu (d\rho ).
}

\par We now restrict the attention to real
$t$ and write $\partial _t$ instead of
$\partial _{{\rm Re\,}t}$. In order to apply
(5.17) in (5.13), we let
$\widetilde{f}(t,\rho )\in C_b^\infty
(I\times {\bf C}^{2n})$ be such that
$\widetilde{f}(t,\cdot )$ is an almost
holomorphic extension of
${{}{\widetilde{f}}_\vert }_{ \tilde{\Lambda
}_t}$ and observe that $\partial
_t\widetilde{f}(t,\rho )={\cal O}(t^\infty
+{\rm dist\,}(\rho ,\Lambda _0)^\infty )$. We
can apply (5.13) and get the following
identity for $t=0$:
\eekv{5.18}
{
\partial _t^2J(t)=4\int_{\Lambda _t}{\partial
^2F\over \partial \overline{p}\partial
p}H_{f_t}(p)H_{f_t}(\overline{p})\mu (d\rho ) }
{
\hskip 5cm +i\int_{\Lambda _t } [H_p({\partial
F\over \partial
p})-H_{\overline{p}}({\partial F\over
\partial \overline{p}})]\partial _tf_t \mu
(d\rho ). }
Comparing with (5.3), we see that the last
term in (5.18) is equal to ${(\partial
_s)}_{s=0}J(t,s)$, with $J(t,s)=\int_{\Lambda
_{t,s}}F(p)\mu (d\rho )$, $\Lambda
_{t,s}=\exp s2{\rm Re\,}(i H_{\partial
_tf_t})(\Lambda _t)$. Clearly this extends to
general
$t$ and the last term can also be written
\ekv{5.19}
{
{1\over i}\int_{\Lambda _t}({\partial F\over
\partial p}H_p-{\partial F\over \partial
\overline{p}}H_{\overline{p}})\partial
_tf_t\mu (d\rho ), }
where $H_p$, $H_{\overline{p}}$ in the last two
formulae denote the complex vector fields on
$\Lambda _t$.

\par Notice that the bracket in the last
integral in (5.18) can also be written as
$$H_p({\partial F\over \partial
p})-H_{\overline{p}}({\partial F\over
\partial \overline{p}})=2{\partial ^2F\over
\partial p\partial \overline{p}}\{
p,\overline{p}\},$$
where $\{ p,\overline{p}\}=H_p(\overline{p})$
denotes the Poisson bracket of $p$ and
$\overline{p}$ for the symplectic structure of
$\Lambda _t$.

\bigskip
\centerline{\bf 6. Continuity and
convergence for the differential of $\Lambda
\mapsto I(\Lambda ,p)$.}
\medskip
\par Let $\Lambda $
be an IR-manifold and $p$ a holomorphic
function  as in section 4. We recall (see
also below) that by Malgrange's preparation
theorem, the  differential form
$d\,{\rm arg\,}p$ on $\Lambda
\setminus p^{-1}(0)$ has $L^1$ coefficients,
and the same holds for the
Hamilton field $H_{{\rm arg\,}p}$. Let
$[d\,{\rm arg\,}p]_{\Lambda _t}$
and $[H_{{\rm arg\,}p}]_{\Lambda _t}$ denote
the corresponding $L^1$ form and $L^1$
vectorfield respectively on $\Lambda $.

We review some calculations for the
differential of $\Lambda \mapsto I(\Lambda
,p)$. Let $t\mapsto \Lambda _t$ be a smooth
deformation of IR-manifolds as in section
4, where $I$ is an open bounded interval
containing $t=0$, and let
$t\mapsto f(t,\rho )\in C_b^\infty (\Lambda
_t;{\bf R})$ be a corresponding generating
family (unique up to a $t$ dependent
constant and extended to be almost
holomorphic in $\rho $), so that $\Lambda
_t=\Phi _t(\Lambda _0)$, with $\partial
_t\Phi _t(\rho )=\widehat{iH_{f_t}}(\Phi
_t(\rho ))$, $\rho \in \Lambda _0$. With
$F_\epsilon (p)={1\over 2}\log
({p\overline{p}+\epsilon ^2\over 1+\epsilon
^2})$, we  obtained in section 4 that
\ekv{6.1}
{
I_\epsilon (\Lambda _t,p)=\int_{\Lambda
_t}F_\epsilon (p) \mu (d\rho )=I(\Lambda
_t,p)+{\cal O}(\epsilon ^{1/m}), }
for some $m\in{\bf N}\setminus\{ 0\}$.
Moreover, we saw that
$$\partial _tI_\epsilon (\Lambda
_t,p)  ={i\over 2}\int_{\Lambda _t}\langle
{p\overline{p}\over p\overline{p}+\epsilon
^2}({dp\over p}-{d\overline{p}\over
\overline{p}} ),H_{f_t}\rangle \mu (d\rho ).
$$
Here ${1\over 2i}({dp\over
p}-{d\overline{p}\over \overline{p}})=d\,{\rm
arg\,}p$ where $p\ne 0$, so that
\eekv{6.2}
{
\partial _tI_\epsilon (\Lambda
_t,p)=-\int_{\Lambda _t}\langle
{p\overline{p}\over p\overline{p}+\epsilon
^2}[d({\rm arg\,}p)]_{\Lambda
_t},H_{f_t}\rangle \mu (d\rho ) }
{
=\int_{\Lambda _t}\langle {p\overline{p}\over
p\overline{p}+\epsilon ^2}[H_{{\rm
arg\,}p}]_{\Lambda _t},df_t\rangle \mu
(d\rho ). }
By the dominated convergence theorem, we got
\ekv{6.3}
{
\partial _tI_\epsilon (\Lambda _t,p)\to
-\int_{\Lambda _t}\langle [d({\rm
arg\,}p)]_{\Lambda _t},H_{f_t}\rangle \mu
(d\rho ) =\int_{\Lambda _t}\langle [H_{{\rm
arg\,}p}]_{\Lambda _t},df_t\rangle \mu
(d\rho ), }
when $\epsilon \to 0$, and since
$\Vert [d({\rm arg\,}p)]_{\Lambda_t }\Vert
_{L^1(\Lambda ,\mu )}$ is locally uniformly
bounded with respect to $t$, we concluded
that $t\mapsto I_\epsilon (\Lambda _t,p)$ is
uniformly Lipschitz with respect to $\epsilon
$, and in particular that $t\mapsto I(\Lambda
_t,p)$ is Lipschitz and that the a.e.
derivative of the last function is given by
\ekv{6.4}
{
{\rm "}\partial _tI(\Lambda
_t,p){\rm "}=-\int_{\Lambda _t}
\langle [d({\rm arg\,}p)]_{\Lambda
_t},H_{f_t}\rangle \mu (d\rho )
=\int_{\Lambda _t}\langle [H_{{\rm
arg\,}p}]_{\Lambda _t},df_t\rangle \mu
(d\rho ). }
In a more fancy way, we can consider an open
connected set $L$ in the family of smooth
IR-manifolds as in section 4, and define the
"distance"between
$\Lambda _0,\Lambda _1\in L$ to be
\ekv{6.5}
{{\rm dist\,}(\Lambda _0,\Lambda _1)=
\inf \int_0^1\Vert df_t\Vert_{L^\infty
(\Lambda _t)}
\mu (d\rho )dt, }
where the infimum is taken over all smooth
curves $[0,1]\ni t\mapsto \Lambda _t\in L$,
that link $\Lambda _0$ to $\Lambda _1$.
(We did not check that this is really a
distance, nor did we consider the problem of
studying the completion of $L$.) We can think
of $L$ as a manifold. The tangent space of $L$
at a point $\Lambda $ is then the space of all
$df$, where $f\in C_b^\infty (\Lambda ;{\bf
R})$. In the same spirit,
we have the differentials
\ekv{6.6}
{
\langle d_\Lambda I_\epsilon (\Lambda
,p),\delta \Lambda \rangle =\int _\Lambda
\langle {p\overline{p}\over
p\overline{p}+\epsilon ^2}[H_{{\rm arg\,
}p}]_\Lambda ,df\rangle \mu (d\rho ), }
\ekv{6.7}
{
\langle "d_\Lambda I(\Lambda
,p)",\delta \Lambda \rangle =\int_\Lambda
\langle [H_{{\rm arg \,}p}]_\Lambda
,df\rangle \mu (d\rho )=-\int_\Lambda
\langle [d\,{\rm arg\,}p]_\Lambda ,H_f\rangle
\mu (d\rho ), }
where $df$ is the differential on $\Lambda $
corresponding to the infinitessimal
variation $\delta \Lambda $. Here "$d_\Lambda
I(\Lambda ,p)$" becomes the almost everywhere
differential of our Lipschitz function
$I(\cdot ,p)$ whenever we restrict $\Lambda $
to vary in a curve or more generally in a
finite dimensional submanifold.

\par We next give sufficient conditions for
the continuity of "$d_\Lambda I(\cdot ,p)$"
at some given point $\Lambda $, and for
having a power law estimate in the
convergence $d_\Lambda I_\epsilon (\Lambda
,p)\to$ \break $ {\rm "}d_\Lambda I(\Lambda
,p){\rm "}$. We start with the
continuity question. Let
$\Lambda _0$ be a fixed IR-manifold (as in
section 4). Let
$x_0\in\Lambda _0$ be a point, where
$p(x_0)=0$ and choose local coordinates
$(x_1,..,x_{2n})$ centered at $x_0$, such
that for some $0<N_0\in{\bf N}$:
$\partial _{x_j}^{N_0}p(0)\ne 0$, $\partial
_{x_j}^k p(0)=0$, $0\le k\le N_0-1$, $1\le
j\le 2n$. Let
$p_\Lambda ={p_\vert}_{\Lambda }$. We are
interested in studying how $d\,{\rm
arg\,}p_\Lambda $ varies when we make a small
variation of
$\Lambda $. For simplicity, we concentrate on
$\partial _{x_1}{\rm arg\,}p_\Lambda $, and
we view $p_\Lambda $ as a small perturbation
in $C^\infty $ of $p_{\Lambda _0}$, by
choosing coordinates $(x_1,..,x_{2n})$ on
$\Lambda _t$ depending smoothly on $t$.

\par By Malgrange's preparation theorem:
\ekv{6.8}
{
p_{\Lambda
_0}(x)=q(x)(x_1^{N_0}+\sum_{1}^{N_0}a_j(x')
x_1^{N_0-j}),}
where $q,\, a_j$ are smooth with $q(0)\ne
0$, $a_j(0)=0$, and we write
$(x_1,x_2,..,x_{2n})=(x_1,x')$. Let $\lambda
_1 (x'),..,\lambda _{N_0}(x')$ be the roots of
the last factor of (6.8) (where the ordering
doesn't matter) so that locally, up to a
multiple of $2\pi $,
\ekv{6.9}
{
{\rm arg\,}p_{\Lambda _0}(x)={\rm
arg\,}(q(x))+\sum_{1}^{N_0}{\rm
arg\,}(x_1-\lambda _j(x')),\ \hbox{ for
}p_{\Lambda _0}(x)\ne 0. } From this we see
that
\ekv{6.10}
{
\int_{-\epsilon _0}^{\epsilon _0}\vert
\partial _{x_1} {\rm arg\,}p_{\Lambda
_0}(x)\vert dx_1\le C<\infty , }
for $\vert x'\vert \le \epsilon _0$, with
$\epsilon _0>0$ fixed and sufficiently small,
and where the integration is restricted to the
points where $p_{\Lambda _0}(x)\ne 0$. Assume,
\ekv{{\rm H}(x_1)}
{
\bigcup_{j=1}^{N_0}\{ x'\in{\bf
R}^{2n-1};\,\vert x'\vert \le \epsilon
_0,\,\lambda _j(x')\in{\bf R}\}\hbox{ is of
Lebesgue measure 0.} } (Equivalently, $\Pi
_{x'}(p_{\Lambda _0}^{-1}(0)\cap([-\epsilon
_0,\epsilon _0]\times B(0,\epsilon _0)))$
should be of measure $0$, where $\Pi _{x'}$
denotes the projection $x\mapsto x'$.) Let
$\Lambda _t$,
$t\in{\rm neigh\,}(0,{\bf R}^k)$ be a smooth
family of IR-manifolds. Then for $t$
sufficiently close to $0$, (6.8) extends:
\ekv{6.11}
{
p_{\Lambda _t}(x)=q(x,t)(x_1^{N_0}+
\sum_1^{N_0}a_j(x',t)x_1^{N_0-j}), }
(choosing the local coordinates to depend
smoothly on $t$) and we want to estimate
\ekv{6.12}
{
\int \vert \partial _{x_1}{\rm
arg\,}p_{\Lambda _0}(x)-\partial _{x_1}{\rm
arg\,}p_{\Lambda _t}(x)\vert dx, }
where we integrate over $\{ x;\,\vert
x_1\vert <\epsilon _0,\, \vert x'\vert
<\epsilon _0,\, p_{\Lambda _0}(x)\ne 0\ne
p_{\Lambda _t}(x)\}$.
The corresponding integral
\ekv{6.13}
{K_t(x')=\int \vert \partial _{x_1}{\rm
arg\,}p_{\Lambda _0}(x)-\partial _{x_1}{\rm
arg\,}p_{\Lambda _t}(x)\vert dx_1}
is a bounded function of $x'$, and converges
to $0$ when $t\to 0$ for $x'$ outside the set
in (H($x_1$)). So under that assumption, the
integral (6.12) converges to $0$
when $t\to 0$.

\par We next find an equivalent form
of (H($x_1$)), which is easier to formulate
globally on $\Lambda _0$: Introduce the
assumption
\eekv{{\rm A}}
{
\hbox{For every smooth hypersurface }\Gamma
\subset]-\epsilon _0,\epsilon _0[\times
B(0,\epsilon _0), } {p_{\Lambda
_0}^{-1}(0)\cap\Gamma
\hbox{ is of (Lebesgue) measure 0 in }\Gamma
.}
\medskip
\par\noindent
\bf Proposition 6.1. \it (A) and (H($x_1$))
are equivalent.
\medskip
\par\noindent
\bf Proof. \rm We first prove
(H($x_1$))$\Rightarrow$(A) and assume that
(A) does not hold. Let $\Gamma \subset
]-\epsilon _0,\epsilon _0[\times B(0,\epsilon
_0)$ be a smooth hypersurface such that
$p^{-1}(0)\cap \Gamma $ is of measure $>0$ in
$\Gamma $. The set
$$
(p^{-1}(0)\cap \Gamma )_\infty =\{ x\in\Gamma
\cap p^{-1}(0);\,{{\rm vol\,}(p^{-1}(0)\cap
\Gamma \cap B(x,\epsilon ))
\over {\rm vol\,}(\Gamma \cap B(x,\epsilon
))} \to 1,\,\epsilon \to 0\}
$$
is of full measure in $\Gamma \cap
p^{-1}(0)$ considered as a subset of
$\Gamma$. We notice that
${p_\vert}_{\Gamma }$
vanishes to infinite order on
$(p^{-1}(0)\cap\Gamma )_\infty $. The set
$$\eqalign{F_\infty =&\{ x\in (\Gamma \cap
p^{-1}(0))_\infty ;\cr &{{\rm vol\,}(\{ y\in
(\Gamma \cap p^{-1}(0))_\infty ;\,{(dx_2\wedge
..\wedge dx_{2n})_\vert }_\Gamma (y)=0\}\cap
B(x,\epsilon ))
\over {\rm vol\,}(B(x,\epsilon )\cap \Gamma
)}\to 1,\, \epsilon \to 0\} }$$
is of full measure in $\{x\in\Gamma \cap
p^{-1}(0);{(dx_2\wedge ..\wedge
dx_{2n})_\vert } _\Gamma (x) =0\}$.
At the points of $F_\infty $, the form
${(dx_2\wedge ..\wedge d_{2n})_\vert }_\Gamma $
vanishes to infinite order as well as
${p_\vert}_{\Gamma }$. Let
$y^0=(y_1^0,..,y_{2n}^0)\in F_\infty $. We
may assume that $\Gamma $ is of the form
$x_{2n}=a(x_1,..,x_{2n-1})$ near $y^0$. Then
${\partial a\over \partial x_1}(x)$ vanishes
to $\infty $ order at
$(y_1^0,..,y_{2n-1}^0)$, so at $y^0$ the
hypersurfaces $\Gamma $ and
$x_{2n}=a(y_1^0,x_2,..,x_{2n})$
are tangent to $\infty $ order. It follows
that
$p(x_1,x_2,..,x_{2n-1},\break
a(y_1^0,x_2,..,x_{2n-1}))$ vanishes to
$\infty $ order at $(y_1^0,..,y_{2n-1}^0)$ in
contradiction with the fact that $\partial
_{x_1}^{N_0}p(y^0)\ne 0$. Consequently
$F_\infty $ is empty and hence $$\{ x\in
\Gamma \cap p^{-1}(0);\,
{(dx_2\wedge ..\wedge
dx_{2n})_\vert}_\Gamma (x)=0\}$$ is of measure
0. In particular $\{ x\in (p^{-1}(0)\cap \Gamma
)_\infty ; {(dx_2\wedge ..\wedge
dx_{2n})_\vert}_{\Gamma }(x)\ne 0\}$ is of
non-vanishing measure, and the same holds for
its $x'$ space projection, which is in
contradiction with (H($x_1$)). We have proved
that (H($x_1$)) implies (A).
\par It remains to prove that (A) implies
(H$(x_1)$). Recall that $\partial
_{x_1}^{N_0}p\ne 0$ in $]-\epsilon
_0,\epsilon _0[\times B(0,\epsilon _0)$.
Without loss of generality, we may assume
that $\partial _{x_1}^{N_0}{\rm Re\,}p\ne 0$
there. Then
$$p^{-1}(0)\subset\cup_{N=1}^{N_0}\Gamma _N,$$
where
$$\Gamma _N=\{ x;\,\partial _{x_1}^N{\rm
Re\,}p(x)\ne 0,\, \partial _{x_1}^{N-1}{\rm
Re\,}p(x)=0\}$$
are smooth hypersurfaces such that
${{\Pi _{x'}}_\vert}_{\Gamma _j}$
is a local diffeomorphism for $j=1,..,N_0$.
Applying (A) to $p^{-1}(0)\cap \Gamma _j$, it
follows that $\Pi _{x'}(p^{-1}(0))$ is of
measure 0.\hfill{$\#$}
\medskip
\par Combining Proposition 6.1 with its
preceding discussion, we get
\medskip
\par\noindent \bf Theorem 6.2. \it Let
$\Lambda _0$ be an IR manifold and assume
\eekv{{\rm A}_{{\rm glob}}}
{\hbox{For every smooth hypersurface }\Gamma
\subset \Lambda _0,}
{p^{-1}(0)\cap\Gamma \hbox{ is of Lebesgue
measure 0 in } \Gamma .} Then  "$d_\Lambda
I(\cdot ,p)$" is continuous at $\Lambda
=\Lambda _0$ in the following sense: Let ${\rm
neigh\,}(0,{\bf R}^k)\ni t\mapsto \Lambda _t$
be a smooth family of IR manifolds as in
section 4, such that
$\Lambda _{t=0}=\Lambda _0$. Let
$f_1(t),..,f_k(t)$ be a corresponding system
of generating functions on $\Lambda _t$, so
that $f_\ell$ is generating for the one
parameter family obtained by varying $t_\ell$,
(freezing all the $t_j$ with $j\ne \ell$) and
hence that
\ekv{6.14}
{
{\rm "}\partial _{t_\ell}I(\Lambda _t,p){\rm
"}=-\int_{\Lambda _t}\langle [d\,{\rm
arg\,}p]_{\Lambda _t},H_{f_\ell (t)}\rangle \mu
(d\rho ). }
Then "$\partial _{t_\ell}I(\Lambda _t,p)$" are
continuous at $t=0$ and consequently
$t\mapsto I(\Lambda _t,p)$ is differentiable
at that point with partial derivatives
(without the " ") given by (6.14).
\rm\medskip

\par In practice, it might be preferable to
use Proposition 6.1, which says that $({\rm
A}_{{\rm glob}})$ holds iff for each $x_0\in
p^{-1}(0)$ there is a (or equivalently for
all) system(s) of smooth local coordinates
$x_1,..,x_{2n}$ centered at that point, such
that $\partial _{x_1}^{N_0}p(0)\ne 0$ for
some $N_0\in{\bf N}\setminus\{ 0\}$ (that we
choose to be minimal), and the roots $\lambda
_j(x')$ of the last factor in (6.8) have
imaginary parts that are
$\ne 0$ on a set of full measure in a
neighborhood of
$x'=0$.

\par We now turn to the question of having a
power law for the convergence $d_\Lambda
I_\epsilon (\Lambda _0,p)\break\to{\rm
"}d_\Lambda I(\Lambda _0 ,p){\rm "}$,
$\epsilon \to 0$. Here we did not find a
nice invariant condition and content ourselves
with
\eeeeekv{{\rm B}}
{
\hbox{For every $x_0\in p^{-1}(0)\cap
\Lambda _0$, we can find $2n$ linearly
independent vectorfields}
}
{
\nu _1,..,\nu _{2n}\hbox{ near $x_0$, such
that $\nu _j^kp(0)=0$, $0\le k\le N_j-1$,
$\nu _j^{N_j}p(0)\ne 0$, for}
}
{
\hbox{some $0<N_j\in{\bf N}$, and such
that if $\Gamma _j$ is a hypersurface passing
through $x_0$,}
}
{\hbox{transversal to $\nu
_j$, then there exist $\epsilon _0,\delta
_0>0$, such that for $0<\epsilon \le \epsilon
_0$,}  }
{{\rm vol\,}(\{ x'\in\Gamma
_j;\, {\rm dist\,}(x',x_0)\le \epsilon _0,\,
\inf_{\vert
t\vert <\epsilon _0}\vert p(\exp t\nu
_j(x'))\vert <\epsilon \} )\le \epsilon
^{\delta _0}.
}
\medskip
\par\noindent \bf Proposition 6.3. \it Under
the assumption (B), there exists $\delta
_1>0$ such that
\ekv{6.15}
{
[d\,{\rm arg\,}p]_{\Lambda
_0}-{p\overline{p}\over \epsilon
^2+p\overline{p}}[d\,{\rm arg\,}p]_{\Lambda
_0}={\cal O}(\epsilon ^{\delta _1}),\hbox{
in }L^1. }
\rm\medskip
\par\noindent \bf Proof. \rm The expression
in (6.15) is equal to
\ekv{6.16}
{
{\epsilon ^2\over p\overline{p}+\epsilon
^2}[d\,{\rm arg\,}p]_{\Lambda _0}. }
The $L^1$ norm of this function over any
region $\vert p\vert \ge {\rm Const.}>0$ is
${\cal O}(\epsilon ^2)$, so we only have to
examine what happens near a point $x_0$,
where $p$ vanishes. If $\nu _1,..,\nu _{2n}$
are the corresponding vectorfields appearing
in (B), it is enough to show that
\ekv{6.17}
{
\Vert {\epsilon ^2\over
p\overline{p}+\epsilon ^2}[\nu _j({\rm
arg\,}p)]_{\Lambda _0}\Vert _{L^1({\rm
neigh\,}(x_0,\Lambda _0))}={\cal O}(\epsilon
^{\delta _1}), } for every $j$, with the
obvious definition of $[\nu _j(\,{\rm
arg\,}p)]_{\Lambda _0}$. Fix a
$j$ and choose local coordinates
$x_1,..,x_{2n}$ centered at $x_0$, such that
$\nu _j=\partial _{x_1}$.

\par Independently of the proof we are
engaged in, it may be of interest to notice
that if we write (6.8) and let $\lambda
_1(x'),..,\lambda _{N_j}(x')$ be the
corresponding roots, (with $N_j=N_0$) then
the volume estimate in (B) implies that for
some new $\delta _0$:
\ekv{6.18}
{
{\rm vol\,}(\{ x'\in B(0,\epsilon _0);\,\vert
{\rm Im\,}\lambda _k(x')\vert \le \epsilon
\})\le{\cal O}(\epsilon ^{\delta _0}), }
for all $k$. Conversely we can go from
(6.18) to the volume estimate in (B) (with a
new $\delta _0$).

\par Rather than quoting (6.8) directly we
shall only use the fact that
\ekv{6.19}
{
\sup_{\vert x'\vert \le \epsilon
_0}\int_{-\epsilon _0}^{\epsilon _0}\vert
\partial _{x_1}{\rm arg\,}p(x)\vert
dx_1\le{\rm Const.} }
Let $0<\epsilon \le\widetilde{\epsilon }\ll
1$. Then, with $q(x')=\inf_{-\epsilon
_0<t<\epsilon _0}\vert p(t,x')\vert$,
$$\eqalign{
&\int_{B(0,\epsilon _0)}\int_{-\epsilon _0}^{
\epsilon _0}{\epsilon ^2\over
p\overline{p}+\epsilon ^2}\vert \partial
_{x_1}{\rm arg\,}p(x)\vert dx_1dx'=\cr
&\int_{q(x')\le\widetilde{\epsilon
}}\int \vert \partial _{x_1}{\rm
arg\,}p(x)\vert  dx_1dx' +{\cal O}(({\epsilon
\over
\widetilde{\epsilon }})^2)\le {\cal
O}(\widetilde{\epsilon }^{\delta
_0}+({\epsilon \over \widetilde{\epsilon
}})^2),
}$$
and choosing $\widetilde{\epsilon }=\epsilon
^\alpha $ with $\alpha =2/(2+\delta _0)$, we
get (6.17) with $\delta _1=2\delta
_0/(2+\delta _0)$.\hfill{$\#$}

\bigskip

\centerline{\bf 7. Minimality to infinite order
of critical points.}
\medskip
\par  Let $\Lambda
_0\subset {\bf C}^{2n}$ be an IR manifold
and $p(x,\xi )$ a holomorphic function as in
section 4. We assume that $\Lambda _0$
is a critical point for the functional
\ekv{7.1}
{I(\Lambda ,p)={1\over 2}\int_\Lambda  \log
(p\overline{p})\mu (d\rho ),}
in the sense that
\ekv{7.2}
{
\langle {\rm "}d_\Lambda I(\Lambda _0,p){\rm
"},\delta \Lambda \rangle =\int_{\Lambda _0}
\langle [H_{{\rm arg\,}p}]_{\Lambda
_0},df\rangle \mu (d\rho )=0, }
for $f\in C_b^\infty (\Lambda ;{\bf R})$
corresponding to the infinitesimal variation
$\Lambda _0+\delta \Lambda $. Equivalently,
\ekv{7.3}
{
{\rm div\,}[H_{{\rm arg\,}p}]_{\Lambda _0}=0,
}
where we recall that
${\rm div\,}[H_{{\rm arg\,}p}]_{\Lambda _0}$
is a distribution of order $\le 1$ on
$\Lambda _0$ with support in $p^{-1}(0)\cap
\Lambda _0$. We assume that $\Lambda _0$
is a regular point in the sense that there
exists $\delta _1>0$ such that
\ekv{7.4}
{
[d\,{\rm arg\,}p]_{\Lambda
_0}-{p\overline{p}\over \epsilon
^2+p\overline{p}}[d\,{\rm arg\,}p]_{\Lambda
_0}={\cal O}(\epsilon ^{\delta _1})\hbox{ in
}L^1, }
and recall that this implies that
\ekv{7.5}
{
\langle d_\Lambda I_\epsilon (\Lambda
_0,p),\delta \Lambda \rangle -\langle {\rm
"}d_\Lambda I(\Lambda _0,p){\rm "},\delta
\Lambda \rangle ={\cal O}(\epsilon ^{\delta
_1}), }
when the form $df$, corresponding to $\delta
\Lambda $ is bounded. We also recall that
(7.4) is a consequence of the property (B) of
the preceding section.
\medskip
\par\noindent \bf
Theorem 7.1. \it Assume (7.3), (7.4) and
let
$${\rm neigh\,}(0,{\bf R}^k)\ni t\mapsto
\Lambda _t$$
be a smooth family of IR manifolds as in
section 4, with
$\Lambda _{t=0}=\Lambda _0$. Then
\ekv{7.6}
{
I(\Lambda _t,p)-I(\Lambda _0,p)\ge -C_N\vert
t\vert ^N, }
for every $N\in {\bf N}$.
\rm\medskip
\par\noindent \bf Proof. \rm Since our
estimates will be uniform w.r.t. additional
parameters, we may assume that $k=1$. The
next result says that up to an error ${\cal
O}(t^\infty )$ it is possible to obtain
$\Lambda _t$ as the result of an autonomous
almost holomorphic flow acting on $\Lambda
_0$, with a real generator (on $\Lambda _0$).
\medskip
\par\noindent\bf Proposition 7.2. \it There
exists $f_t(\rho )=f(t,\rho )\in C_b^\infty
({\rm neigh\, }(0,{\bf R})\times {\bf
C}^{2n})$, real-valued on $\Lambda _0$ and
almost holomorphic in
$\rho $ on $\Lambda _0$, such that
\ekv{7.7}
{
{\rm dist\,}(\Lambda _t,\exp
(t\widehat{iH_{f_t}})(\Lambda _0))={\cal
O}(t^\infty ). }
Moreover $df_t$ is unique on $\Lambda _0$  mod
${\cal O}(t^\infty )$.\rm\medskip
\par\noindent \bf Proof. \rm We first
consider the situation locally near a point
$(x_0,\xi _0)\in\Lambda _0$. After a complex
canonical transformation, we may assume that
$\Lambda _t$ is given by
\ekv{7.8}
{
\xi ={2\over i}{\partial H_t\over \partial
x}(x), }
where $H_t(x)$ is real and smooth in
$(t,x)$ and strictly plurisubharmonic in $x$.
Let
$f\in C^\infty (\Lambda _0;{\bf R})$ and
extend $f$ to a smooth function on ${\bf
C}^{2n}$, almost holomorphic at $\Lambda _0$.
Then up to ${\cal O}(t^\infty )$,
$\widetilde{\Lambda }_t=\exp
(t\widehat{iH_f})(\Lambda _0)$ is given by
\ekv{7.9}
{
\xi ={2\over i}{\partial G_t\over \partial
x}(x), }
where $G_t(x)=G(t,x)$ solves the
Hamilton--Jacobi problem
\ekv{7.10}
{\partial _tG(t,x)-{\rm Re\,}f(x,{2\over
i}\partial _xG(t,x))=0,\ G(0,x)=H_0(x).}
(Here we use that up to infinite order at
$\Lambda _0$,we have $\widehat{iH_f}=H_{{\rm
Re\,}f }^{{\rm Im\,}\sigma }=H_{-{\rm
Re\,}f}^{-{\rm Im\,}\sigma }$.) Moreover, it
follows that
\ekv{7.11}
{
{\rm Re\,}f(x,{2\over i}\partial
_xG(t,x))=f(x,{2\over i}\partial
_xG(t,x))+{\cal O}(t^\infty ), }
so (7.10) gives
\ekv{7.12}
{
\partial _tG(t,x)-f(x,{2\over i}\partial
_xG(t,x))={\cal O}(t^\infty ),\
G(0,x)=H_0(x). } Notice that
\ekv{7.13}
{
G(t,x)=H_0(x)+tf(x,{2\over i}\partial
_xH_0(x))+{\cal O}(t^2). }

\par If we let $f=f_s$ depend smoothly on a
real parameter $s$, then differentiating
(7.12), we get ${\partial
_sG_\vert}_{t=0}=0$,
\ekv{7.14}
{\partial _t\partial _s G-{2\over i}\partial
_\xi f_s(x,{2\over i}\partial
_xG(t,x))\cdot \partial _x\partial _sG(t,x)
=(\partial _sf_s)(x,{2\over
i}\partial _xG(t,x))+{\cal O}(t^\infty ), }
and it follows that
\ekv{7.15}
{
\partial _sG=t\partial _sf_s(x,{2\over
i}\partial _xH_0)+{\cal O}(t^2). }
So if we replace $f(x,\xi )$ by $f(x,\xi
)+\epsilon k(x,\xi )$, then $G(t,x)$ is
replaced by
\ekv{7.16}
{
G(t,x)+t\epsilon k(x,{2\over i}\partial _x
H_0)+{\cal O}(t^2\epsilon )+{\cal
O}(t\epsilon ^2). }

\par It is now clear how to get the local
existence of a suitable $f_s$ in (7.7) by
successive approximations. We start by
choosing (cf. (7.13)) $f^{(2)}(x,\xi )$ so
that
\ekv{7.17}
{
G(t,x)=H(t,x)+{\cal O}(t^2).
}
Then (with $\epsilon =s$), we try a corrected
\ekv{7.18}
{f_s^{(3)} =f^{(2)}(x,\xi )+sk(x,\xi ).
}
Then (with $s=t$), $G(t,x)$ is replaced by
$$G(t,x)+t^2k(x,{2\over i}\partial
_xH_0(x))+{\cal O}(t^3),$$
so there is a unique choice of $k$ so that
for the new $G(t,x)$, we have
\ekv{7.19}
{
G(t,x)=H(t,x)+{\cal O}(t^3).
}
Now we take $\epsilon =s^2$ and try a
corrected
$$f_s^{(4)}(x,\xi )=f^{(3)}(x,\xi
)+s^2k(x,\xi ),$$
with a new $k$. Then $G(t,x)$ in (7.19) is
replaced by
\ekv{7.20}
{
G(t,x)+t^3k(x,{2\over i}\partial _xH_0(x))
+{\cal O}(t^4),
}
and there is a unique $k$ so that the new $G$
satisfies
\ekv{7.21}
{
G(t,x)=H(t,x)+{\cal O}(t^4).
}
Continuing this way we get the local
existence of $f_s$.

\par We next look at the local uniqueness.
Let $f_s(x,\xi )$, $g_s(x,\xi )$ be two
functions which are almost holomorphic in
$(x,\xi )$ at $\Lambda _0$, smooth in
$(s,x,\xi )$ and real-valued for $(x,\xi
)\in\Lambda _0$. Let $k\in{\bf N}$ be such
that $f_s(x,\xi )-g_s(x,\xi )=s^kh(x,\xi
)+{\cal O}(s^{k+1})$ for some smooth function
$h$. Then the two varieties
$\exp (t\widehat{iH_{f_t}})(\Lambda _0)$,
$\exp (t\widehat{iH_{g_t}})(\Lambda _0)$ are
given (mod. ${\cal O}(t^\infty )$) by $\xi
={2\over i}\partial _xH(t,x)$ and $\xi
={2\over i}\partial _xG(t,x)$, with
$$H_t(x)-G_t(x)=t^{k+1}h(x,{2\over i}\partial
_xH_0(x))+{\cal O}(t^{k+2}),$$
and if $H=G$, it follows that $h(x,{2\over
i}\partial _xH_0(x))=0$, and hence (since
$\Lambda _0$ is maximally totally real), that
$h={\cal O}({\rm dist\,}(\cdot ,\Lambda
_0)^\infty )$. (That we have uniqueness only
up to a $t$ dependent constant is due to the
fact that $H_t(x)$ is unique only up to such
a constant.)\par \hfill{$\#$}
\medskip

\par As already used in the proof above, we
may replace $\exp
(t\widehat{iH_{f_t}}(\Lambda _0))$ in the last
proposition by $\widetilde{\Lambda }_{t,t}$,
where $\widetilde{\Lambda }_{t,s}=\exp
s(H_{{\rm Re\,}f_t}^{{\rm Im\,}\sigma
})(\Lambda _0 )$. According to (5.18), we
have
\eekv{7.22}
{
\partial _s^2\int_{\widetilde{\Lambda
}_{t,s}}F(p)\mu (d\rho
)=4\int_{\widetilde{\Lambda }_{t,s}}{\partial
F\over \partial \overline{p}\partial p}(p)
H_{\widetilde{f}_{t,s}}(p)
H_{\widetilde{f}_{t,s}}(\overline{p})\mu
(d\rho )}
{\hskip 4cm +i\int_{\widetilde{\Lambda }_{t,s}}
[H_p({\partial F\over \partial
p})-H_{\overline{p}}({\partial F\over
\partial \overline{p}})]\partial
_s\widetilde{f}_{t,s}\mu (d\rho ), }
where $\widetilde{f}_{t,s}$ denotes an almost
holomorphic extension of ${\rm Re\,}f_t$ from
$\widetilde{\Lambda }_{t,s}$ and all the
Hamilton fields are in the sense of the
symplectic \mfld{}
$\widetilde{\Lambda}_{t,s}$. Now
$f_t$ is almost holomorphic on $\Lambda _0$ and
$f_t={\rm Re\,}f_t+{\cal O}(s^\infty )$ on
$\widetilde{\Lambda }_{t,s}$. Consequently,
$\nabla _{s,\rho }(\widetilde{f}_{t,s}(\rho
)-f_t(\rho ))={\cal O}(s^\infty )$ on
$\widetilde{\Lambda }_{t,s}$ and in particular
\ekv{7.23}
{
\partial _s\widetilde{f}_{t,s}={\cal
O}(s^\infty )\hbox{ on }\widetilde{\Lambda
}_{t,s}. }

\par We shall use this in (7.22) with
$F(\rho )={1\over 2}\log ({\epsilon
^2+p\overline{p}\over \epsilon ^2+1})$. Since
${\partial ^2F\over \partial
\overline{p}\partial p}\ge 0$ (as we shall
review in the beginning of section 8), the
first term to the right will be
$\ge 0$. The second derivatives of $F$ are
${\cal O}(\epsilon ^{-2})$ on some compact
subset of
$\widetilde{\Lambda }_{t,s}$ and uniformly
\bdd{} outside, so the last term in (7.22) is
${\cal O}(s^M/\epsilon ^2)$ for every
$M\in{\bf N}$. It follows that
\ekv{7.24}
{I_\epsilon (\widetilde{\Lambda
}_{t,t},p)-I_\epsilon (\Lambda _0,p)\ge
t(\partial _s)_{s=0}I_\epsilon
(\widetilde{\Lambda }_{t,s},p)-{\cal
O}_M(1){t^M\over \epsilon ^2}.}
Here
\eekv{7.25}
{
(\partial _s)_{s=0}I_\epsilon
(\widetilde{\Lambda }_{t,s},p)=(\partial
_s)_{s=0}(I_\epsilon (\widetilde{\Lambda
}_{t,s},p)-I_0(\widetilde{\Lambda }_{t,s},p)) }
{
=-\int_{\Lambda _0}\langle [d\,{\rm
arg\,}p]_{\Lambda _0}-{p\overline{p}\over
\epsilon ^2+p\overline{p}}[d\,{\rm
arg\,}p]_{\Lambda
_0},H_{\widetilde{f}_{t,0}}\rangle \mu
(d\rho )={\cal O}(\epsilon ^{\delta _1}), }
according to (7.3), (7.4),
 and using this in (7.24) together with the
fact that $I_\epsilon (\widetilde{\Lambda
}_{t,t},p)-I_\epsilon (\Lambda _t,p)={\cal
O}(t^\infty )$ (since ${\rm
dist\,}(\widetilde{\Lambda }_{t,t},\Lambda
_t)={\cal O}(t^\infty )$ and $\Lambda \mapsto
I_\epsilon (\Lambda ,p)$ is uniformly
Lipschitz), we get
\ekv{7.26}
{
I_\epsilon (\Lambda _t,p)-I_\epsilon (\Lambda
_0,p)\ge -{\cal O}(\epsilon ^{\delta
_1})\vert t\vert -{\cal O}_M(1){t^M\over
\epsilon ^2}. }
Finally we use that $I_\epsilon (\Lambda
,p)-I(\Lambda ,p)={\cal O}(\epsilon ^{1/m})$
(see Lemma 4.2), to get
\ekv{7.27}
{
I(\Lambda _t,p)-I(\Lambda _0,p)\ge -{\cal
O}_M( \epsilon ^{1/m}+\epsilon ^{\delta
_1}+{t^M\over \epsilon ^2}). }
Choosing first $\epsilon =t^K$ with $K$ very
large, and then $M$
sufficiently large, we obtain (7.6) with $N$
as large as we like, and the proof of Theorem
7.1 is complete.\hfill{$\#$}\medskip

\par Assume that $p(\rho )=p(z,\rho )$
depends holomorphically on $z\in\Omega
\subset{\bf C}$. It is clear that $I(\Lambda
,p(z))$ is subharmonic in $z$ for every fixed
$\Lambda $. We end this section by giving a
completely formal argument which indicates
that if we have a minimizer $\Lambda
=\Lambda (z)$ to $\Lambda \mapsto I(\Lambda
,p(z))$, then $I(\Lambda (z),p(z))$ will also
be subharmonic. The argument is formal but it
should be possible to turn it into a proof
whenever we have a sufficiently good control
over the variational problem:

For $\Lambda $ close to $\Lambda (0)$, we
have
$$\Lambda =\Lambda _f=\exp
(\widehat{iH_f)}(\Lambda (0)),$$
to infinite order at $\Lambda (0)$, where
$f\in C_b^\infty (\Lambda (0);{\bf R})$ is
small and we also let $f$ denote an almost
holomorphic extension. If we allow
$f$ to be complex valued, we still get an
IR-\mfld , close to
$\Lambda (0)$, and we see that $I(\Lambda
_f,p(z))=I(f,z)$ is plurisubharmonic in
$(z,f)$. At $z=0$, $f=0$, the Hessian (that
we assume exists) of $I(f,z)$ is a quadratic
form in $(\Re f,z)$ and the
plurisubharmonicity means that
$$\pmatrix{
{1\over 4}I''_{ff} &{1\over 2}I''_{fz}\cr
{1\over 2}I''_{\overline{z}f}
&I''_{\overline{z}z}}\ge 0.$$ Here
the subscript $f$ indicates ordinary (real)
derivatives with respect to $f$, while
$z$ and $\overline{z}$ indicate holomorphic
and anti-holomorphic derivatives (as in the
Levi form). In other terms:
\ekv{7.28}
{
0\le {1\over 4}(I''_{ff}\phi \vert \phi )+
{1\over 2}(I''_{fz}\zeta \vert \phi )+
{1\over 2}(I''_{\overline{z}f}\phi \vert
\zeta )+ (I''_{\overline{z}z}\zeta \vert \zeta
),  } for all tangent vectors $(\phi ,\zeta
)$, where we use standard sesquilinear scalar
products.

\par  Assume (for
instance after an arbitrarily small
convexification in
$f$) that $I''_{ff}>0$. For $z$ close to $0$,
let
$f(z)$ be the real function with $\Lambda
(z)=\Lambda _{f(z)}$, so that
$I'_f(f(z),z)=0$. Differentiating, we get
$0=I''_{ff}f'_z+I''_{fz}$,
$f'_z=-(I''_{ff})^{-1}I''_{fz}$. On the
other hand,
$(I(f(z),z))'_{\overline{z}}=
I'_{\overline{z}}(f(z),z)$, so that
$$(I(f(z),z))_{\overline{z}z}''=
I_{\overline{z}z}''(f(z),z)+
I''_{\overline{z}f}(f(z),z)f_z'(z)=
I''_{\overline{z}z}-I''_{ \overline{z} f}(
I''_{ff})^{-1}I''_{fz}.$$
Choosing $\phi
=-2(I''_{ff})^{-1}I''_{fz}\zeta $ in (7.28),
we see that this quantity is $\ge 0
$.
\bigskip

\centerline{\bf 8. The codimension 2 case.}
\medskip
\par We start by recalling a very
classical formula and its proof. Let
$\epsilon >0$. Then for $p\in {\bf C}$:
\ekv{8.1}
{
\partial _p\partial _{\overline{p}}\log
(\epsilon ^2+p\overline{p})=\partial _p
{p\over \epsilon ^2+p\overline{p}}={1\over
\epsilon
^2+p\overline{p}}-{p\overline{p}\over
(\epsilon ^2+p\overline{p})^2}={\epsilon
^2\over (\epsilon ^2+p\overline{p})^2}. }
This is a non-negative function which
tends to 0 for $z\ne 0$. If $L(dp)$
denotes the Lebesgue measure on ${\bf
C}$, we have
\eekv{8.2}
{
\int_{\bf C} {\epsilon ^2\over (\epsilon
^2+p\overline{p})^2}L(dp)=2\pi
\int_0^\infty {\epsilon ^2\over (\epsilon
^2+r^2)^2}rdr= }
{\pi \int_0^\infty {\epsilon ^2\over
(\epsilon ^2+t)^2}dt=\pi [-{\epsilon
^2\over (\epsilon ^2+t)}]_{t=0}^\infty
=\pi .}
We conclude that $\log (\epsilon
^2+p\overline{p})$ is subharmonic for
$0\le \epsilon \le 1$, and that
\ekv{8.3}
{\partial _p\partial _{\overline{p}}\log
(p\overline{p})=\pi \delta .}

\par Let $\Lambda \subset{\bf C}^{2n}$ be
an IR \mfld , and let $p$ be
a holomorphic function as in section 4.
Recall that
\ekv{8.4}
{I_\epsilon (\Lambda ,p)={1\over
2}\int_\Lambda \log({\epsilon
^2+p\overline{p}\over \epsilon ^2+1})\mu
(d\rho ),\ I_0=I.}

\par If $p=p(\rho ,z)$ depends
holomorphically on a complex parameter
$z$, we get for $\epsilon >0$:
\eekv{8.5}
{
\partial _z\partial
_{\overline{z}}I_\epsilon (\Lambda
,p)={1\over 2}\int_\Lambda \partial
_p\partial _{\overline{p}}(\log (\epsilon
^2+p\overline{p}))\,\partial _z
p\,\overline{\partial _z p}\,\mu (d\rho )
}  {
={1\over 2}\int_\Lambda {\epsilon ^2\over
(\epsilon ^2+p\overline{p})^2}\partial
_zp\overline{\partial _zp}\,\mu (d\rho ), }
so $I_\epsilon (\Lambda ,p)$ is a
subharmonic function of $z$. This also
holds for the limiting case $\epsilon =0$.

\par We shall investigate the effect of
small variations of $\Lambda $. Let
$\Lambda =\Lambda _0$, where $]-\delta
,\delta [\ni t\mapsto \Lambda _t$ is a
smooth deformation as in section 4. Let
$\widehat{H_{if_t}}$
be the corresponding deformation field,
$f_t\in C_b^\infty (\Lambda _t;{\bf
R})$ and $f_t$ also denotes a corresponding
almost holomorphic extension in $\rho $.
Let $\kappa _t:\Lambda _0\to \Lambda _t$ be
the corresponding flow, so that $\kappa _t$
is symplectic. Repeating earlier
calculations,
\ekv{8.6}
{I_\epsilon (\Lambda _t,p)=I_\epsilon
(\Lambda _0,p_t),\hbox{ where }p_t=p\circ
\kappa_t.}
Since
\ekv{8.7}
{
\partial _tp_t=iH_{f_t}(p)\circ \kappa _t,
}
we get
\eeekv{8.8}
{
\partial _tI_\epsilon (\Lambda
_t,p)=}{{1\over 2}\int_{\Lambda
_0}\partial _p(\log (\epsilon
^2+p_t\overline{p}_t))\,iH_{f_t}(p)\circ
\kappa _t\,\mu (d\rho )+{1\over
2}\int_{\Lambda _0}\partial
_{\overline{p}}(\log (\epsilon
^2+p_t\overline{p}_t))\,\overline{iH_{f_t}p}
\circ
\kappa _t\,\mu (d\rho ) } {=\Re
\int_{\Lambda _0}\partial _p(\log (\epsilon
^2+p_t\overline{p}_t))\, iH_{f_t}(p)\circ
\kappa_t\, \mu (d\rho ).}
This can also be written
\ekv{8.9}
{\partial _tI_\epsilon (\Lambda
_t,p)=\Re\int_{\Lambda
_t}\partial _p(\log (\epsilon
^2+p\overline {p}))\, iH_{f_t}(p)\,\mu
(d\rho ).}
Take $t=0$ for simplicity, and use that
$H_{f_t}(p)=-H_p(f_t)$, to get
\ekv{8.10}
{
(\partial _t)_{t=0}I_\epsilon (\Lambda
_t,p)=-\Re i\int_\Lambda
\partial _p(\log (\epsilon
^2+p\overline{p}))H_p(f_0)\mu (d\rho ),}
where $H_p$ is viewed as a differential
operator on $\Lambda $, namely the
Hamilton field of the restriction of $p$
to $\Lambda $ with respect to the
restriction of $\sigma $ to $\Lambda $.
In order not to make the notations too
heavy, we shall often not distinguish
between functions on ${\bf C}^{2n}$ and
their restrictions to $\Lambda $. Since
the transpose of $H_p$ with respect to
the symplectic volume form $\mu (d\rho )$
is equal to $-H_p$, we get
\eekv{8.11}
{(\partial _t)_{t=0}I_\epsilon (\Lambda
_t,p)=\Re i\int_\Lambda
f_0H_p(\partial _p(\log (\epsilon
^2+p\overline{p})))\, \mu (d\rho )}
{=\Re i\int_\Lambda
f_0(\partial _{\overline{p}}\partial
_p(\log (\epsilon ^2+p\overline{p})))\{
p,\overline{p}\}\mu (d\rho
)=2\int_\Lambda  f_0{\epsilon ^2\over
(\epsilon ^2+p\overline{}p)^2}{i\over 2}\{
p,\overline{p}\} \mu (d\rho ).}
Here $\{ p,\overline{p}\}=\{ p_\Lambda
,\overline{p_\Lambda }\}$ denotes the
Poisson bracket for the symplectic
\mfld{} $\Lambda $.
\par Throughout the remainder of this
section, we shall assume with
$p=p_\Lambda $:
\ekv{{\rm H}}
{dp ,\, d\overline{p}
\hbox{ are independent at all points of
}\Lambda
\cap p^{-1}(0).}
This implies that $\Sigma :=\Lambda \cap
p^{-1}(0)$ is a smooth (possibly empty)
submanifold of codimension 2 in $\Lambda
$. Since
$\Lambda $ is symplectic, it is
orientable: We say that the volume form
$\sigma ^n/n!$ is positive on $\Lambda $
(and hence identified with the volume
density $\mu (d\rho )$). We then have a
natural induced orientation on $\Sigma $
depending also on $p$: We say that a
$2n-2$ form $\alpha $ on $\Sigma $ is
positive at some point in $\Sigma $ if
$\alpha \wedge {i\over 2}dp\wedge
d\overline{p}$ is a positive multiple of
$\sigma ^n/n!$ at that point. (Here we
notice that ${i\over 2}dp\wedge
d\overline{p}=d\Re p\wedge d\Im p$.) We
define the Liouville measure on $\Sigma $
to be the density alias positive $2n-2$
form
$\lambda _{p,0}$ on $\Sigma $, such that
\ekv{8.12}{\lambda _{p,0}\wedge{i\over
2}dp\wedge d\overline{p}={\sigma ^n\over
n!}.}
In the next lemma we continue to write
$\{ p,\overline{p}\}=\{ p_\Lambda
,\overline{p_\Lambda } \}$.
\medskip
\par\noindent \bf Lemma 8.1. \it For
$n\ge 2$, we have
\ekv{8.13}
{{i\over 2}\{ p,\overline{p}\} \lambda
_{p,0}={({\sigma ^{n-1}\over
(n-1)!})_\vert}_{\Sigma }.}
In particular, it follows that $\int_\Gamma
{i\over 2}\{ p,\overline{p}\} \lambda
_{p,0}(d\rho )=0$, for every connected
component $\Gamma
$ of
$\Sigma $.\rm\medskip
\par\noindent \bf Proof. \rm We know that
$\Sigma $ is symplectic at a point $\rho
\in\Sigma $ (i.e. ${\sigma
_\vert}_{\Sigma }$ is non-degenrate at
$\rho $) iff $\{ p,\overline{p}\}\ne 0$,
and we know that ${{\sigma
^{n-1}}_\vert}_{\Sigma }$ is non-vanishing
precisely at the points where $\Sigma $
is symplectic. Consequently both members
of (8.13) vanish at the points where
$\Sigma $ is not symplectic. On the other
hand, near a point where $\Sigma $ is
symplectic, we can choose symplectic
coordinates $x_1,..,x_n,\xi _1,..,\xi _n$
on $\Lambda $, so that $\Sigma $ is given
by $x_n=\xi _n=0$. Since $\lambda _{p,0}$
is a form of maximal degree in $(x',\xi
')=(x_1,..,x_{n-1},\xi _1,..,\xi
_{n-1})$, we get from (8.12):
$$\eqalign{{\sigma ^n\over n!}=\lambda
_{p,0}\wedge {i\over 2}dp\wedge
d\overline{p}&=\lambda _{p,0}\wedge
{i\over 2}({\partial p\over \partial \xi
_n}d\xi _n+{\partial p\over \partial
x_n}dx_n)\wedge ({\partial
\overline{p}\over \partial \xi _n}d\xi
_n+{\partial \overline{p}\over \partial
x_n}dx_n)\cr &={i\over 2}\{
p,\overline{p}\}
\lambda _{p,0}\wedge d\xi _n\wedge
dx_n.}$$
Now use that $${\sigma ^n\over n!}=d\xi
_1\wedge dx_1\wedge d\xi _2\wedge
dx_2\wedge ..\wedge d\xi _n\wedge dx_n,$$
to conclude that
$$\hskip 1cm {i\over 2}\{ p,\overline{p}\}
\lambda _{p,0}=d\xi _1\wedge dx_1\wedge
..\wedge d\xi _{n-1}\wedge
dx_{n-1}={({\sigma  ^{n-1}\over
(n-1)!})_\vert}_{\Sigma }.\hskip 2cm \#$$
\medskip
\par Next we consider the limits of (8.5),
(8.11), when $\epsilon \to 0$. It is easy
to see that the contribution from the
region where $\vert p\vert \ge \delta $
is $0$, for every fixed $\delta >0$. On
the other hand, in the region $\vert
p\vert <\delta $, we introduce $\Sigma
_w:=p^{-1}(w)$, for $w\in D(0,\delta
):=\{\omega \in{\bf C};\,\vert \omega \vert
<\delta\} $, and the corresponding
Liouville form $\lambda _{p,w}$. Then
$$\lambda _{p,w}\wedge {i\over
2}dp\wedge d\overline{p}={1\over
n!}\sigma ^n,$$
at the points of $\Sigma _w$, and in view
of (8.5):
\eekv{*}
{
{1\over 2}\int_{\Lambda \cap
p^{-1}(D(0,\delta ))}{\epsilon ^2\over
(\epsilon ^2+p\overline{p})^2}\partial
_zp\overline{\partial _zp}\,\mu (d\rho
)=}{{1\over 2}\int_{D(0,\delta )} {\epsilon
^2\over (\epsilon
^2+w\overline{w})^2}(\int_{\Sigma
_w}\partial _zp\overline{\partial
_zp}\,\lambda _{p,w}(d\rho ))\wedge
d\Re w\wedge d\Im w .}
As we saw in the beginning of this
section,
$${\epsilon ^2\over (\epsilon
^2+w\overline{w})^2}\to \pi \delta
_{w=0},$$
so the expression (*) converges to
\ekv{8.14}
{
{\pi \over 2}\int_\Sigma \partial
_zp\overline{\partial _zp}\lambda
_{p,0}(d\rho ). }
We conclude that
\ekv{8.15}
{
\partial _z\partial
_{\overline{z}}I(\Lambda ,p)={\pi \over
2}\int_{\Lambda \cap p^{-1}(0)}\partial
_zp\overline{\partial _zp}\lambda
_{p,0}(d\rho ), }
first in the sense of distributions, then
in the classical sense, since the RHS in
this last equation is smooth.

\par The discussion also applies to
(8.11), and shows that
\ekv{8.16}
{
{(\partial _t)_\vert}_{t=0}I(\Lambda
_t,p)=2\pi \int_{\Lambda \cap
p^{-1}(0)}f_0 {i\over 2}\{
p,\overline{p}\} \lambda _{p,0}(d\rho )
=2\pi \int_{\Lambda \cap
p^{-1}(0)}f_0{1\over (n-1)!}\sigma ^{n-1}.}
We conclude that $\Lambda =\Lambda _0$ is
a critical point for the functional
$\Lambda \mapsto I(\Lambda ,p)$ iff for
$p={p_\vert}_{\Lambda }$, we have on
$\Lambda $:
\ekv{8.17}{p(\rho )=0\Rightarrow
\{p,\overline{p}\}=0.}
Combining (8.1) and the earlier
discussion with (5.18) and the subsequent
remark, we get
$$
\partial _t^2I(\Lambda _t,p)=2\pi
\int_{\Lambda _t\cap
p^{-1}(0)}H_pf_t\overline{H_pf_t}\lambda
_{p,0}(d\rho )+2\pi \int_{\Lambda _t\cap
p^{-1}(0)}(\partial _tf_t){i\over 2}\{
p,\overline{p}\} \lambda _{p,0}(d\rho ).
$$
If we extend the definition of $I(\Lambda
_t,p)$ to complex $t$ by almost
holomorphic extension of the flow $\kappa
_t$, a simpler and more direct
computation shows that
\ekv{8.18}
{
\partial _{\overline{t}}\partial
_tI(\Lambda _t,p)={\pi \over
2}\int_{\Lambda _t\cap
p^{-1}(0)}H_pf_t\overline{H_pf_t}\lambda
_{p,0}(d\rho ),\hbox{ for }t\hbox{ real}. }
In this last identity, we let $\partial
_t$ and $\partial _{\overline{t}}$ denote
holomorphic and antiholomorphic
derivatives.

\par Let $\Lambda $ be
critical for $I(\cdot ,p)$ in the sense
that we have (8.17). We shall next see
that if $n=2$ and if we make an
infinitesimal change $p\mapsto p+\delta
p$, then there is a corresponding
infinitesimal change $\Lambda \mapsto
\exp (iH_{\delta f})(\Lambda
)=\widetilde{\Lambda }$ of $\Lambda $,
such that $\widetilde{\Lambda }$ is
critical for $I(\cdot ,p+\delta p)$:
\medskip
\par\noindent \bf Theorem 8.2. \it Let
$n=2$, let $\Lambda $ be an IR-\mfld{}
and $p$ a holomorphic function as in
section 4. Assume (H) and (8.17). Let
$q$ be holomorphic and ${\cal O}(\langle
(x,\xi )\rangle ^m)$, $m<-2n$, in some
tubular \neigh{} of $\Lambda $ and put
$p_z(\rho )=p(\rho )+zq(\rho )$,
$z\in{\rm neigh\,}(0,{\bf R})$. Then
there exists $f\in C_b^\infty (\Lambda
;{\bf R})$ such that if $\Lambda _z=\exp
(z\widehat{iH_f})(\Lambda )$, is the
corresponding IR-deformation, we have
$$\{ p_z,\overline{p}_z\}={\cal
O}(z^2)\hbox{ on }p_z^{-1}(0)\cap \Lambda
_z.$$
Here the Poisson bracket is the one
given by the symplectic form on $\Lambda
_z$.\rm\medskip

\par In a forthcoming paper, we shall
show, by using non-linear
$\overline{\partial }$ equations, a
stronger version of the above result,
namely that the ${\cal O}(z^2)$ can be
replaced by $0$, for small $z$. We think
that the proof below has some
independent interest and reveals some
intereresting structures. It is based on
the use of second order elliptic
operators rather than $\overline{\partial
}$ type operators.

\par In proving the theorem, it will be
convenient to use the terminology of
infinitesimal variations $\delta p=zq$,
$z\to 0$, and $\delta f=zf$, $z\to 0$, so
that our calculations will be modulo
errors ${\cal O}(z^2)$.

\par From (8.17) and (H), we see that
\ekv{8.19}
{
\{
p,\overline{p}\}=\overline{a}p-a
\overline{p}, }
for some smooth function $a$ which is
uniquely determined on $\Sigma $ (and
where the Poisson bracket is the one of
the symplectic \mfld{} $\Lambda $). For a
general infinitesimal change
$\delta p$, of $p$, we get
\eekv{8.20}
{
\{ p+\delta p,\overline{p+\delta
p}\}=\overline{a}p-a\overline{p}+\{\delta
p,\overline{p}\}+\{p,\delta
\overline{p}\} =}{\overline{a}(p+\delta
p)-a(\overline{p}+\delta
\overline{p})+(H_p+a)(\delta
\overline{p})-(H_{\overline{p}}+
\overline{a})(\delta p), }
so (on $\Lambda $)
\ekv{8.21}
{\{ p+\delta p,\overline{p+\delta p} \}
=(H_p +a)(\delta \overline{p})
-(H_{\overline{p}}+\overline{a})(\delta
p)\hbox{ on }(p+\delta p)^{-1}(0).}
If we apply an infinitesimal change
$\Lambda \mapsto\widetilde{\Lambda }= \exp
iH_{\delta f}(\Lambda )$, we get
$${(p+\delta p)_\vert
}_{\widetilde{\Lambda }}\simeq p+\delta
p+{iH_{\delta f}(p)_\vert }_\Lambda
={(p+\delta p-iH_p(\delta f))_\vert
}_{\Lambda },$$
so we get
\eeekv{8.22}
{
{\{ p+\delta p,\overline{p+\delta p}\}
_\vert }_{\widetilde{\Lambda }\cap
(p+\delta p)^{-1}(0)}\simeq }{{ \{ p+\delta
p-iH_p(\delta f),\overline{p+\delta
p-iH_p(\delta f)}\}_\vert }_{\Lambda
\cap (p+\delta p-iH_p(\delta f))^{-1}(0)}
=}
{
(H_p+a)(\delta
\overline{p}+iH_{\overline{p}}\delta
f)-(H_{\overline{p}}+\overline{a})(\delta
p-iH_p\delta f), }
where we used (8.21) with a new
infinitesimal change in the last
step. Since the last expression is already
infinitely small, we can ignore the effect
of the infinitesimal displacement of the
zero set of $p$ and consider that it is
given on
$p^{-1}(0)$. To have it equal to 0 amounts
to solving on
$p^{-1}(0)$ the equation
$$i((H_p+a)H_{\overline{p}}+
(H_{\overline{p}}+\overline{a}
)H_p)\delta f=-(H_p+a)\delta
\overline{p}+(H_{\overline{p}}+\overline{
a})\delta p,$$
or in other terms
\ekv{8.23}
{
\Re (H_p+a)H_{\overline{p}}\delta f= -\Im
(H_p+a)\delta \overline{p}\hbox{ on
}p^{-1}(0)\cap\Lambda . }
Here $(H_p+a)H_{\overline{p}}$ is an
elliptic second order operator on the
leaves of the bicharacteristic foliation
of $\Sigma =p^{-1}(0)\cap \Lambda $, so
when $n=2$, it is an elliptic operator on
$\Sigma $. Let us notice that it is real
on $\Sigma $:
$$2i\Im
(H_p+a)H_{\overline{p}}=(H_p+a)H_{
\overline{p}}-(H_{\overline{p}}+
\overline{a})H_p=[H_p,
H_{\overline{p}}]-(\overline{a}H_p-aH_{
\overline{p}})=0,$$
where in the last step we used (8.19) and
the fact that
$[H_p,H_{\overline{p}}]=H_{\{ p,
\overline{p}\} }$. Since $\delta f$ is
real the problem (8.23) reduces to
\ekv{8.24}
{
-(H_p+a)H_{\overline{p}}(\delta f)=\Im
(H_p+a)\delta \overline{p}\hbox{ on
}\Sigma . }
\medskip
\par\noindent \bf Lemma 8.3. \it For every
smooth function $u$
on $\Sigma $ and every component  $\Gamma
$ of $\Sigma $, we have
\ekv{8.25}
{
\int_\Gamma (H_p+a)u\, \lambda _{p,0}(d\rho
)=0. }
In particular, if we replace $u$ by a
product $u\overline{v}$, we see that
\ekv{8.26}
{(H_p+a)^*=-H_{\overline{p}},}
where the star indicates that we take
the complex adjoint in
$L^2(\Gamma ,\lambda _{p,0}(d\rho
))$.\medskip
\par\noindent\bf Proof. \rm We compute
the Lie derivative ${\cal
L}_{H_p}(\lambda _{p,0})$: First we
recall that ${\cal L}_{H_p}\sigma =0$
and hence ${\cal L}_{H_p}\sigma ^n=0$.
Secondly, at the points of $\Sigma $, we
get
$${\cal L}_{H_p}(dp\wedge
d\overline{p})=dp\wedge d\{
p,\overline{p}\} = dp\wedge
(\overline{a}dp
-ad\overline{p})=-adp\wedge
d\overline{p}.$$
Still at the points of $\Sigma $, we
apply ${\cal L}_{H_p}$ to the relation
$$\lambda _{p,0}\wedge {i\over
2}dp\wedge d\overline{p}={1\over
n!}\sigma ^n,$$
and get
$$0={\cal L}_{H_p}(\lambda
_{p,0})\wedge {i\over 2}dp\wedge
d\overline{p}+\lambda _{p,0}\wedge
{i\over 2}{\cal L}_{H_p}(dp\wedge
d\overline{p})={\cal L}_{H_p}(\lambda
_{p,0})\wedge {i\over 2}dp\wedge
d\overline{p}-a\lambda _{p,0}\wedge{i\over
2}dp\wedge d\overline{p}.$$
It follows that
\ekv{8.27}
{{\cal L}_{H_p}(\lambda _{p,0})=a\lambda
_{p,0}.}
To get (8.25), we write
$$0=\int_\Gamma {\cal L}_{H_p}(u\lambda
_{p,0}(d\rho ))=\int_\Gamma
(H_p(u)\lambda _{p,0}+u{\cal
L}_{H_p}\lambda _{p,0})=\int_\Gamma
(H_p+a)u\,\lambda _{p,0}(d\rho ).$$
\hfill{$\#$}
\medskip

\par The operator in the LHS of (8.24) can
be written in different forms:
\ekv{8.28}
{
-(H_p+a)H_{\overline{p}}=(H_{\overline{p}}
)^*H_{\overline{p}}=H_p^*H_p. }
Here the first equality follows from
(26) and second one from the fact that
the operator is real. In particular, our
operator is formally self-adjoint. If
$u\in C^\infty (\Gamma )$, where $\Gamma
$ is a component of $\Sigma $, and
$(H_p+a)H_{\overline{p}}u=0$, we get
$$\Vert H_{\overline{p}}u\Vert
^2_{L^2(\Gamma ,\lambda _{p,0})}=\Vert
H_pu\Vert ^2_{L^2(\Gamma ,\lambda
_{p,0})}=0,$$
So $H_{\Re p}u=H_{\Im p}u=0$. This means
that $u$ is constant on every
bicharacteristic leaf and hence on the
whole component $\Gamma $, since $n=2$. Let
us now check that the RHS in (8.24) is
orthogonal to
$C^\infty (\Gamma )\cap {\rm
Ker\,}(-(H_p+a)H_{\overline{p}})$ also
for general $n$: The latter space is the
complexification of its maximal
subspace of real functions, and for
such a real function $u$, we get $$(\Im
(H_p+a)\delta
\overline{p}
\vert u)=\Im ((H_p+a) \delta \overline{p}
\vert u)=-\Im (\delta \overline{p} \vert
H_{
\overline{p}}u)=0.$$
Returning to the case $n=2$, our
operator is elliptic and
essentially self-adjoint. The previous
discussion implies that for every $\delta
p\in C^\infty (\Gamma )$, the equation
(8.24) has a solution $\delta f\in C^\infty
(\Gamma )$, which is unique up to a
constant. Theorem 8.2
follows.\hfill{$\#$}
\medskip

\par Let $p_z$ be as in Theorem 8.2, now
with $z\in {\rm neigh\,}(0,{\bf C})$ and
assume that $p_z$ is holomorphic in $z$. We
want to choose
$z_j\mapsto f_j, j=1,2$, so that with $z=z_1+iz_2$,
$\Lambda _z$ is critical for $I(\cdot
,p_z)$, where $\Lambda _w=\exp
i(w_1H_{f_1}+w_2H_{f_2})(\Lambda
_0)$, and where $f_j$ are real and smooth
on $\Lambda _0$ and the Hamilton flow is
defined to infinite order at $\Lambda
_0\times \{ w=0\}$. We will only consider
the infinitesimal solution of the problem,
that we obtain from (8.24), so that for
$z=w=0$:
$$-(H_p+a)H_{\overline{p}}f_1={1\over
2i}(H_p+a) \overline{{\partial p\over
\partial \Re z}}-{1\over
2i}(H_{\overline{p}}+\overline{a}){\partial
p\over \partial \Re z},$$
$$-(H_p+a)H_{\overline{p}}f_2={1\over
2i}(H_p+a)\overline{{\partial p\over
\partial \Im z}}-{1\over
2i}(H_{\overline{p}}+
\overline{a}){\partial p\over \partial
\Im z}.$$
Equivalently we have for $z=w=0$:
\eekv{8.29}
{-(H_p+a)H_{\overline{p}}f_w=-{1\over
2i}(H_{\overline{p}}+\overline{a})(
{\partial p\over \partial z}),}
{-(H_p+a)H_{\overline{p}}f_{\overline{w}}
={1\over
2i}(H_p+a)(\overline{
{\partial p\over \partial z}}),}
where $f_w={1\over 2}(f_1+{1\over i}f_2)$.
We then know that for the function
$I(\Lambda _w,p_z)$,
\ekv{8.30}
{
\nabla _z((\nabla _wI)(\Lambda
_z,p_z))=0,\hbox{ at }z=0, }
or equivalently that
$$((\nabla _z+\nabla _w)\nabla
_wI)(\Lambda _w,p_z)=0,\hbox{ at }z=w=0,$$
which can be expanded in terms of
holomorphic and antiholomorphic
derivatives as:
\eekv{8.31}
{
(\partial _z\partial _w+\partial
_w^2)I=0,\ (\partial _z\partial
_{\overline{w}} +\partial _w\partial
_{\overline{w}})I=0,}
{
(\partial _{\overline{z}}\partial
_w+\partial _{\overline{w}}\partial
_w)I=0,
\ (\partial _{\overline{z}}\partial
_{\overline{w}}+\partial
^2_{\overline{w}})I=0,\hbox{ for }z=w=0. }
>From these relations we see that
$\partial _z\partial
_{\overline{w}}I=\partial
_{\overline{z}}\partial _wI$ is real (at
$z=w=0$) and that:
\ekv{8.32}
{
\partial _{\overline{z}}\partial
_z(I(\Lambda _z,p_z))=(\partial
_{\overline{z}}+\partial _{\overline{w}})
(\partial _z+\partial _w)I(\Lambda
_w,p_z)=\partial _{\overline{z}}\partial
_zI(\Lambda _w,p_z)-\partial
_{\overline{w}}\partial _wI(\Lambda
_w,p_z). }
Since
$${\partial \over\partial
\overline{w}}{\partial \over \partial
w}={1\over 4}(({\partial \over \partial
w_1})^2+({\partial \over \partial
w_2})^2),$$
we have for $z=w=0$:
$$({\partial \over \partial
\overline{w}}{\partial \over \partial
w})I(\Lambda _w,p_z)={1\over 4}(({\partial
\over \partial w_1})^2I_1(w_1)+({\partial
\over \partial w_2})^2I_2(w_2)),$$
$$
I_j(w_j)=\cases{I(\Lambda
_{(w_1,0)},p_0),\, j=1\cr I(\Lambda
_{(0,w_2)},p_0),\, j=2.}$$
In the definition of $I_j(w_j)$ we can
let $w_j$ be complex and we notice
that to infinite order at $w_j=0$, we
have $I_j(w_j)\equiv I_j(\Re w_j)$. Hence
for $w_j=0$, $z=0$:
$$\eqalign{
&({\partial \over \partial
w_1})^2I_1(w_1)=(({\partial \over
\partial \Re w_1})^2+({\partial \over
\partial \Im w_1})^2)I_1(w_1)=\cr &
4{\partial \over \partial
\overline{w}_1}{\partial \over \partial
w_1}I_1(w_1)=4{\pi \over
2}\int_{p^{-1}(0)}H_pf_1\overline{H_pf_1}\lambda
_{p,0}(d\rho ), }$$
where we abused the notation since
${\partial \over \partial w_1}$,
${\partial \over \partial
\overline{w}_1}$ in the third member
denote holomorphic and antiholomorphic
derivatives, and where we used (8.18) in
the last step.  The same
calculation is valid for $j=2$, and we
get for $z=w=0$:
$$\partial _{\overline{w}}\partial
_wI={\pi \over
2}\sum_{j=1}^2\int_{p^{-1}(0)}H_pf_j
\overline{H_p f_j}\lambda _{p,0}(d\rho
)=-{\pi \over
2}\sum_{j=1}^2\int_{p^{-1}(0)}((
H_{\overline{p}}+\overline{a})H_p
f_j)\overline{f_j}\lambda _{p,0}(d\rho
).$$
Here we used again that the adjoint of
$H_p$ for $\lambda _{p,0}$ is
$-(H_{\overline{p}}+\overline{a})$. Now
recall that
$(H_{\overline{p}}+\overline{a}) H_p$ is a
real selfadjoint operator, to get
$$\sum_{j=1}^2\int_{p^{-1}(0) }
((H_{\overline{p}}+\overline{a})H_pf_j)
\overline{f}_j\lambda
_{p,0}(d\rho
)=4\int_{p^{-1}(0)}((H_{\overline{p}}+
\overline{a})H_pf_w)\overline{f_w}\lambda
_{p,0}(d\rho ),$$ and
\ekv{8.33}{\partial _{\overline{w}}\partial
_wI= -2\pi
\int_{p^{-1}(0)}((H_{\overline{p}}
+\overline{a} )H_pf_w)\overline{f_w}
\lambda _{p,0}(d\rho )=2\pi
\int_{p^{-1}(0)}H_pf_w\,\overline{H_pf_w}
\lambda _{p,0}(d\rho ).}

 Since
$(H_p+a)H_{\overline{p}}$ is a real
operator, the first equation in (8.29), can
also be written
\ekv{8.34}
{(H_{\overline{p}}+\overline{a})(
H_pf_w-{1\over 2i}\partial _zp)=0.}
Here the orthogonal of the image of $H_p$
is the kernel of
$-H_p^*=(H_{\overline{p}}+\overline{a})$,
so (8.34) says that
\ekv{8.35}
{
H_p(f_w)=\Pi ({1\over 2i}\partial _zp),
}
where $\Pi $ is the orthogonal projection
onto the image of $H_p$. Then (8.33)
becomes
\ekv{8.36}
{
\partial _{\overline{w }}\partial _wI=
{\pi \over 2}\int_{ p^{-1}(0)}\vert \Pi
\partial _zp\vert ^2\lambda _{p,0} (d\rho
). }

\par Using this and (8.15) in (8.32), we
finally get for $z=0$:
\ekv{8.37}
{
\partial _{\overline{z}}\partial
_z(I(\Lambda _z,p_z))= {\pi \over
2}\int_{p^{-1}(0)}\vert
\partial _zp\vert ^2\lambda _{p,0}(d\rho
)-{\pi \over 2}\int_{p^{-1}(0)}\vert \Pi
\partial _zp\vert ^2\lambda _{p,0}(d\rho
), }
which can also be written
\ekv{8.38}{\partial _{\overline{z}}\partial
_z(I(\Lambda _z,p_z))= {\pi \over
2}\int_{p^{-1}(0)}\vert
(I-\Pi )\partial _zp\vert ^2\lambda
_{p,0}(d\rho ).}
This formula gives some possible
indication about the distribution of
zeros of\break $\det p_z^w(x,hD_x)$.
\bigskip
\centerline{\bf 9. The case when $p_{\Lambda
_0}$ is of principal type.}
\medskip
\par Let $\Lambda _0$ be an IR-\mfld{} and $p$
holomorphic as in section 4. Write
$p_0=p_{\Lambda _0}$. Assume,
\ekv{9.1}
{p_0(\rho )=0\Rightarrow dp_0(\rho )\ne 0.}
Let $\rho _0\in p_0^{-1}(0)$. We can find
local symplectic coordinates $(x,\xi )$ on
${\rm neigh\,}(\rho _0,\Lambda _0)$, centered
at $\rho _0$, such that with $x=(x',x_n)$,
$\xi =(\xi ',\xi _n)$:
\ekv{9.2}
{
p_0=q(x,\xi )(\xi _n+ir(x,\xi ')),
}
with $q,r$ smooth, $q(0)\ne 0$, $r$
real valued with $r(0)=0$. We want to study
\ekv{9.3}
{
\langle {\rm "}d_\Lambda I(\Lambda _0,p){\rm
"},\delta \Lambda \rangle = \int \langle
H_{{\rm arg\,}p_0},df\rangle \mu (d\rho
)=-\int \langle d\,{\rm arg\,}p_0,H_f\rangle
\mu (d\rho ), }
with the integrals over $\Lambda _0$,
where we recall that $d\,{\rm arg\,}p_0$,
$H_{{\rm arg\,}p_0}$
extend from $\Lambda _0\setminus p_0^{-1}(0)$
to $L^1$ functions on $\Lambda _0$. Here $df$
is the generator of the infinitesimal
deformation $\delta \Lambda $ of $\Lambda
_0$. We restrict the attention to $f\in
C_0^\infty ({\rm neigh\,}(\rho _0,\Lambda
_0);{\bf R})$, so that (9.2) is applicable
with a non-vanishing $q$. Since
$$\int \langle d\,{\rm arg\,}q,H_f\rangle \mu
(d\rho )=0,$$
we get
\ekv{9.4}
{\langle {\rm "}d_\Lambda I(\Lambda _0,p){\rm
"},\delta \Lambda \rangle =\int \langle
H_{{\rm arg\,}p_1},df\rangle \mu (d\rho
)=-\int \langle d\,{\rm arg\,}p_1,H_f\rangle
\mu (d\rho ),}
where $p_1=\xi _n+ir(x,\xi ')$. We can view
(9.4) as the limit when $\epsilon \to 0$ of
\eeekv{9.5}
{
\langle d_\Lambda J_\epsilon (\Lambda
_0,p_1),\delta \Lambda \rangle ={\rm
Re\,}i\int_{\Lambda _0}{\overline{p_1}\over
p_1 \overline{p_1}+\epsilon ^2}H_fp_1\, \mu
(d\rho ) }
{
=-{\rm Re\,}i\int_{\Lambda _0}{\overline{p_1}
\over p_1\overline{p_1}+\epsilon
^2}H_{p_1}f\, \mu (d\rho )={\rm
Re\,}i\int_{\Lambda
_0}H_{p_1}({\overline{p_1}\over
p_1\overline{p_1}+\epsilon ^2})f\,\mu (d\rho )
}
{
={\rm Re\,}i\int_{\Lambda _0}{\partial
\over \partial
\overline{p_1}}({\overline{p_1}\over
p_1\overline{p_1}+\epsilon
^2})H_{p_1}(\overline{p_1})f\,\mu (d\rho ). }
Since $r$ is independent of $\xi _n$, we see
that
\ekv{9.6}
{H_{p_1}\overline{p_1}=H_{p_1}^{(n)}
\overline{p_1},
}
where $H_{p_1}^{(n)}$ denotes the $(\partial
_{x_n},\partial _{\xi _n})$ component of
$H_{p_1}$, i.e. the Hamilton field of $p_1$,
in the variables $(x_n,\xi _n)$ with $(x',\xi
')$ as parameters. We can insert this in the
last expression of (9.5) and run the
computation there backwards with
$H^{(n)}_{p_1}$, $H^{(n)}_f$ instead of
$H_{p_1}$, $H_f$, to get
\ekv{9.7}
{
\langle d_\Lambda J_\epsilon (\Lambda
_0,p),\delta \Lambda \rangle ={\rm
Re\,}i\int
{\overline{p_1}\over
p_1\overline{p_1}+\epsilon ^2}H_f^{(n)}p_1\,
\mu (d\rho ).  }
Here we can integrate first in $(x_n,\xi
_n)$ and then in $(x',\xi ')$. The earlier
discussion then applies to the
$(x_n,\xi _n)$ integral, making use of the
fact that $d^{(n)}{\rm arg\,}p_1:=
d_{(x_n,\xi _n)}{\rm arg\,}p_1$ is uniformly in
$L^1$, when $(x',\xi ')$ varies. (Cf the remark on
parameter dependence giving the last part of
Proposition 4.4, where $t$ can be replaced by
$(x',\xi ')$.) We get
\ekv{9.8}
{
\langle {\rm "}d_\Lambda I(\Lambda _0,p){\rm
"},\delta \Lambda \rangle =\int \langle
H_{{\rm arg\,}p_1}^{(n)},d^{(n)}f\rangle \mu
(d\rho )=-\int \langle d^{(n)}{\rm
arg\,}p_1,H_f^{(n)}\rangle \mu (d\rho ), }
where we recall that in view of the uniform
boundedness of $d^{(n)}{\rm arg\,}p_1$ in
$L^1$, we have
\ekv{9.9}
{
{\rm Re\,}i\iint {\overline{p_1}\over
p_1\overline{p_1}+\epsilon
^2}H_f^{(n)}p_1\, dx_nd\xi _n\to
\iint \langle H_{{\rm
arg\,}p_1}^{(n)},d^{(n)}f\rangle dx_nd\xi
_n,\ \epsilon \to 0, }
uniformly in $(x',\xi ')$, since the LHS in
(9.9) is equal to
\ekv{9.10}
{
-\iint {p_1\overline{p_1}\over
p_1\overline{p_1}+\epsilon ^2}\langle
d^{(n)}{\rm arg\,}p_1,H_f^{(n)}\rangle
dx_nd\xi _n. }
Now recall that the LHS of (9.9) is also equal
to
\eekv{9.11}
{
{\rm Re\,}i\iint {\partial \over \partial
\overline{p_1}}({\overline{p_1}\over
p_1\overline{p_1}+\epsilon ^2})\{
p_1,\overline{p_1}\} f\,dx_nd\xi _n=}
{{\rm
Re\,} \iint {\epsilon ^2\over
(p_1\overline{p_1}+\epsilon ^2)^2}i\{
p_1,\overline{p_1}\} f\, dx_nd\xi _n=2\iint
{\epsilon ^2\over (\epsilon ^2+\xi
_n^2+r(x,\xi ')^2)^2}(\partial _{x_n}r)f\,
dx_nd\xi _n.}
As we already know in arbitrary dimension
the limit only depends on the behaviour of
$f$ near the points $(x_n,\xi _n)=(x_n,0)$,
with $r(x,\xi ')=0$.

\par Assume for some fixed $(x',\xi ')$, that
\ekv{9.12}
{\hbox{The zeros of }x_n\mapsto r(x,\xi
')\hbox{ are all of finite order.}}
Assume that $x_n=0$ (say) is a zero of order
$m$ and assume that ${\rm supp}_{x_n,\xi
_n}f$ is contained in a small neighborhood of
$(0,0)$. We have
\ekv{9.13}
{
r=x_n^mu(x_n),\ u(0)\ne 0,
}
and consider 4 different cases.
\smallskip
\par\noindent \it 1) $m$ is odd and $u(0)>0$.
\rm Then $r$ considered as a function of $x_n$
is invertible near
$0$ and using that $\partial _{x_n}r\,
dx_n=dr$, the integral in (9.11) becomes
\ekv{9.14}
{
\iint {\epsilon ^2\over (\epsilon ^2+\xi
_n^2+y^2)^2}f(r^{-1}(y),\xi _n)\, dyd\xi
_n\to \pi f(0,0). }
\smallskip
\par\noindent \it 2) $m$ is odd and $u(0)<0$.
\rm  Now the map $r$ reverses the orientation
and we get the limit $-\pi f(0,0)$.
\smallskip
\par\noindent
\it 3) $m$ is even and $u(0)>0$. \rm The
restrictions $r_\pm$ to a neighborhood of $0$
on the positive and negative half axies
respectively are now invertible  and we can
cut the integral into two, which simplifies to
$$\eqalign{
\iint_{y\ge 0}{\epsilon ^2\over (\epsilon
^2+\xi _n^2+y^2)^2}(&f(r_+^{-1}(y),\xi
_n)-f(r_-^{-1}(y),\xi
_n))dyd\xi _n\cr &\to {\pi \over
2}(f(0,0)-f(0,0))=0 .}$$
\smallskip
\par\noindent
\it 4) $m$ is even and $u(0)<0$. \rm We get a
change of sign compared to the previous case,
and hence we still have the limit 0.
\smallskip

\par Summing up, we get
\medskip
\par\noindent \bf Proposition 9.1. \it
\par\noindent (A) Assume for a given $(x',\xi
')$, that (9.12) holds and define the
index $\iota (x_n)$ of a zero $x_n$ of
$r(x',\cdot ,\xi ')$ to be $+1$ if $r$
changes sign from $-$ to $+$, $-1$ if $r$
changes sign from $+$ to $-$, and to be 0 if
we have no change of sign which happens
precisely when $m$ is even. Then
\ekv{9.15}
{
\iint {\epsilon ^2\over (\epsilon ^2+\xi
_n^2+r(x,\xi ')^2)^2}(\partial
_{x_n}r)f(x_n,\xi _n)dx_nd\xi _n\to
\sum_{x_n;\,r(x,\xi ')=0}\pi \iota (x_n)
f(x_n,0). }
\smallskip
\par\noindent (B) Assume that (9.12) holds for
almost all $(x',\xi ')$. Then
\ekv{9.16}
{
\langle {\rm "}d_\Lambda I(\Lambda _0,p){\rm
"},\delta \Lambda \rangle =2\pi \iint
(\sum_{x_n; \,r(x,\xi ')=0}\iota (x_n)f(x,\xi
',0)) dx'd\xi '. }
\rm
\medskip
\par It is easy to construct an example of a
$p_0$ of the form (9.2) which does not satisfy
(9.12), such that $r$ has infinitely many zeros
in $x_n$ and for which $\langle {\rm
"}d_\Lambda I(\Lambda _0,p){\rm "},\delta
\Lambda
\rangle $ \it is a distribution of order 1, but
not a Radon measure, acting on $f$.\rm

\par The codimension 2 case,
i.e. the case when $dp_0$, $d\overline{p_0}$
are linearly independent at every point of
$p_0^{-1}(0)$, is equivalent to the case when
$r(x,\xi ')=0\Rightarrow dr(x,\xi ')\ne 0$,
and in this case (9.12) does hold for almost
all $(x',\xi ')$. In this case we can get
(9.16) more directly from (8.16). Indeed,
consider the situation locally and assume
that $p_0(\rho _0)=0$ and that
\ekv{9.17}
{d\,\Re p_0,\, d\,\Im p_0\hbox{ are linearly
independent at }\rho _0.}
Choose local symplectic coordinates $(x,\xi
)$ centered at $\rho _0$, so that
\ekv{9.18}
{p_0=\xi _n+ir(x,\xi ).}
Notice that $r$ is now allowed to depend on
$\xi _n$. Then by (9.17):
\ekv{9.19}
{d_{(x,\xi ')}r(0)\ne 0.}
Let us first assume that
\ekv{9.20}
{\{ \Re p_0,\Im p_0\} (\rho _0)\ne 0,}
or equivalently that $\partial _{x_n}r(0)\ne
0$. Then near $\rho _0$, the surface
$p_0^{-1}(0)$ is given by $\xi _n=0$,
$r(x,\xi ',0)=0$ and can be parametrized by
$x_n=x_n(x',\xi ')$. The Liouville measure
becomes
\ekv{9.21}
{
L(d(x',\xi '))={dx'd\xi '\over \vert \partial
_{x_n}r\vert }, }
and hence (8.16) becomes
$$\langle d_\Lambda I(\Lambda _0),\delta
\Lambda \rangle =2\pi \iint f(x',x_n(x',\xi
'),\xi ',0) \mu (d(x',\xi ')), $$
where
\ekv{9.22}
{
\mu =\{ \Re p_0,\Im p_0\} L(d(x',\xi
'))={\partial _{x_n}r\over \vert \partial
_{x_n}r\vert }dx'd\xi '={\rm sgn\,}(\partial
_{x_n}r)dx'd\xi '. }
recall from Lemma 8.1, that this density is
also given by
$${1\over (n-1)!}{{\sigma ^{n-1}}_\vert
}_{p_0^{-1}(0)},$$
with a suitable choice of orientation on
$p_0^{-1}(0)$. Since $\xi _n=0$ on
$p_0^{-1}(0)$, we see that even when (9.20)
is not fulfilled (but still under the
assumption (9.17)), the absolute value of the
density appearing in (8.16) is bounded by
\ekv{9.23}
{
{\vert dx'd\xi '_\vert}_{r^{-1}(0)\cap \xi
_n^{-1}(0)}\vert . }
We already know that this density vanishes
precisely when $\partial _{x_n}r= 0$, and
those points are precisely the ones where the
projection
$$r^{-1}(0)\cap \xi _n^{-1}(0)\ni (x,\xi
')\mapsto (x',\xi ')\in{\bf R}^{2(n-1)}$$
is not a local diffeomorphism. The Lebesgue
measure on $r^{-1}(0)\cap \xi _n^{-1}(0)$ is
locally equivalent to the Liouville measure
$L$, given by (9.21), and by (9.22), we have
\ekv{9.24}
{
\vert \mu \vert =\vert \partial _{x_n}r\vert
L. }
We conclude that
\ekv{9.25}
{
\int_{\{ (x,\xi ',0)\in r^{-1}(0)\cap{\rm
neigh\,}(0);\vert \partial _{x_n}r\vert \le
\epsilon \}}\vert \mu \vert (d(x,\xi '))\to 0,\
\epsilon \to 0, }
and this means that the points in
$p_0^{-1}(0)$, where $\partial _{x_n}r(x,\xi
',0)=0$ can be neglected in the (local) formula
for $\langle d_\Lambda I(\Lambda _0,p),\delta
\Lambda \rangle $. Further, by Sard's theorem
(here in an easy case), the set
\ekv{9.26}
{
\{ (x',\xi ')\in{\rm neigh\,}(0,{\bf
R}^{2(n-1)});\exists x_n\in{\rm
neigh\,}(0,{\bf R}),\, r(x,\xi ',0)=0,\,
\partial _{x_n}r(x,\xi ',0)=0\} }
is of Lebesgue measure zero.

\par Summing up, under the assumption
(9.17), and for $f\in C_0^\infty ({\rm
neigh\,}(\rho _0,\Lambda _0))$, we have,
using the coordiates in (9.18):
\ekv{9.27}
{
\langle d_\Lambda I(\Lambda _0,p),\delta
\Lambda \rangle =2\pi \iint dx'd\xi '
\sum_{x_n;\, r(x,\xi ',0)=0}\iota (x,\xi
')f(x,\xi ',0), }
where $\iota (x,\xi ')={\rm sgn\,}\partial
_{x_n}r(x,\xi ',0)$, if $\partial _{x_n}r(x,\xi
',0)\ne 0$ and $\iota (x,\xi ')=0$
otherwise. Here we use the fact that the set
(9.26) is of measure 0. Notice that if
$(x',\xi ')$ is not in that set, and we
enumerate the zeros $x_n=x_n^j(x',\xi ')$ of
$r(x,\xi ',0)$ in increasing order for $j$ in
a subinterval of ${\bf Z}$, then we may assume
that ${\rm sgn\,}\partial _{x_n}r=(-1)^j$,
and we get
$$\eqalign{
\vert \sum_j \iota (x',x_n^j(x',\xi '),\xi
')f(x',x_n^j,\xi ',0)\vert &=\vert \sum_k
(f(x',x_n^{2k},\xi ',0)-f(x',x_n^{2k-1},\xi
',0))\vert \cr &\le C\sup_t \vert \partial
_{x_n}f(x',t,\xi ',0)\vert , }$$
so we have a locally uniform bound on the RHS
even when approaching the set (9.26) where the
number of zeros may tend to infinity. (The
need to use
$\vert
\nabla f\vert $ looks a little strange since
we are working under assumptions that imply
that our differential is a Radon measure with
respect to $f$.)

\par We next want to globalize the formula
(9.27) and make the global assumptions
\ekv{9.28}
{
p_0(\rho )=0\Rightarrow d\,\Re p_0(\rho ),\,\,
d\,\Im p_0(\rho )\hbox{ are independent,}
}
\ekv{9.29}
{
\Re p_0(\rho )=0\Rightarrow d\,\Re p_0(\rho
)\ne 0. }
If we choose local symplectic coordinates
such that $\Re p_0=\xi _n$, then the Liouville
measure on $({\rm Re\,}p_0)^{-1}(0)$ is
$dx'd\xi 'dt$ where $t=x_n$ corresponds to
the time of the Hamilton flow of $\Re p_0$. A
contribution from a zero in (9.27): $\iota
(x',x_n^j,\xi ')f(x',x_n^j,\xi ',0)$ can
be expressed as
$$2\pi \iiint_{t\in I_{x',\xi '}}{1\over
T(x',\xi ')}\iota (x',x_n^j,\xi
')f(x',x_n^j,\xi ',0) dtdx'd\xi ',$$
where $I_{x',\xi '}$
is some interval of length $T(x',\xi ')>0$.

\par We shall use this without refering to
(9.18) first in the case when the length of
the time interval is fixed and then for a
special choice of variable length. First in
the case of a fixed length we have
\ekv{9.30}
{
\langle d_\Lambda I(\Lambda _0,p),\delta
\Lambda \rangle ={2\pi \over 2T}\int_{(\Re
p_0)^{-1}(0)}\sum_{t\in ]-T,T[;\, \Im
p_0(\Phi _t(\rho ))=0}\iota (\Phi _t(\rho
))f(\Phi _t(\rho ))\,\lambda _{\Re
p_0=0}(d\rho ), }
where $\Phi _t(\rho )=\exp tH_{\Re p_0}(\rho
)$, $\lambda _{\Re p_0=0}$ is the Liouville
measure on $(\Re p_0)^{-1}(0)$ and $\iota
(\rho )$ is equal to $\pm 1$ if $\pm\Im
p_0(\Phi _s(\rho ))$ has a simple zero at
$s=0$ with a change of sign from $-$ to $+$
and is 0 otherwise. Notice that for $\rho \in
(\Re p_0)^{-1}(0)$ away from some set of
measure 0, all the zeros of ${\bf R}\ni t
\mapsto \Im p_0(\Phi _t(\rho ))$ are simple.
\medskip
\par\noindent \bf Lemma 9.2. \it The set of
points $\rho \in (\Re p_0)^{-1}(0)$, such that
$\Im p_0(\Phi _t(\rho ))$ is $<0$ for some
$t\in{\bf R}$ and $\ge 0$ for all sufficiently
large positive $t$, is of measure 0.
\medskip
\par\noindent \bf Proof. \rm Let $\Omega $ be
the set in question, so that $\Omega $ is a
union of non-closed trajectories. Let
$$\Omega
_+=\{\rho \in\Omega ;\, \Im p_0(\Phi _t(\rho
))\ge 0,\hbox{ for }t\ge 0\} .$$
Then
$$\Phi _t(\Omega _+)\to\cases{ \emptyset,\ t\to
+\infty, \cr \Omega ,\ t\to -\infty  }$$
so by the dominated convergence theorem,
$$\lambda _{\Re p_0=0}(\Phi _t(\Omega _+))\to
\cases{0,\ t\to +\infty ,\cr \lambda _{\Re
p_0=0}(\Omega ), \ t\to -\infty .}$$
But $\Phi _t$ is measure preserving so
$\lambda _{\Re p_0=0}(\Phi _t(\Omega _+))$ is
independent of $t$ and we conclude that
$\lambda _{\Re p_0 =0}(\Omega )=0$.\hfill{$\#$}
\medskip
\par We have 4 variants of the lemma
 since we may replace $\Im p_0$ by $-\Im
p_0$ and $t$ by $-t$. It follows that for
$\rho $ outside some set of Liouville
measure 0, we are in one of the following
three cases:
\smallskip
\par\noindent 1) $\Im p_0(\Phi _t(\rho ))\ge
0$ for all $t\in{\bf R}$,
\smallskip
\par\noindent 2) $\Im p_0(\Phi _t(\rho ))\le
0$ for all $t\in{\bf R}$,
\smallskip
\par\noindent 3) $t\mapsto \Im p_0(\Phi
_t(\rho ))$ has infinitely many changes of
sign both when $t\to +\infty $ and when $t\to
-\infty $. Moreover each zero is simple.
\smallskip
\par Here the last sentence in 3) could be
added because of the observation
prior to Lemma 9.2.
\par Define a function $F$ on $(\Re
p_0)^{-1}(0)$ in the following way: Put
$F(\rho )=0$ in the cases 1) or 2) above or if
$t\mapsto \Im p_0(\Phi _t(\rho ))$ has at
least one zero which is not simple.
\par In case 3)
let $t(\rho )\le 0<s(\rho )$ be the zeros of
$\Im p_0(\Phi _t(\rho ))$ which are closest to
$0$ in the sense that $\Im p_0(\Phi _t(\rho
))\ne 0$ for $t(\rho )<t<s(\rho )$. If
${{\Im p_0(\Phi _\cdot (\rho
))}_\vert}_{]t(\rho ),s(\rho )[}>0$, we put
$F(\rho )=0$ and in the opposite case we put
\ekv{9.31}
{
F(\rho )={1\over s(\rho )-t(\rho )}(f(\Phi
_{s(\rho )}(\rho ))-f(\Phi _{t(\rho )}(\rho
)))={1\over s(\rho )-t(\rho )}\int_{t(\rho
)}^{s(\rho )}(H_{\Re p_0}f)(\Phi _t(\rho ))dt.
} Alternatively we could take
\ekv{9.32}
{
\widetilde{F}(\rho )=H_{\Re p_0}(f)(\rho )
}
in the last case and $\widetilde{F}(\rho
)=F(\rho )=0$ in the other cases.

\par With this choice of $F$ or with $F$
replaced by $\widetilde{F}$, we have
\ekv{9.33}
{
\langle d_\Lambda I(\Lambda _0,p),\delta
\Lambda \rangle =2\pi \int_{(\Re
p_0)^{-1}(0)}F(\rho )\,\lambda _{\Re
p_0=0}(d\rho ). }

\par Before continuing the main discussion,
we need an auxiliary result:
\medskip
\par\noindent \bf Lemma 9.3. \it Let $M$ be
a compact smooth \mfld{} and let $v$ be a
smooth non-vanishing \vf{} on $M$. If $f_0\in
C^\infty (M;{\bf R})$, then we can find $f\in
C^\infty (M;{\bf R})$ arbitrarily close to
$f_0$ in the $C^\infty $ topology, such that
\ekv{9.34}
{
v(f)(x)=0\Rightarrow d(v(f))(x)\ne 0,
}
for every $x\in M$.
\medskip
\par\noindent \bf Proof. \rm
We can cover $M$ by finitely many open sets
$\Omega _j$, $j=1,..,N$ such that for each
$j$ there is a real function $g_j\in C^\infty
(M) $ with $vg_j=1$ in $\overline{\Omega
}_j$. Consider $f_1=f_0-\epsilon_1 g_1$. In
$\overline{\Omega }_1$ we have
$v(f_1)=v(f_0)-\epsilon _1$, and we can
choose $\epsilon _1$ arbitrarily small, so
that $\epsilon _1$ is not a critical value of
$v(f_0)$ on $\overline{\Omega }_1$. Then
there exists $\delta _1>0$, such that
\ekv{9.35}
{\vert \nabla (vf_1)\vert +\vert vf_1\vert \ge
\delta _1}
in $\overline{\Omega }_1$. Define
$f_2=f_1-\epsilon _2g_2$, so that
$v(f_2)=v(f_1)-\epsilon _2$ in
$\overline{\Omega }_2$. We choose $\epsilon
_2$ very small and not equal to any critical
value of $v(f_1)$ on $\Omega _2$. Then
(9.35) holds in $\overline{\Omega }_1$ after
replacing
$(f_1,\delta _1)$ by $(f_2,\delta _1/2)$ and there is
some $\delta _2>0$ such that
\ekv{9.36}
{\vert \nabla (vf_2)\vert +\vert vf_2\vert \ge
\delta _2,}
on $\overline{\Omega _1\cup \Omega _2}$.
Continuing this procedure we get the Lemma
after $N$ steps.
 \hfill{$\#$}
\medskip

\par We next consider the case when $p_0$ is
of real principal type. Let $\Lambda _0$ be
an IR \mfld{ } and $p$ a holomorphic function
as in section 4. Let $p_0=p_{\Lambda _0}$. We
now assume
\ekv{9.37}
{
p_0\hbox{ is real-valued}.
}
\ekv{9.38}
{
dp_0\ne 0\hbox{ on }p_0^{-1}(0),
}
and that $p_0^{-1}(0)\ne \emptyset$.
Let ${\rm neigh\,}(0,{\bf R})\ni
t\mapsto\Lambda _t$ be a smooth deformation
of IR \mfld s as in section 4, with $\Lambda
_{t=0}=\Lambda _0$. Let $f_t\in C_b^\infty
(\Lambda _t;{\bf R})$ be a corresponding
smooth family of generators. Then (7.3),
(7.4) are satisfied, so Theorem 7.1 applies
and shows that
\ekv{9.39}
{
I(\Lambda _t,p)\ge I(\Lambda _0,p)-C_N\vert
t\vert ^N, }
for every $N\in{\bf N}$. $\Lambda _0$ is
therefore a critical \mfld{ } in a
generalized sense and we shall now see that
the derivative of $t\mapsto I(\Lambda _t,p)$
has a jump discontinuity at $t=0$ if $f_0$ is
appropriately chosen.

\par Apply Lemma 9.3 with $M=p_0^{-1}(0)$,
$v=H_{p_0}$ and conclude that there exists
$f=f_0\in C_b^\infty (\Lambda _0;{\bf R})$
such that:
\ekv{9.40}
{
{{d\, H_{p_0}f}_\vert }_{M}\ne 0\hbox{
whereever }H_{p_0}f=0. }
>From now on we will assume that $f_0$ has
this property.

\par Identifying $\Lambda _t$ and $\Lambda
_0$ by means of the symplectic map $\kappa
_{0,t}:\Lambda _0\to\Lambda _t$ of section
1, we can view $p_t$ as a function on
$\Lambda _0$. By Taylor expansion, we get
\ekv{9.41}
{
p_t(\rho )=p_0(\rho )+itH_fp_0(\rho )+{\cal
O}(t^2)\hbox{ in }C_b^\infty . }
Hence
\ekv{9.42}
{
\Re p_t=p_0+{\cal O}(t^2),\ {1\over t}\Im
p_t=-H_{p_0}f+{\cal O}(t). }
For $t\ne 0$, we put
\ekv{9.43}
{
\Sigma _t=\{ \rho \in\Lambda _0;\, p_t(\rho
)=0\}, }
and for $t=0$:
\ekv{9.44}
{
\Sigma _0=\{ \rho \in\Lambda _0;\, p_0(\rho
)=0,\, H_{p_0}f(\rho )=0\} . }
>From (9.40,42) it follows that $\Sigma _t$,
$\Sigma _0$ are smooth compact submanifolds
of $\Lambda _0$ of codimension 2, and that
${\rm dist\,}(\Sigma _0,\Sigma _t)={\cal
O}(t)$, in the natural sense. We notice that
$\Sigma _t=\widetilde{p}_t^{-1}(0)$, where
$\widetilde{p}_t:=\Re p_t+i{\Im p_t\over t}$,
for $t\ne 0$ and
$\widetilde{p}_0=p_0-iH_{p_0}f$. Then
$\widetilde{p}_t-\widetilde{p}_0={\cal O}(t)$
in $C_b^\infty $ and $d\Re \widetilde{p}_t$,
$d\Im \widetilde{p}_t$ are linearly
independent on $\Sigma _t$. For $t\ne 0$,
we have
\ekv{9.45}
{
{i\over 2}\{ p_t,\overline{p_t}\}\lambda
_{p_t=0}=({\rm
sgn\,}t)\{\widetilde{p}_t,\overline{
\widetilde{p}_t }\}\lambda _{\tilde{p}_t=0}}
on $\Sigma _t$. Further
$${i\over 2}\{ \widetilde{p}_t,\overline{
\widetilde{p}_t}\}=\{\Re \widetilde{p}_t,\Im
\widetilde{p}_t\} =\{\Re \widetilde{p}_0,\Im
\widetilde{p}_0\} +{\cal
O}(t)=-H_{p_0}^2f+{\cal O}(t),$$
and we get for $t\ne 0$ (cf. (8.16)):
\eekv{9.46}
{
\partial _tI(\Lambda _t,p)=({\rm sgn\,}t)2\pi
\int_{\Sigma _t}f_t{i\over 2}\{
\widetilde{p}_t,\overline{\widetilde{p}_t}\}
\lambda _{\tilde{p}_t=0}(d\rho )}
{\hskip 2cm = ({\rm
sgn\,}t) 2\pi \int_{\Sigma _0}f_0{i\over 2}
\{\widetilde{p}_0,\overline{\widetilde{p}_0}
\} \lambda _{\tilde{p}_0=0}(d\rho )+{\cal
O}(t). }
>From this and (9.39) we see that
\ekv{9.47}
{
2\pi \int_{\Sigma _0}f_0{i\over 2}\{
\widetilde{p}_0,\overline{\widetilde{p}_0}\}
\lambda _{\tilde{p}_0=0}(d\rho )\ge 0. }

\par The following proposition shows that most
of the time we have strict inequality and
hence that $\partial _tI(\Lambda _t,p)$ has a
jump discontinuity at $t=0$. The proof also
gives a more direct explanation of (9.47).
\medskip
\par\noindent \bf
Proposition 9.4. \it Assume that there is a
point $\rho _0\in\Sigma _0$ where ${i\over
2}\{ \widetilde{p}_0,\overline{\widetilde{p}
_0}\}=-H_{p_0}^2f\ne 0$. Then we have strict
inequality in (9.47).\medskip
\par\noindent
\bf Proof. \rm The expression (9.47) if
formally equal to $\langle d_\Lambda
I(\Lambda _0,\widetilde{p}_0),\delta \Lambda
\rangle $, with $\delta \Lambda $ generated
by $f_0$. (This is only formal, since
$\widetilde{p}_0$ does not in general have a
holomorphic extension.) The discussion
starting at (9.17) applies and we have (9.33)
with $F$ defined there with $f$ equal to
$f_0$ and with $p_0$ replaced by
$\widetilde{p}_0$. The points  $\rho \in
p_0^{-1}(0)$ of type 3) are the ones for
which $t\mapsto -(H_{p_0}f)(\Phi _t(\rho ))$
has only simple zeros, infinitely many near
both $t=+\infty $ and $t=-\infty $. The
assumption in the proposition implies (and is
in fact equivalent to) the fact that the
points of type 3) form a set of measure $>0$.
Let
$\rho $ be point of type 3) for which $F(\rho
)$ maybe $\ne 0$ i.e. for which the closest
zeros $t(\rho )\le 0<s(\rho )$ of $\Im
\widetilde{p}_0(\Phi _t(\rho
))=-(H_{p_0}f)(\Phi _t(\rho ))$ are such that
$-(H_{p_0}^2f)(\Phi _t(\rho ))$ is $<0$ for
$t=t(\rho )$ and $>0$ for $t=s(\rho )$. Then
$t=t(\rho ), s(\rho )$ are subsequent local
extrema of the function $t\mapsto f(\Phi
_t(\rho ))$ which is strictly increasing
between the two points. Consequently $f(\Phi
_{s(\rho )}(\rho ))-f(\Phi _{t(\rho )}(\rho
))>0$. We conclude that $F(\rho )$ in (9.31)
is $\ge 0$ with strict inequality on a set of
measure $>0$. As already noticed,
\ekv{9.48}
{
2\pi \int_{\Sigma _0}f_0{i\over 2}\{
\widetilde{p}_0,\overline{\widetilde{p}_0}\}
\lambda _{\tilde{p}_0=0}(d\rho
)=2\pi \int_{p_0^{-1}(0)}F(\rho )\lambda
_{p_0=0}(d\rho ). }
and the proposition follows.\hfill{$\#$}
\medskip

\par At least in the case when $f=f_0$
extends to a bounded holomorphic function in
a tube and $\Lambda _t=\exp iH_f(\Lambda
_0)$, there is a more general way of
detecting a jump discontinuity at $t=0$ of
$\partial _tI(\Lambda _t,p)$, when
$p_0^{-1}(0)$ contains a real hypersurface,
even without assuming that $p_0$ is of
principal type. This can be done by examining
the second derivative and the Levi form with
respect to $t$, and we hope to develop this
point of view in some future work. Another
possible approach would be to look for jump
disconituities in $d\,{\rm arg\,}p_t$ as a
function of $t$.
\bigskip
\centerline{\bf 10. Examples}
\medskip
\par We consider two simple examples with
$\Lambda $ fixed and with $p$ depending on a
complex spectral parameter. Let $\Lambda _0$
and $p$
be as in section 4, and assume that $p=p(z)$
depends holomorphically on $z\in \Omega $
where $\Omega \subset\subset{\bf C}$ is open.
Recall that
$$
I(\Lambda _0,p(z))={1\over 2}\int_{\Lambda
_0}\log (p(z)\overline{p(z)})\mu (d\rho ). $$
Using (8.3), we get
\ekv{10.1}{\partial _z\partial
_{\overline{z}}({1\over 2}\log (p(\rho
,z)\overline{p(\rho ,z)}))=\pi
\delta (p(\rho ,z))\partial _zp(\rho
,z))\overline{\partial _zp(\rho ,z)},}
near simple zeros of $z\mapsto p(\rho ,z)$.
We want to discuss the special case when
$p(\rho ,z)=p(\rho )-z$, with $p(\rho )$ as
in section 4, but since this symbol tends to
$1-z$ rather than $1$ when $\rho \to \infty $,
we need to consider a modified symbol,
corresponding to the estimate of a relative
determinant. Let
$\widetilde{p}(\rho )\in S(\Lambda ,1)$
with
\ekv{10.2}{p(\rho )-\widetilde{p}(\rho )={\cal
O}(\langle \rho \rangle ^m),\ m<-2n.}
Assume that $\Omega \subset\subset {\bf
C}\setminus \{ 1\}$ and that
\ekv{10.3}
{
\widetilde{p}(\rho )\not\in\overline{\Omega
},\ \rho \in\Lambda _0. }
Consider
\ekv{10.4}
{
I(z)=I(\Lambda_0 ,{p-z\over
\widetilde{p}-z})={1\over 2}\int_{\Lambda
_0}(\log
((p-z)\overline{(p-z)}-\log((\widetilde{p}-z)
\overline{(\widetilde{p}-z)}))\mu (d\rho ), }
for $z\in \Omega $.
Then,
\ekv{10.5}
{
\partial _z\partial _{\overline{z}}I(z)=\pi
\int_{\Lambda _0}\delta (p(\rho )-z)\mu
(d\rho ).  }
This means that if $\phi \in C_0^\infty
(\Omega )$, then
\ekv{10.6}
{
\int_{{\bf C}}\phi (z)(\partial _z\partial
_{\overline{z}}I(z))L(dz)=\pi \int_{\Lambda_0}
\phi (p(\rho ))\mu (d\rho ). }
In other words $\partial _z\partial
_{\overline{z}}I(z)$ can be described in
$\Omega $
as the direct image under $p$ of the
symplectic volume on $\Lambda _0$:
\ekv{10.7}
{
\partial_z\partial _{\overline{z}}I(z)=\pi
p_{*}(\mu ). }

\par If we assume that $p_0=p_{\Lambda _0}$
is real-valued and that $\Omega \cap {\bf
R}\subset ]-\infty ,0[$, then the measure
$\partial _z\partial _{\overline{z}}I(z)$ on
$\Omega $ is supported in $\Omega \cap
{\bf R}$, and is given there by the Stieltjes
measure $\pi dV(E)$, with
\ekv{10.8}
{
V(E)=\int_{\rho ;\, p(\rho )\le E}\mu (d\rho
). }

\par For our second example, we drop the
assumption that $p_0=p_{\Lambda _0}$ be real
and assume instead that $\Omega $  is a small
neighborhood of 0 in ${\bf C}$, that
$p_0^{-1}(0)$ consists of precisely one point
$\rho _0$ and that $\vert p_0(\rho )\vert
\sim{\rm dist\,}(\rho ,\rho _0)^2$, for $\rho
\in {\rm neigh\,}(\rho _0,\Lambda _0)$. Let
$\nu =p_{0\,*}(\mu )$, so that $\partial
_z\partial _{\overline{z}}I(z)=\pi \nu $. We
observe that $\nu (D(0,r))\sim r^n$.

\par Let $\widehat{p}$ be a second
holomorphic function with the same
properties and the same point $\rho _0$, and
assume that
\ekv{10.9}
{
\widehat{p}_0(\rho )-p_0(\rho )={\cal O}({\rm
dist\,}(\rho ,\rho _0)^{2N_0}), }
for some $N_0>1$. (Actually the discussion is
valid with a comparison function
$\widehat{p}_0\in C_b^\infty (\Lambda _0)$
which does not necessarily have a holomorphic
extension $\widehat{p}$.) Let $V=V_r\subset
{\bf C}\setminus \{ 0\}$ with $V\subset
D(0,r)$, $0<r\le 1$. If $z=p_0(\rho )\in V$
for some $\rho \in \Lambda _0$, then ${\rm
dist\,}(\rho ,\rho _0)={\cal O}(r^{1/2})$
and it follows that $\widehat{p}_0(\rho )\in
V+D(0,Cr^{N_0})$. Consequently,
\ekv{10.10}
{
\nu (V)\le \widehat{\nu }(V+D(0,Cr^{N_0})),
} where $\widehat{\nu
}=\widehat{p}_{0\,*}(\mu )$. The preceding
argument is symmetric in $p$, $\widehat{p}$,
so we also have
\ekv{10.11}
{
\widehat{\nu }(V)\le \nu (V+D(0,Cr^{N_0})).
}
The last two estimates express that $\partial
_z \partial _{\overline{z}}I(z)$ and
$\partial _z\partial
_{\overline{z}}\widehat{I}(z)$ (where
$\widehat{I}(z)$ is defined as in (10.4),
with $\widehat{p}$ replacing $p$) are close
to each other near $z=0$.

\par The last example is motivated by the
study of resonances for the semi-classical
Schr{\"o}dinger operator that are generated by
a non-degenerate critical point of the
potential, and we refer to [HeSj], [Sj3],
[BrCoDu] for more details. After a suitable
complex scaling or application of the theory
of [HeSj], we can reduce the study of
resonances  of  a Schr{\"o}dinger operator
$-h^2\Delta +V(x)$ in ${\bf R}^n$ in a
neighborhood of 0 (say), to that of the zeros
of a relative determinant $\det
(P-z)(\widetilde{P}-z)^{-1}$, where $P$ is a
realization of $-h^2\Delta +V$ in a suitable
space that is associated to an IR-\mfld{
} $\Lambda _0$ and with principal symbol
$p$, and $\widetilde{P}$ is a \schr{ }
with a potential $\widetilde{V}$ such that
$\widetilde{V}-V$ is sufficiently short
range and such that $\widetilde{P}$ has no
resonances in some neighborhood of 0.
(Actually, the property $p(\rho )\to 1$,
$\rho \to \infty $, must be replaced by a
more general ellipticity property, but that
does not affect the validity of the
discussion above.) Then as we saw in section
3, though we now have to appeal to the
[HeSj] theory which is similar but slightly
heavier than the one we explained in that
section, we obtain
\ekv{10.12}
{
\log \vert \det
((P-z)(\widetilde{P}-z)^{-1})\vert \le (2\pi
h)^{-n}(I(z)+o(1)),\ h\to 0. }
The symbol $\widehat{p}$ may be obtained by
using a quadratic approximation of $V$ at
the critical point and then we can take
$N_0=3/2$, or it can be obtained by using a
(more elaborate) Birkhoff normal form in
which case we may have larger values of
$N_0$ sometimes even $N_0=+\infty $.

\par As a special case, we may take $n=2$,
and assume that for some local smooth
symplectic coordinates $(x,\xi )$, centered
at $\rho _0$, we have
$$
\widehat{p}_0(\rho )=f({1\over 2}(\xi
_1^2+x_1^2),{1\over 2}(\xi _2^2+x_2^2)),
$$ near $\rho _0$ with $f\in C^\infty ({\rm
neigh\,}(0,{\bf C}))$ with
$$f(\iota_1,\iota_2)=\mu _1\iota_1-i\mu
_2\iota_2 +{\cal O}(\iota^2),$$
with $\mu _1,\mu _2>0$. Then near $0$:
$\widehat{\nu }=2\pi f_*L$, where $L$ is the
Lebesgue measure on the 4th quadrant in ${\bf
R}^2$.

\bigskip
\centerline{\bf References}
\medskip
\par\noindent
[AgCo] J. Aguilar, J.M. Combes, \it
A class of analytic perturbations for one-body\break Schr{\"o}dinger
Hamiltonians, \rm Comm. Math. Phys., 22(1971), 269--279.
\smallskip

\par\noindent
[BaCo] E. Balslev, J.M. Combes, \it Spectral properties of
  many-body  Schr{\"o}dinger operators with dilatation-analytic
  interactions, \rm
Comm. Math. Phys., 22(1971), 280--294.
\smallskip

\par\noindent
[BoSj] L. Boutet de Monvel, J. Sj{\"o}strand, \it Sur la singularit{\'e}
des noyaux de Bergman et de Szeg{\"o}, \rm
Journ{\'e}es {\'e}quations aux D{\'e}riv{\'e}es Partielles de Rennes (1975),
pp. 123--164. Ast{\'e}risque, No. 34-35,
Soc. Math. France, Paris, 1976.
\smallskip

\par\noindent
[BrCoDu] P. Briet, J.M. Combes, P. Duclos, \it On the location of
resonances  for\break Schr{\"o}dinger operators in the
semiclassical limit. I. Resonance free domains, \rm
 J. Math. Anal. Appl.  126(1)(1987), 90--99.
\smallskip

\par\noindent
[Bu] N. Burq, \it D{\'e}croissance de l'{\'e}nergie locale de l'{\'e}quation des ondes
pour le probl{\`e}me ext{\'e}rieur et absence de r{\'e}sonances au voisinage du
r{\'e}el, \rm Acta Math., 180(1998), 1--29.
\smallskip

\par\noindent
[Da] E.B. Davies, \it Semi-classical states for non-self-adjoint
Schr{\"o}dinger operators,\break\rm Comm. Math. Phys., 200(1999),
35--41.
\smallskip

\par\noindent
[DiSj] M. Dimassi, J. Sj{\"o}strand, \it Spectral asymptotics in the
semi-classical limit, \rm London Math. Soc. Lecture Notes Series 269,
Cambridge University Press 1999.

\par\noindent
[GoKr] I.C. Gohberg, M.G. Kre\u\i n, \it Introduction to the theory of
linear nonselfadjoint operators, \rm Translations of Mathematical
Monographs,  Vol. 18 American
Mathematical Society, Providence, R.I. 1969.
\smallskip

\par\noindent
[GuZw] L. Guillop{\'e}, M. Zworski, \it
Scattering asymptotics for Riemann surfaces, \rm Ann. of Math.
145 (1997), no. 3, 597--660.
\smallskip

\par\noindent
[HeSj] B. Helffer, J. Sj{\"o}strand,
\it  R{\'e}sonances en limite
  semi-classique, \rm  M{\'e}m. Soc. Math. France (N.S.) No. 24-25,
(1986).
\smallskip

\par\noindent [H{\"o}1] L. H{\"o}rmander, \it On the singularities of
partial  differential equations, \rm 1970 Proc. Internat. Conf. on
Functional Analysis and Related Topics (Tokyo, 1969) pp. 31--40 Univ. of
Tokyo Press, Tokyo.
\smallskip

\par\noindent [H{\"o}2] L. H{\"o}rmander, \it On the existence and the
regularity  of solutions of linear pseudo-differential
equations, \rm Enseignement Math. (2) 17 (1971), 99--163.
\smallskip

\par\noindent
[H{\"o}We] L. H{\"o}rmander, J. Wermer, \it Uniform approximation on
compact sets in ${\bf C}^n$, \rm Math. Scand., 23(1968), 5--21.
\smallskip

\par\noindent
[LeRo] G. Lebeau, L. Robbiano, \it Contr{\^o}le exact de l'{\'e}quation de la
chaleur, \rm Comm. PDE, 20(1995), 335--356.
\smallskip

\par\noindent
[Ma] B. Malgrange, \it  Ideals of differentiable functions, \rm Tata
Institute  of Fundamental Research Studies in
Mathematics, No. 3 Tata Institute of Fundamental Research, Bombay;
Oxford University Press, London 1967.
\smallskip

\par\noindent
[MaMa] A.S. Markus, V.I. Matsaev, \it
Comparison theorems for spectra of linear
operators, and spectral asymptotics, \rm
Trans. Moscow. Math. Soc. 1984(1),
139--187, and Trudy Moskov. Matem. O.
45(1984).
\smallskip

\par\noindent
[Me] R. Melrose, \it Polynomial bound on the distribution of
 poles in scattering by an obstacle, \rm Proc. Journ{\'e}es e.d.p. St
 Jean de Monts, 1984, Soc. Math. de France.
\smallskip

\par\noindent
[MeSj] A. Melin, J. Sj{\"o}strand, \it
Fourier integral operators with complex
valued phase functions, \rm Springer
Lecture Notes in Math., no 459.
\smallskip

\par\noindent
[Sj1] J. Sj{\"o}strand, \it Singularit{\'e}s analytiques microlocales, \rm
Ast{\'e}risque, 95(1982).
\smallskip

\par\noindent
[Sj2] J. Sj{\"o}strand, \it Function spaces associated to global
I-Lagrangian  manifolds, \rm pages 369-423 in Structure of solutions of
differential equations, Katata/Kyoto, 1995, World Scientific 1996.
\smallskip

\par\noindent
[Sj3] J. Sj{\"o}strand, \it Semiclassical resonances generated by a
non-degenerate critical point, \rm Springer LNM, 1256, 402-429.
\smallskip

\par\noindent
[Sj4] J. Sj{\"o}strand, \it A traceformula and review of some estimates
for
resonances, \rm p.377-437 in Microlocal Analysis and Spectral Theory,
NATO ASI Series C, vol.490, Kluwer 1997.
\smallskip

\par\noindent
[Sj5] J. Sj{\"o}strand, \it Geometric
bounds on the density of resonances for semiclassical  problems, \rm
Duke
Math. J., 60(1)(1990), 1-57.
\smallskip

\par\noindent
[SjZw] J. Sj{\"o}strand, M. Zworski, \it Complex
scaling and the distribution of scattering poles, \rm Journal of the
AMS,
4(4)(1991), 729-769.
\smallskip

\par\noindent
[Vo] G. Vodev, \it Sharp bounds on the number of scattering
poles  for perturbations of the Laplacian, \rm
Comm. Math. Phys. 146(1)(1992), 205--2
\smallskip

\par\noindent
[Ze] M. Zerzeri, \it Majoration du nombre de r{\'e}sonances
 pr{\`e}s de l'axe r{\'e}el pour une perturbation, {\`a} support compact,
 abstraite, du laplacien, \rm preprint, Universit{\'e} de Paris Nord,
 1999, Comm. PDE, to appear.
\smallskip

\par\noindent
[Zw1] M. Zworski, \it Poisson formulae for resonances, \rm
S{\'e}m.  {\'e}quations aux D{\'e}riv{\'e}es
Partielles, 1996--1997, Exp. No. XIII,  Ecole Polytech., Palaiseau,
1997.
\smallskip

\par\noindent
[Zw2] M. Zworski, \it A remark on a paper by E.B. Davies, \rm Proc.
AMS, to appear.
\smallskip

\par\noindent
[Zw3] M. Zworski, \it Sharp polynomial bounds on the number of
scattering poles, \rm Duke Math. J. 59(2)(1989), 311--323
\smallskip
\end